\documentclass[leqno,12pt]{amsart}

\usepackage[T1]{fontenc}
\usepackage[utf8]{inputenc}
\usepackage[english]{babel}
\usepackage[top=1.5in,bottom=1.3in,left=1.1in,right=1.1in,marginparwidth=1.5cm]{geometry}
\usepackage{times}
\usepackage{amssymb}
\usepackage{microtype}
\usepackage{xurl}
\usepackage{tikz-cd}
\usepackage{booktabs}
\usepackage{needspace}
\usepackage{aliascnt}
\usepackage[colorlinks=true,linkcolor=blue,citecolor=blue,urlcolor=blue,pagebackref=true]{hyperref}
\usepackage[capitalise,nameinlink]{cleveref}

\newcommand{\C}{\mathbb C}
\newcommand{\Q}{\mathbb Q}
\newcommand{\R}{\mathbb R}
\newcommand{\Z}{\mathbb Z}
\newcommand{\cD}{\mathcal D}
\newcommand{\cL}{\mathcal L}
\newcommand{\cN}{\mathcal N}
\newcommand{\cP}{\mathcal P}
\newcommand{\BG}{\operatorname{BG}}
\newcommand{\Kum}{\operatorname{Kum}}
\newcommand{\DR}{\operatorname{DR}}
\newcommand{\dbar}{\bar\partial}
\newcommand{\rank}{\operatorname{rank}}
\newcommand{\tr}{\operatorname{tr}}
\newcommand{\Hilb}{\operatorname{Hilb}}

\newcommand{\Sym}{\operatorname{Sym}}
\newcommand{\IC}{\operatorname{IC}}

\numberwithin{equation}{section}
\theoremstyle{plain}
\newtheorem{Thm}{Theorem}[section]
\newtheorem{ThmIntr}{Theorem}

\newaliascnt{Prop}{Thm}
\newtheorem{Prop}[Prop]{Proposition}
\aliascntresetthe{Prop}
\newaliascnt{Lem}{Thm}
\newtheorem{Lem}[Lem]{Lemma}
\aliascntresetthe{Lem}
\newaliascnt{Cor}{Thm}
\newtheorem{Cor}[Cor]{Corollary}
\aliascntresetthe{Cor}
\newaliascnt{Conj}{Thm}
\newtheorem{Conj}[Conj]{Conjecture}
\aliascntresetthe{Conj}
\theoremstyle{definition}
\newaliascnt{Def}{Thm}
\newtheorem{Def}[Def]{Definition}
\aliascntresetthe{Def}
\newaliascnt{Ex}{Thm}
\newtheorem{Ex}[Ex]{Example}
\aliascntresetthe{Ex}
\newaliascnt{Alg}{Thm}
\newtheorem{Alg}[Alg]{Algorithm}
\aliascntresetthe{Alg}
\theoremstyle{remark}
\newaliascnt{Rem}{Thm}
\newtheorem{Rem}[Rem]{Remark}
\aliascntresetthe{Rem}

\crefname{Thm}{Theorem}{Theorems}
\crefname{ThmIntr}{Theorem}{Theorems}
\crefname{Prop}{Proposition}{Propositions}
\crefname{Lem}{Lemma}{Lemmas}
\crefname{Cor}{Corollary}{Corollaries}
\crefname{Conj}{Conjecture}{Conjectures}
\crefname{Def}{Definition}{Definitions}
\crefname{Ex}{Example}{Examples}
\crefname{Alg}{Algorithm}{Algorithms}
\crefname{Rem}{Remark}{Remarks}
\crefname{equation}{}{}
\creflabelformat{equation}{#2(#1)#3}

\author[Nikon Kurnosov]{Nikon Kurnosov}
\address{London Institute for Mathematical Sciences, 21 Albemarle Street, London W1S 4DQ, United Kingdom}
\email{nk@lims.ac.uk}
\author[\'Angel David R\'ios Ortiz]{\'Angel David R\'ios Ortiz}
\address{Instituto de Matem\'aticas, Universidad Nacional Aut\'onoma de M\'exico, Unidad Cuernavaca, Av. Universidad s/n, Col. Lomas de Chamilpa, 62210 Cuernavaca, Morelos, M\'exico}
\email{david.rios@im.unam.mx}
\title[Hodge theory of Bogomolov--Guan manifolds]
{Hodge theory and topology of Bogomolov--Guan manifolds}

\begin{document}
\raggedbottom
\begin{abstract}
Bogomolov--Guan manifolds are simply connected, irreducible holomorphic symplectic manifolds that are not K\"ahler, yet they echo hyperk\"ahler geometry in striking ways. We show that their Fr\"olicher spectral sequence degenerates at $E_1$ and that their Hodge numbers are stable under small deformations. The general fibers of their natural Lagrangian fibration are rationally null-homologous, and the Beauville--Bogomolov--Fujiki form is degenerate. Despite this, in every dimension, the Euler characteristics agree with those of generalized Kummer manifolds. In dimensions four and six, we compute the full Hodge diamonds and all Chern and Pontryagin numbers. Surprisingly, these characteristic numbers are the same as the generalized Kummer ones, even though the Betti numbers differ and the Hodge diamonds are asymmetric.
\end{abstract}
\maketitle
\tableofcontents
\section{Introduction}

Irreducible holomorphic symplectic manifolds occupy a central position in compact K\"ahler geometry, where they are also known as hyper-Kähler manifolds. In the simply connected case, they form, together with Calabi--Yau manifolds, the fundamental building blocks of compact K\"ahler manifolds with vanishing first Chern class, by the Beauville--Bogomolov decomposition theorem \cite{Beauville1983}. Their geometry is governed by a remarkable interaction between Hodge theory, the topology of the underlying manifold, and the quadratic structure on second cohomology. The Fr\"olicher spectral sequence degenerates at the first page, the Hodge numbers satisfy Hodge symmetry, and $H^2$ carries the nondegenerate Beauville--Bogomolov--Fujiki form, whose Fujiki relation controls all top-degree intersections of degree-two classes.

Lagrangian fibrations provide one of the clearest manifestations of this structure. If $X$ is a hyper-K\"ahler manifold of dimension $2n$ admitting a Lagrangian fibration $\pi:X\to\mathbb P^n$, then the pullback $h=\pi^*c_1(\mathcal O_{\mathbb P^n}(1))$ is isotropic for the Beauville--Bogomolov--Fujiki form, while
\[
h^n\neq0
\]
is the cohomology class of a general Lagrangian fiber. Thus in the K\"ahler setting the quadratic form, the intersection theory of $X$, and the topology of its Lagrangian fibers are tightly linked.

The purpose of this paper is to understand what becomes of this picture once the K\"ahler hypothesis is removed. 

Outside the Kähler category the situation is strikingly less understood. There is no analogue of the Beauville--Bogomolov decomposition theorem \cite[Remark~2.9]{AnthesCattaneoRollenskeTomassini2018}, and the \textit{Bogomolov--Guan manifolds} \cite{GuanII1995,GuanIII1995,Bogomolov1996} remain the only known series of compact simply connected non-Kähler holomorphic symplectic manifolds. Beyond these examples there is at present no classification theory comparable to the Kähler case. In particular, it is not known whether the Bogomolov--Guan construction exhausts the possible deformation types, nor is there a general boundedness result for compact simply connected non-Kähler holomorphic symplectic manifolds. Already in complex dimension four, there is no known reason to expect only finitely many deformation types. 

Following \cite[Definition~1.1]{RiosOrtiz2026BGNonformal}, we say that a manifold is of $\BG_n$-type if it is deformation equivalent to a standard representative $Q_{n,t}$ of complex dimension $2n$, with $t\in\mathbb Z\setminus\{0\}$. For $n=1$ the construction gives a K3 surface, so throughout the paper we assume $n\ge2$.

Bogomolov--Guan manifolds are particularly striking because the failure of K\"ahlerness does not simply destroy the familiar hyper-K\"ahler package. On the contrary, substantial parts of it survive. Their Fr\"olicher spectral sequence still degenerates at $E_1$, they admit a Hodge decomposition, and their second cohomology carries a Beauville--Bogomolov--Fujiki-type quadratic form. At the same time, their topology is genuinely non-K\"ahler: as shown in \cite[Theorems~1.2 and~1.3]{RiosOrtiz2026BGNonformal}, every $\BG_n$-manifold is nonformal, and its odd cohomology already violates the constraints imposed by the $\partial\bar\partial$-lemma.

One of the aims of the present paper is to show that this tension is already visible in the interaction between the Beauville--Bogomolov--Fujiki form and the natural Lagrangian fibration. If $h$ denotes the pullback of the hyperplane class from the base, then
\[
h\neq0,
\qquad
h^2=0.
\]
This is impossible for the corresponding class on a hyper-K\"ahler Lagrangian fibration: there one has $h^n\neq0$, represented by the class of a general fiber. For Bogomolov--Guan manifolds, by contrast, the fiber class is homologically trivial.

We show that this phenomenon is not accidental. The class $h$ generates the radical of the Beauville--Bogomolov--Fujiki form. Thus the failure of nondegeneracy of the quadratic form is detected precisely by the class defining the Lagrangian fibration. Moreover, the same class enters the Massey-product construction responsible for the nonformality of Bogomolov--Guan manifolds. In this way, three apparently different manifestations of non-K\"ahler geometry---degeneracy of the Beauville--Bogomolov--Fujiki form, null-homology of the Lagrangian fibers, and nonformality of the rational homotopy type---are linked by a single cohomology class.

This gives a rather sharp picture of how Bogomolov--Guan geometry departs from the hyper-K\"ahler case. The Hodge decomposition and much of the degree-two formalism survive, but the quadratic form acquires a radical, the fibers of the natural Lagrangian fibration disappear in homology, and higher rational homotopy detects a nontrivial obstruction invisible to ordinary cohomology alone.

Alongside these results in arbitrary dimension, we determine the low-dimensional Hodge-theoretic and numerical invariants explicitly. In complex dimensions four and six we compute the complete Hodge diamonds, determine the relevant characteristic numbers, and compare them with those of the corresponding generalized Kummer manifolds. These computations make precise both how closely the Bogomolov--Guan examples continue to resemble hyper-K\"ahler manifolds and where that analogy necessarily breaks down. We now turn to the precise statements of our main results.

\subsection{Hodge theory}

Every $\BG_n$-type manifold is nonformal and has $b_3=3$ \cite[Theorems~1.2 and~1.3]{RiosOrtiz2026BGNonformal} (nonformality was first proved for $n=2$ and conjectured in general by Guan \cite{Guan2025Nonformal}). In particular, the underlying smooth manifold admits no complex structure satisfying the $\partial\bar\partial$-lemma. Nevertheless, the standard complex structures satisfy the following.

\begin{ThmIntr}\label{thm:main-degeneration}
For every $n\ge2$ and $t\ne0$, the Frölicher spectral sequence of $Q_{n,t}$ degenerates at $E_1$. The same holds for every sufficiently small deformation of $Q_{n,t}$, and all its Hodge numbers are locally constant there.
\end{ThmIntr}

The proof uses the projective semismall resolution in Guan's construction. Its strata are indexed by partitions of $n+1$, and the decomposition theorem expresses both de Rham and Dolbeault cohomology in terms of the corresponding sectors. These sectors are finite quotients of nilmanifolds with abelian complex structure, on which invariant forms compute both cohomologies. An equivariant isomorphism between the total Dolbeault and de Rham complexes gives equality in every Fr\"olicher inequality. It proves degeneration before the individual Hodge numbers are computed, see \cref{sec:cohomology}.

This comparison preserves total degree, but not the bigrading or complex conjugation. Accordingly, \cref{thm:main-degeneration} separates $E_1$-degeneration from both formality and Hodge symmetry in this simply connected holomorphic symplectic setting. For example, the standard fourfold has $h^{2,1}=1$ and $h^{1,2}=2$. Together with the triviality of the canonical bundle, the theorem also yields a cohomological proof of the smoothness of the Kuranishi space, previously established by Kurnosov--Verbitsky \cite[Corollary~5.10]{KurnosovVerbitsky2019}. For the standard fourfolds and sixfolds, this space has dimension four, see \cref{cor:unobstructed-deformations}.

\subsection{Lagrangian fibrations}\label{subsec:intro-reduction}

Write $m=n+1$. Guan's construction gives a resolution $f\colon Q_{n,t}\to Y=M/S_m$, which is projective in the complex analytic sense and semismall \cite[Lemma~5.6]{RiosOrtiz2026BGNonformal}. Here $M=M_{n,t}$ is a compact complex nilmanifold fibered over $B_E=\ker(E^m\xrightarrow{\Sigma}E)$ for an elliptic curve $E$. The map to $B_E/S_m\simeq|m\cdot0_E|\simeq\mathbb P^n$ induces a morphism $\ell\colon Q_{n,t}\to\mathbb P^n$.

In Bogomolov's description, let $S\to E$ be the primary Kodaira surface, with elliptic fiber $F$, and let $W\subset S^{[m]}$ be the fiber over $0$ of the summation map $S^{[m]}\to\operatorname{Sym}^mE\to E$, where $S^{[m]}$ is the Douady space of length-$m$ subschemes. Then $W_F=W/F$ is a compact complex holomorphic symplectic orbifold of dimension $2n$, obtained by the diagonal action of $F$. The induced map $\Pi\colon W_F\to|m\cdot0_E|$ and a finite map $p_W\colon Q_{n,t}\to W_F$ of degree $m^2$ give the algebraic reduction $\Phi=\Pi\circ p_W$ \cite{Bogomolov1996,BKKY2022,Kuznetsova2022Thesis}.

The nilmanifold $M_{n,t}$ is distinct from the precursor $N_n=(S^m|_{B_E})/F_\Delta$: there is an \'etale cover $M_{n,t}\to N_n$ of degree $m^2$. Their fundamental groups are torsion-free two-step nilpotent groups, each a central extension of $\Z^{2n}$ by $\Z^{2n}$, and $\pi_1(M_{n,t})$ embeds in $\pi_1(N_n)$ with index $m^2$. The resolution $Q_{n,t}$ is simply connected. The two bundles and the covering are described in \cref{sec:preliminaries}.

We prove that $\ell$ is a Lagrangian fibration with connected fibers of pure dimension $n$ and is itself an algebraic reduction. After an automorphism of the base, the two descriptions fit into the commutative diagram
\[
\begin{tikzcd}[column sep=small]
& Q_{n,t} \arrow[dl,"f"'] \arrow[dr,"p_W"] & \\
Y=M/S_m \arrow[dr,"\ell_Y"'] & & W_F \arrow[dl,"\Pi"] \\
& {|m\cdot0_E|\simeq\mathbb P^n}. &
\end{tikzcd}
\]
The construction and comparison are given in \cref{sec:lagrangian}. In particular, the hyperplane pullback is independent of the choice of description.

\begin{ThmIntr}\label{thm:main-lagrangian}
Let $Q=Q_{n,t}$, with $n\ge2$, and let $g\colon Q\to B$ be a nonconstant surjective holomorphic map with connected fibers onto a normal projective variety. Then there is a unique isomorphism $\psi\colon\mathbb P^n\to B$ with $g=\psi\circ\ell$. In particular, $\dim B=n$ and $g$ is a Lagrangian fibration. Under this identification, put $h=g^*c_1(\mathcal O_{\mathbb P^n}(1))$. Then $h\ne0$ in $H^2(Q,\R)$, whereas $h^2=0$ in $H^4(Q,\R)$. Moreover:
\begin{enumerate}
\item A general fiber of $g$, and the inverse image of a general codimension-two linear subspace of $\mathbb P^n$, are nonempty compact complex submanifolds that are null-homologous over $\Q$. The manifold $Q$ is not of Fujiki class~$\mathcal C$.
\item The class $h$ lies in the radical of every real quadratic form on $H^2(Q,\R)$ satisfying a Fujiki relation. It spans the radical of the Beauville--Bogomolov--Fujiki form, whose real inertia is $(2,3,1)$.
\item Every class in $H^2(Q,\C)$ restricts to zero on a general fiber of $g$.
\end{enumerate}
Here the inertia records the numbers of positive, negative, and zero directions, respectively, with the BBF form normalized to be positive on the real part of the holomorphic symplectic class.
\end{ThmIntr}

The uniqueness assertion uses the algebraic dimension $a(Q)=n$ and the equidimensionality of $\ell$. Together they imply that every holomorphic map from $Q$ to a projective variety factors through $\ell$. A nonconstant morphism from $\mathbb P^n$ onto a projective variety is finite onto its image; connectedness of the fibers of $g$ and normality of $B$ then make the induced map $\mathbb P^n\to B$ an isomorphism. Thus the conclusion determines both the map and its base; see \cref{thm:unique-projective-fibration}.

The identity $h^2=0$ has a direct geometric explanation. The divisors in $|m\cdot0_E|$ containing a point $z\in E$ form a hyperplane. Its preimage in $M$ is the sum of the divisors on which one coordinate lies over $z$. These divisors have the same rational cohomology class, while the divisors obtained by fixing the same coordinate over two distinct points of $E$ are disjoint. Their common class therefore has square zero. The argument is made explicit in \cref{thm:h-square,rem:geometric-square}. The Fujiki relation then forces $h$ into the radical: for every degree-two class $\alpha$, the polynomial $\int_Q(\alpha+s h)^{2n}$ has degree at most one in $s$, whereas it is a nonzero constant times the $n$-th power of $q(\alpha+s h)$.

On an irreducible hyper-K\"ahler manifold, the BBF form is nondegenerate of signature $(3,b_2-3)$ \cite{Beauville1983,Fujiki1987}; the class of a Lagrangian fiber cannot vanish. Guan had already proved that the BBF form of his examples has a one-dimensional radical \cite{GuanIII1995,Guan2023Survey,Guan2025Nonformal}. Our computation determines its full inertia and identifies its generator geometrically. It also explains why the nondegeneracy conjectured in \cite[Conjecture~5.11]{KurnosovVerbitsky2019}, and expected in \cite[Section~1.1]{BKKY2022}, fails.

The same class occurs in the topology of these manifolds. If $b$ is the central class in the Massey product of \cite[Section~3]{RiosOrtiz2026BGNonformal}, then $h=\kappa f^*[ib]$ for a positive constant $\kappa$. The annihilator of $h$ in $f^*H^2(Y,\C)$ is spanned by the pullbacks of the three classes entering that product and is a maximal totally isotropic subspace for the BBF form; see \cref{prop:massey,cor:isotropic}. Thus a single degree-two class links the vanishing of the fiber class, the radical of the BBF form, and the obstruction to formality.

\subsection{Characteristic numbers}

Let $\Kum_n(A)$ denote the generalized Kummer manifold of complex dimension $2n$ associated with an abelian surface $A$ \cite{Beauville1983}. Despite the preceding differences in cohomological geometry, its characteristic numbers provide the comparison for our numerical results.

\begin{ThmIntr}\label{thm:main}
Every $\BG_n$-type manifold $X$, with $n\ge2$, satisfies
\[
e(X)=(n+1)^3\sigma_1(n+1)=e(\Kum_n(A)),\qquad
\chi(\mathcal O_X)=n+1,
\]
where $\sigma_1(m)=\sum_{d\mid m}d$.
The standard $\BG_2$ and $\BG_3$ manifolds have the Hodge diamonds in \cref{thm:hodge,thm:six-hodge}. Their nonzero Chern numbers are
\[
\begin{aligned}
\BG_2:&\qquad c_2^2=756,\quad c_4=108,\\
\BG_3:&\qquad c_2^3=30208,\quad c_2c_4=6784,\quad c_6=448.
\end{aligned}
\]
These characteristic numbers hold throughout the corresponding deformation classes. For $n=2,3$, a $\BG_n$-type manifold and a $\Kum_n$-type manifold have the same Chern and Pontryagin numbers and represent the same unitary cobordism class.
\end{ThmIntr}

The generalized Kummer values are given in \cite{Gulbrandsen2007,Nieper2002}. The comparison does not extend to cohomology: $b_2(\BG_n)=6$, whereas $b_2(\Kum_n(A))=7$, and the standard $\BG_2$ and $\BG_3$ Hodge diamonds fail Hodge symmetry. The sector calculations determine the Hodge diamonds, and Hirzebruch--Riemann--Roch then determines all Chern numbers in these two dimensions. Their deformation invariance does not require the individual Hodge numbers to remain constant throughout a deformation class.

This agreement motivates \cref{ques:chern}, which predicts equality of all Chern, equivalently Pontryagin, numbers in every dimension. It is also consistent with the interpretation announced by Abasheva of the standard $\BG_n$ manifold as a Shafarevich--Tate twist of the natural Lagrangian fibration on $\Kum_n(E\times F)$ \cite[Remark~6.3.13]{Abasheva2026v5}. Under this announced identification, our computation places it in the case where $b_2$ drops by one and every degree-two class restricts trivially to a smooth fiber. The twist description does not by itself identify tangent-bundle characteristic classes; our proofs of the Euler formulas and the low-dimensional Chern numbers use the quotient construction directly.

\subsection{Organization of the paper}

\cref{sec:preliminaries,sec:cohomology} establish the covering, sector decomposition, and \cref{thm:main-degeneration}. \cref{sec:lagrangian,sec:bbf} study the Lagrangian class and the BBF form, proving \cref{thm:main-lagrangian}. \cref{sec:euler} computes both Euler characteristics in arbitrary dimension. \cref{sec:fourfold,sec:sixfold} give the fourfold and sixfold Hodge calculations. For the fourfold case we derive $e(Q)=9\cdot12=108$ from the cusp fibers of the Lagrangian fibration alone. \cref{sec:chern} determines the characteristic numbers and cobordism classes. Finally, \cref{sec:conjecture} formulates the Kummer conjecture and describes the additional characteristic data needed in higher dimensions.

\section{Semismall stratification of standard BG manifolds}\label{sec:preliminaries}

We fix a standard $\BG_n$-type manifold $Q=Q_{n,t}$, with $n\geq2$, and put $m=n+1$. The purpose of this section is to describe the geometric model underlying Guan's construction in a form adapted to the cohomological calculations below. We first construct the holomorphic $S_m$-equivariant covering of the nilmanifold appearing before resolution and make explicit its type-$A_n$ root-lattice structure. We then describe the stabilizer loci for the $S_m$-action, the resulting partition strata of the quotient $Y$, and their relation with the projective semismall resolution $f:Q\to Y$. These strata, together with their normalizations and top-cohomology local systems, are precisely the geometric contributions that enter the decomposition theorem and whose cohomology will be computed in \cref{sec:cohomology}. Unless a deformation is explicitly considered, all statements concerning Dolbeault cohomology refer to the standard complex structure.

We begin by recalling the precursor description of Guan manifolds following \cite[Section~5.1]{RiosOrtiz2026BGNonformal}. Let $E$ be an elliptic curve, let $t\in\Z\setminus\{0\}$, and let $\cL$ be a line bundle of degree $mt$ on $E$. Let $\cL^\times$ be the complement of the zero section. Choose $\gamma\in\C^\times$ with $|\gamma|>1$ and set $F=\C^\times/\langle\gamma\rangle$. Then $S:=\cL^\times/\langle\gamma\rangle$, with $\gamma$ acting by multiplication on the fibers, is a primary Kodaira surface and has a natural structure as a principal $F$-bundle over $E$. Let $n=m-1$ and define
\[
B_E:=\ker(E^m\xrightarrow{\Sigma}E),\qquad T_F:=\ker(F^m\xrightarrow{\Sigma}F).
\]
Thus $B_E\simeq E^n$ and $T_F\simeq F^n$. If $R$ denotes the restriction of $S^m$ to $B_E$, the diagonal action of $F$ on $R$ is free, and
\[
N_n:=R/F_\Delta
\]
is a principal $F^m/F_\Delta$-bundle over $B_E$. It is a compact complex nilmanifold by \cite[Claim~4.6]{AKV2026}.

It will be useful to make explicit the type-$A_n$ lattice structure underlying these two constructions. We first define
\[
Q_n:=Q(A_n)=\{(a_1,\ldots,a_m)\in\Z^m:\textstyle\sum_i a_i=0\},\qquad P_n:=P(A_n)=\Z^m/\Z(1,\ldots,1).
\]
Here $Q_n$ is the root lattice of type $A_n$, while $P_n$ is naturally identified with its weight lattice. The symmetric group $S_m=W(A_n)$ acts on both lattices by permuting the coordinates, and the natural map $Q_n\to P_n$ is injective with cokernel $\Z/m\Z$. If $\Lambda_F=H_1(F,\Z)$, then
\[
H_1(T_F,\Z)=Q_n\otimes\Lambda_F,\qquad H_1(F^m/F_\Delta,\Z)=P_n\otimes\Lambda_F.
\]
Consequently, the inclusion $T_F\subset F^m$ followed by the quotient $F^m\to F^m/F_\Delta$ is the isogeny induced by $Q_n\hookrightarrow P_n$. Its kernel is $F[m]$ and, since $\Lambda_F$ has rank two, its degree is $m^2$. This isogeny is $S_m$-equivariant.

The following construction globalizes this fiberwise isogeny, as in \cite[Proposition~5.3]{RiosOrtiz2026BGNonformal}.

\begin{Lem}\label{lem:holomorphic-cover}
There is a principal $T_F$-bundle $M_{n,t}\to B_E$ and an étale $S_m$-equivariant map $\rho:M_{n,t}\to N_n$ over $B_E$ of degree $m^2$. For $1\leq i\leq m$ there exist invariant smooth forms of type $(1,0)$, denoted by $x_i,y_i$, satisfying
\[
\sum_i x_i=\sum_i y_i=0,\qquad dx_i=0,\qquad dy_i=\sum_jx_j\bar x_j-mx_i\bar x_i.
\]
The forms $x_1,\ldots,x_n,y_1,\ldots,y_n$ form a coframe.
\end{Lem}

\begin{proof}
Choose a line bundle $\cL_0$ of degree $t$ and an isomorphism $\phi:\cL_0^{\otimes m}\simeq\cL$. Such a choice exists because multiplication by $m$ on $\mathrm{Pic}^0(E)$ is surjective. For the coordinate projections $p_i:B_E\to E$, set $\cN:=\bigotimes_i p_i^*\cL_0$ and $\cP_i:=p_i^*\cL\otimes\cN^{-1}$. Using $\phi$, we obtain a trivialization
\[
\bigotimes_i\cP_i\simeq\left(\bigotimes_i p_i^*\cL\right)\otimes\cN^{-m}\simeq\left(\bigotimes_i p_i^*\cL_0\right)^{\otimes m}\otimes\cN^{-m}\simeq\mathcal O_{B_E}.
\]

Let $\widetilde M_{n,t}$ be the locus of tuples $(v_1,\ldots,v_m)$, with $v_i\in(\cP_i)_b\setminus\{0\}$ over a common point $b\in B_E$, whose tensor product is $1$ under this trivialization. The root lattice $Q_n$ acts on $\widetilde M_{n,t}$ by
\[
(a_1,\ldots,a_m)\cdot(v_1,\ldots,v_m)=(\gamma^{a_1}v_1,\ldots,\gamma^{a_m}v_m),
\]
and we define $M_{n,t}:=\widetilde M_{n,t}/Q_n$.

To identify its fibers, choose local frames $e_i$ of $\cP_i$ satisfying $\bigotimes_i e_i=1$ and write $v_i=z_ie_i$. Then $\widetilde M_{n,t}$ is locally the product of the base with
\[
H:=\{(z_1,\ldots,z_m)\in(\C^\times)^m:z_1\cdots z_m=1\}.
\]
The quotient $\C^\times\to F$ induces a surjective homomorphism $H\to T_F$, given as 
\[
(z_1,\ldots,z_m)\mapsto([z_1],\ldots,[z_m]),
\]
whose kernel is precisely $\gamma^{Q_n}$. Hence $H/\gamma^{Q_n}\simeq T_F$. If $e_i'=g_ie_i$ is another compatible collection of frames, then $\prod_i g_i=1$ and the corresponding fiber coordinates differ by translation by $([g_1^{-1}],\ldots,[g_m^{-1}])\in T_F$. Thus $M_{n,t}\to B_E$ is a holomorphic principal $T_F$-bundle.

The permutation action of $S_m=W(A_n)$ preserves $Q_n$. It also permutes the factors of $\widetilde M_{n,t}$ and preserves both $\cN=\bigotimes_i p_i^*\cL_0$ and the tensor-product trivialization above. It therefore descends to an action of $S_m$ on $M_{n,t}$.

We now construct $\rho$. Let $U\subset B_E$ be an open set over which $\cN$ is trivial and choose a nowhere-vanishing frame $s$ of $\cN|_U$. Since $\cP_i\otimes\cN=p_i^*\cL$, each $v_i\in(\cP_i)_b$ determines
\[
u_i:=v_i\otimes s(b)\in(p_i^*\cL)_b.
\]
Taking the class of each $u_i$ modulo $\langle\gamma\rangle$ gives a holomorphic map $\widetilde M_{n,t}|_U\to R$. The $Q_n$-action does not change this point of $R$, since replacing $v_i$ by $\gamma^{a_i}v_i$ replaces $u_i$ by $\gamma^{a_i}u_i$. Hence the map descends to $M_{n,t}|_U\to R$.

We are left to show that this map does not depend on the frame. If the frame is replaced by $s'=fs$, with $f\in\mathcal O_U^\times$, then $u_i'=v_i\otimes s'=fu_i$ for every $i$. Thus changing the trivialization of $\cN$ multiplies all coordinates by the same element of $F$, so the resulting maps to $R$ differ by the diagonal action of $F_\Delta$. After passing to $N_n=R/F_\Delta$, this dependence disappears. The local maps therefore glue to a holomorphic map
\[
\rho:M_{n,t}\longrightarrow N_n
\]
over $B_E$. The construction is symmetric in the coordinates, hence $\rho$ is $S_m$-equivariant.

On every fiber, $\rho$ is precisely the isogeny $T_F\to F^m/F_\Delta$ induced by the inclusion $Q_n\hookrightarrow P_n$. By the discussion above, its kernel is $F[m]$ and its degree is $m^2$. Since the map on the base is the identity and its differential along every fiber is an isomorphism, $\rho$ is étale of degree $m^2$.

Finally, choose the invariant Kodaira coframe $x,\zeta$ with $dx=0$ and $d\zeta=x\bar x$, after rescaling $\zeta$, which is possible since $t\neq0$. Denote by $x_i,\zeta_i$ their pullbacks from the $i$-th factor of $S^m$. Since $x$ is pulled back from $E$ and $R$ lies over $B_E$, we have $\sum_i x_i=0$. Define $y_i:=\sum_j\zeta_j-m\zeta_i$. These forms are invariant under the diagonal $F$-action and vanish on its tangent direction, so they descend to $N_n$. Moreover,
\[
\sum_i y_i=0,\qquad dx_i=0,\qquad dy_i=\sum_jx_j\bar x_j-mx_i\bar x_i.
\]
The $x_i$ span the cotangent directions of $B_E$ with the single relation $\sum_i x_i=0$. On the vertical cotangent space, the map $(\zeta_1,\ldots,\zeta_m)\mapsto(y_1,\ldots,y_m)$ has kernel the diagonal line and image the hyperplane $\sum_i y_i=0$. Hence the descended $y_i$ span the vertical cotangent space with the single relation $\sum_i y_i=0$. Therefore $x_1,\ldots,x_n,y_1,\ldots,y_n$ form a coframe on $N_n$, and their pullbacks by $\rho$ give the required invariant coframe on $M_{n,t}$.
\end{proof}

\begin{Rem}[Fundamental groups]
Keep the notation $Q_n=Q(A_n)$, $P_n=P(A_n)$ and $\Lambda_F=H_1(F,\Z)$, and put $\Lambda_B=\pi_1(B_E)\simeq\Z^{2n}$. The two principal torus bundles give central extensions
\[
\begin{aligned}
1&\longrightarrow Q_n\otimes\Lambda_F
  \longrightarrow\pi_1(M_{n,t})\longrightarrow\Lambda_B\longrightarrow1,\\
1&\longrightarrow P_n\otimes\Lambda_F
  \longrightarrow\pi_1(N_n)\longrightarrow\Lambda_B\longrightarrow1.
\end{aligned}
\]
Indeed, $\pi_2(B_E)=0$, so the fiber fundamental groups inject, and the structure groups act by translations, so the monodromy on them is trivial. These extensions are torsion-free and two-step nilpotent; they are nonabelian for $n\ge2$ and $t\ne0$.

Their integral extension classes can be read directly from the line bundles in \cref{lem:holomorphic-cover}. Write $\eta_i=p_i^*[\mathrm{pt}]\in H^2(B_E,\Z)$ and $\eta=\sum_j\eta_j$. Choose the primitive element $\lambda\in\Lambda_F$ in the circle direction of $\C^\times\to F$, with sign chosen so that the characteristic class of $\cL^\times/\langle\gamma\rangle$ is $c_1(\cL)\otimes\lambda$. Then the extension class for $M_{n,t}$ is
\[
c_M=t\,(m\eta_1-\eta,\ldots,m\eta_m-\eta)\otimes\lambda
\in H^2(B_E,Q_n\otimes\Lambda_F),
\]
since $c_1(\cP_i)=t(m\eta_i-\eta)$ and these classes sum to zero. The class for $N_n$ is
\[
c_N=mt\sum_i\eta_i\otimes[e_i]\otimes\lambda
\in H^2(B_E,P_n\otimes\Lambda_F),
\]
where $[e_i]$ is the image of the $i$-th standard basis vector of $\Z^m$ in $P_n$. The map $Q_n\hookrightarrow P_n$ takes $c_M$ to $c_N$. Since $B_E$ is a torus, these classes specify the central extensions, and hence the two fundamental groups.

The cover $\rho$ induces the identity on $\Lambda_B$ and the inclusion $Q_n\otimes\Lambda_F\hookrightarrow P_n\otimes\Lambda_F$ on the kernels. Consequently,
\[
1\longrightarrow\pi_1(M_{n,t})\xrightarrow{\rho_*}\pi_1(N_n)
\longrightarrow F[m]\longrightarrow1,
\qquad F[m]\simeq(\Z/m\Z)^2.
\]
Equivalently, $\rho$ is the regular cover defined by the translation action of $F[m]$. Both nilmanifolds are aspherical. This description concerns the manifolds before the finite symmetric quotient and resolution; the resulting manifold $Q_{n,t}$ is simply connected by the Guan--Bogomolov construction.
\end{Rem}

We henceforth write $M=M_{n,t}$ and $Y=Y_{n,t}=M/S_m$. The form $\sigma=\sum_i x_i y_i$ is closed, since $\sum_i x_i=0$ and $x_i^2=0$. After eliminating $x_m$ and $y_m$, its coefficient matrix between the $x$- and $y$-families is $I_n+\mathbf 1\mathbf 1^t$, of determinant $m$. Thus $\sigma$ is nondegenerate. Being closed and of type $(2,0)$, it is holomorphic, and it is $S_m$-invariant. The construction gives a compact simply connected holomorphic symplectic manifold $Q=Q_{n,t}$ and a resolution $f:Q\longrightarrow Y$ \cite{GuanII1995,GuanIII1995,Bogomolov1996}.

\begin{Lem}\label{lem:projectivesemismall}
The resolution $f:Q\to Y$ is projective and semismall.
\end{Lem}

\begin{proof}
This is \cite[Lemma~5.6]{RiosOrtiz2026BGNonformal}.
\end{proof}

We now describe the strata of $Y$ in terms of the stabilizers of the $S_m$-action on $M$. Let
\[
V=\{(z_1,\dots,z_m)\in\C^m:\textstyle\sum_i z_i=0\}
\]
be the standard representation of $S_m$. By \cref{lem:holomorphic-cover}, at a point of $M$ with stabilizer $H\subseteq S_m$ the tangent representation is $(V\oplus V)|_H$.

The stabilizers that occur are Young subgroups. Let $\lambda=(\lambda_1,\ldots,\lambda_l)$ be a partition of $m$, choose a decomposition $\{1,\dots,m\}=B_1\sqcup\cdots\sqcup B_l$ with $|B_a|=\lambda_a$, and let $H_\lambda\subseteq S_m$ be the subgroup preserving each block. Thus $H_\lambda=\prod_a S_{\lambda_a}$, and different choices of blocks give conjugate subgroups. Set
\[
F_\lambda:=M^{H_\lambda},\qquad G_\lambda:=N_{S_m}(H_\lambda)/H_\lambda,\qquad r_\lambda:=m-l.
\]
The normalizer $N_{S_m}(H_\lambda)$ consists of permutations sending blocks to blocks, possibly exchanging blocks of equal size, so $G_\lambda\cong\prod_k S_{a_k}$, where $a_k$ is the number of parts of $\lambda$ equal to $k$. For a connected component $C$ of $F_\lambda$, let $G_{\lambda,C}\subseteq G_\lambda$ denote its stabilizer.

\begin{Prop}\label{prop:partition-sectors}
Every stabilizer of the $S_m$-action on $M$ is a Young subgroup. Every connected component of $F_\lambda$ contains a dense open subset of points with stabilizer exactly $H_\lambda$. The normalizations of the connected strata of $Y$ of type $\lambda$ are
\[
Z_{\lambda,C}=C/G_{\lambda,C},
\]
with one representative $C$ for each $G_\lambda$-orbit of connected components. Each $Z_{\lambda,C}$ has finite quotient singularities and dimension $2(l-1)$. The corresponding stratum of $Y$ is relevant for $f$, with $r=r_\lambda$, and its top-cohomology local system is the constant Tate local system $\Q(-r_\lambda)$.
\end{Prop}

\begin{proof}
We first determine the stabilizers of the $S_m$-action on $M$. Let $x\in M$ be represented, as in \cref{lem:holomorphic-cover}, by a tuple $(v_1,\ldots,v_m)$ over $b=(b_1,\ldots,b_m)\in B_E$. A permutation fixing $x$ must preserve the base coordinates $b_i$. If $b_i=b_j$, then $(\cP_i)_b$ and $(\cP_j)_b$ are canonically identified, since they have the same $\cL|_{b_i}$ factor and the same $\cN_b^{-1}$ factor. Under this identification, a permutation fixes $x$ precisely when it preserves the classes of the corresponding $v_i$ modulo $\langle\gamma\rangle$. Indeed, if $v_{\sigma(i)}=\gamma^{a_i}v_i$, then the fact that both tuples have tensor product $1$ gives $\gamma^{\sum_i a_i}=1$, hence $\sum_i a_i=0$. Thus the stabilizer consists of arbitrary permutations within the blocks on which both the base coordinates and the fiber classes agree, and is therefore a Young subgroup.

Fix $\lambda$ and a connected component $C$ of $F_\lambda=M^{H_\lambda}$. Since fixed loci of finite groups on complex manifolds are smooth, $C$ is smooth. At every $x\in C$,
\[
T_xC=(T_xM)^{H_\lambda}\simeq(V\oplus V)^{H_\lambda}.
\]
The subspace $V^{H_\lambda}$ consists of vectors constant on each block $B_a$, subject to the single relation $\sum_a\lambda_a c_a=0$, and hence has dimension $l-1$. Therefore $\dim C=2(l-1)$.

A generic vector in $V^{H_\lambda}$ takes distinct values on different blocks, so its stabilizer is exactly $H_\lambda$. Consequently, if $g\notin H_\lambda$, then $C\cap M^g$ is a proper analytic subset of $C$. Removing the finitely many such subsets gives a dense open subset $C^\circ\subset C$ whose points have stabilizer exactly $H_\lambda$.

Let $q:M\to Y=M/S_m$ be the quotient map. The restriction of $q$ induces a finite map $C/G_{\lambda,C}\to Y$. If $x,x'\in C^\circ$ and $x'=gx$ for some $g\in S_m$, then
\[
gH_\lambda g^{-1}=\operatorname{Stab}(x')=H_\lambda,
\]
so $g\in N_{S_m}(H_\lambda)$. Since $x'\in C\cap gC$ and $C,gC$ are connected components of $F_\lambda$, we have $gC=C$. Hence $C/G_{\lambda,C}\to Y$ is generically one-to-one onto the closure of $q(C^\circ)$. Since $C/G_{\lambda,C}$ is normal, it is the normalization of that closure. Two components determine the same stratum precisely when they lie in the same $G_\lambda$-orbit.

It remains to identify the fiber of $f$ over this stratum. The local analytic model of $f$ is the generalized Kummer resolution associated with the partition $\lambda$. Its fiber is
\[
\prod_{a=1}^l\Hilb^{\lambda_a}(\C^2)_0,
\]
where $\Hilb^k(\C^2)_0$ is the punctual Hilbert scheme of length $k$ subschemes supported at the origin. By \cite{Briancon1977}, each factor is irreducible of dimension $\lambda_a-1$, so the fiber is irreducible of dimension
\[
\sum_a(\lambda_a-1)=m-l=r_\lambda.
\]
On the other hand, $\dim Y=2(m-1)$ and $\dim Z_{\lambda,C}=2(l-1)$, so the corresponding stratum has codimension $2(m-l)=2r_\lambda$. Since $f$ is semismall by the preceding lemma, this equality shows that the stratum is relevant.

Finally, the fiber is irreducible of dimension $r_\lambda$, so its top cohomology is one-dimensional and generated by its fundamental class. The monodromy on this line is trivial, since there is only one top-dimensional irreducible component, and the generator has Hodge type $(r_\lambda,r_\lambda)$. Hence the top-cohomology local system is $\Q(-r_\lambda)$.
\end{proof}

\begin{Def}\label{def:partition-sector}
The \emph{partition sector} associated with $\lambda$ is $Z_\lambda=F_\lambda/G_\lambda$. It is the disjoint union
\[
Z_\lambda\simeq\bigsqcup_{[C]\in\pi_0(F_\lambda)/G_\lambda} Z_{\lambda,C},\qquad Z_{\lambda,C}=C/G_{\lambda,C}.
\]
Thus its connected sectors are indexed by the orbits of connected components of $F_\lambda$ under $G_\lambda$.
\end{Def}

Components exchanged by the normalizer give a single connected sector, and a component preserved by a subgroup contributes its quotient by that subgroup. A partition may therefore contribute several connected sectors. Their number depends on the global fixed locus and the normalizer action, and cannot be recovered from the tangent representation alone.

\begin{Ex}
For $m=3$, consider the partition $\lambda=(2,1)$. Here $H_\lambda=\langle(12)\rangle$ is its own normalizer, so $G_\lambda$ is trivial. Each connected component of $F_\lambda=M^{(12)}$ therefore gives a connected sector of dimension two and shift $r_\lambda=1$.

For $\lambda=(3)$, we have $H_\lambda=S_3$ and again $G_\lambda=1$, but the fixed locus is zero-dimensional. As computed in \cref{sec:fourfold}, it consists of $81$ points. Thus this single partition contributes $81$ connected sectors, each a point with shift $r_\lambda=2$. 
\end{Ex}

We conclude this section by bounding the connected components of fixed loci.

\begin{Lem}\label{lem:component-bound}
In the root-lattice bundle of \cref{lem:holomorphic-cover}, $F_\lambda$ is nonempty and has at most $\delta_\lambda^4$ connected components, where $\delta_\lambda=\gcd(\lambda_1,\ldots,\lambda_l)$. In particular, it is connected if $\delta_\lambda=1$.
\end{Lem}
\begin{proof}
The fixed locus on the base and the fixed subgroup of the structure group are
\[
\mathcal B_\lambda=\ker\left(E^l\xrightarrow{\sum_a\lambda_a}E\right),\qquad
T_\lambda=\ker\left(F^l\xrightarrow{\sum_a\lambda_a}F\right).
\]
Let $\delta_\lambda=\gcd(\lambda_1,\ldots,\lambda_l)$. Using the Smith normal form there exists $U\in\mathrm{GL}_l(\Z)$ such that
\[
(\lambda_1,\ldots,\lambda_l)U=(\delta_\lambda,0,\ldots,0).
\]
The matrix $U$ induces an automorphism of $E^l$, under which the homomorphism $(z_a)\mapsto\sum_a\lambda_a z_a$ becomes $(w_1,\ldots,w_l)\mapsto\delta_\lambda w_1$. Hence
\[
\mathcal B_\lambda\simeq E[\delta_\lambda]\times E^{l-1}.
\]
The same argument on $F^l$ gives $T_\lambda\simeq F[\delta_\lambda]\times F^{l-1}$. Since both $E[\delta_\lambda]$ and $F[\delta_\lambda]$ have order $\delta_\lambda^2$, each of $\mathcal B_\lambda$ and $T_\lambda$ has exactly $\delta_\lambda^2$ connected components.

The projection $M\to B_E$ restricts to a locally trivial $T_\lambda$-torsor $F_\lambda\to\mathcal B_\lambda$. Indeed, over $\mathcal B_\lambda$ the bundles $\cP_i$ are canonically identified within each block. Locally choose one nonzero frame $e_a$ for each block. After multiplying one of them by a local $\lambda_l$-th root, we may arrange that the weighted tensor product $\bigotimes_a e_a^{\otimes\lambda_a}$ is $1$. Repeating $e_a$ in every coordinate of the block $B_a$ then gives an $H_\lambda$-fixed local section. Once such a section is chosen, the $H_\lambda$-fixed points in each fiber form a translate of $T_F^{H_\lambda}=T_\lambda$.

Thus, over each connected component of $\mathcal B_\lambda$, the total space has at most $\delta_\lambda^2$ connected components, since monodromy can only permute the $\delta_\lambda^2$ connected components of $T_\lambda$. As $\mathcal B_\lambda$ itself has $\delta_\lambda^2$ connected components,
\[
|\pi_0(F_\lambda)|\leq\delta_\lambda^4.
\]
The local $H_\lambda$-fixed sections above show that $F_\lambda$ is nonempty. If $\delta_\lambda=1$, then both $\mathcal B_\lambda$ and $T_\lambda$ are connected, hence so is $F_\lambda$.
\end{proof}

\section{Cohomology and Frölicher degeneration}\label{sec:cohomology}

We now pass from the sectors of \cref{def:partition-sector}, constructed in \cref{sec:preliminaries}, to the cohomology of $Q$. We first explain the additive decomposition into sector contributions. We then compare the Dolbeault and de Rham complexes on each sector, obtaining $E_1$-degeneration in every dimension. Finally, we recover the low-degree Betti numbers from the identity and transposition sectors.

\subsection{De Rham complexes and the Hodge filtration}

We follow Fu's approach to the cohomology of Hilbert schemes of compact complex surfaces \cite[Sections~2--4]{Fu2024}. The basic idea is to apply the decomposition theorem at the level of Hodge modules and then pass to the associated graded pieces of their filtered de Rham complexes. The rational part of the decomposition controls ordinary cohomology, while the graded de Rham complexes recover holomorphic differential forms. We recall precisely the pieces of this formalism that will be used below.

\begin{Def}\label{def:holomorphic-de-rham}
Let $X$ be a complex manifold of dimension $d$. Its \emph{holomorphic de Rham complex} is
\[
\Omega_X^\bullet
=
\bigl[
\mathcal O_X
\xrightarrow{d}
\Omega_X^1
\xrightarrow{d}
\cdots
\xrightarrow{d}
\Omega_X^d
\bigr],
\]
placed in cohomological degrees $0,\ldots,d$. It carries the decreasing filtration by holomorphic degree
\[
F^p\Omega_X^\bullet
:=
\bigl[
0\longrightarrow\cdots\longrightarrow0
\longrightarrow\Omega_X^p
\xrightarrow{d}
\cdots
\xrightarrow{d}
\Omega_X^d
\bigr].
\]
Thus $\operatorname{gr}_F^p\Omega_X^\bullet\simeq\Omega_X^p[-p]$.
\end{Def}

The elementary relation between this complex and the usual cohomology of $X$ will be used repeatedly.

\begin{Lem}\label{lem:de-rham-frolicher}
Let $X$ be a complex manifold. Then the holomorphic de Rham complex is a resolution of the constant sheaf $\C_X$. Consequently,
\[
\mathbb H^k(X,\Omega_X^\bullet)\simeq H^k(X,\C),
\]
and the spectral sequence associated with the filtration of \cref{def:holomorphic-de-rham} is
\[
E_1^{p,q}=H^q(X,\Omega_X^p)
\Longrightarrow H^{p+q}(X,\C).
\]
This is the Frölicher spectral sequence. No Kähler hypothesis on $X$ is required.
\end{Lem}

\begin{proof}
The first assertion is the holomorphic Poincaré lemma. Taking hypercohomology gives the identification with singular cohomology, while the spectral sequence of the filtered complex has
\[
E_1^{p,q}
=
\mathbb H^{p+q}
\bigl(X,\operatorname{gr}_F^p\Omega_X^\bullet\bigr)
=
H^q(X,\Omega_X^p).
\qedhere
\]
\end{proof}

We denote by $\mathcal D_X$ the sheaf of holomorphic differential operators on $X$. A Hodge module combines a rational complex of sheaves with a compatible filtered $\mathcal D_X$-module; its filtered de Rham complex is the part of this data used to recover differential forms.

We next recall the de Rham complex appearing in Hodge-module theory. This should be distinguished from the complex of \cref{def:holomorphic-de-rham}, although for the constant Hodge module the two are canonically related.

\begin{Def}\label{def:Dmodule-de-rham}
Let $X$ be a complex manifold of dimension $d$, and let $(\mathcal M,F_\bullet)$ be a filtered right $\mathcal D_X$-module. Its \emph{de Rham complex} is
\[
\DR(\mathcal M)
:=
\bigl[
\mathcal M\otimes\wedge^dT_X
\longrightarrow
\mathcal M\otimes\wedge^{d-1}T_X
\longrightarrow\cdots\longrightarrow
\mathcal M\otimes T_X
\longrightarrow
\mathcal M
\bigr],
\]
placed in degrees $-d,\ldots,0$. The filtration on $\mathcal M$ induces an increasing filtration
\[
F_k\DR(\mathcal M)
=
\bigl[
F_{k-d}\mathcal M\otimes\wedge^dT_X
\longrightarrow\cdots\longrightarrow
F_{k-1}\mathcal M\otimes T_X
\longrightarrow
F_k\mathcal M
\bigr],
\]
and we write $\operatorname{gr}_k^F:=F_k/F_{k-1}$.
\end{Def}

For a Hodge module $\mathcal M$, we use $\DR(\mathcal M)$ to denote the de Rham complex of its underlying filtered right $\mathcal D$-module. The compatibility with ordinary holomorphic forms is contained in the following basic example.

\begin{Prop}\label{prop:constant-hodge-de-rham}
Let $X$ be a complex manifold of dimension $d$. The constant Hodge module $\Q_X^H[d]$ has rational complex $\Q_X[d]$ and underlying filtered right $\mathcal D_X$-module $(\omega_X,F_\bullet)$, with
\[
F_{-d}\omega_X=\omega_X,
\qquad
F_{-d-1}\omega_X=0.
\]
Its de Rham complex satisfies
\[
\DR(\omega_X)\simeq\Omega_X^\bullet[d],
\]
and, for every integer $p$,
\[
\operatorname{gr}_{-p}^F\DR(\Q_X^H[d])
\simeq
\Omega_X^p[d-p].
\]
In particular, the individual sheaves of holomorphic forms can be recovered from the associated graded pieces of the filtered de Rham complex.
\end{Prop}

\begin{proof}
This is the standard description of the constant Hodge module in the right $\mathcal D$-module convention as given in \cite[Section~2.2.1]{Fu2024}. Under the canonical identification $\DR(\omega_X)\simeq\Omega_X^\bullet[d]$, the filtration $F_{-p}$ corresponds to the truncation
\[
\bigl[
\Omega_X^p\longrightarrow\Omega_X^{p+1}
\longrightarrow\cdots\longrightarrow\Omega_X^d
\bigr],
\]
and taking the associated graded gives $\Omega_X^p[d-p]$.
\end{proof}

\begin{Rem}\label{rem:dr-conventions}
Our shift convention is $K[s]^j=K^{j+s}$, so a sheaf $\mathcal F[s]$ is concentrated in degree $-s$. The negative filtration index in \cref{prop:constant-hodge-de-rham} comes from the fact that Hodge modules are written with an increasing filtration $F_\bullet$, whereas the classical Hodge filtration on $\Omega_X^\bullet$ is decreasing.
\end{Rem}

The sectors appearing in our decomposition need not be smooth, so we also recall the corresponding statement for finite quotient singularities.

For a complex space $Z$, we write $d_Z=\dim_\C Z$.

\begin{Def}\label{def:reflexive-de-rham}
Let $Z$ be a normal complex space, let $j:Z_{\mathrm{reg}}\hookrightarrow Z$ be its smooth locus, and define
\[
\Omega_Z^{[p]}
:=
j_*\Omega_{Z_{\mathrm{reg}}}^p.
\]
These are the sheaves of \emph{reflexive differential forms}. If $Z=C/G$ is the quotient of a complex manifold $C$ by a finite group $G$, with quotient map $\pi:C\to Z$, then
\[
\Omega_Z^{[p]}
\simeq
(\pi_*\Omega_C^p)^G.
\]
For a space with finite quotient singularities, exterior differentiation on the smooth locus extends to a complex
\[
\Omega_Z^{[\bullet]}
=
\bigl[
\mathcal O_Z
\longrightarrow
\Omega_Z^{[1]}
\longrightarrow\cdots\longrightarrow
\Omega_Z^{[d_Z]}
\bigr],
\]
called the \emph{reflexive de Rham complex}.
\end{Def}

A rational local system is a sheaf of finite-dimensional rational vector spaces which is locally constant. Its fibers can be identified along paths, although continuation around a loop may give a nontrivial automorphism, called monodromy.

\begin{Def}\label{def:intersection-complex}
Let $W$ be a pure-dimensional complex space of dimension $e$, let $j:U\hookrightarrow W$ be a smooth dense open subset, and let $\mathcal L$ be a rational local system on $U$. The \emph{intersection complex} is
\[
\IC_W(\mathcal L):=j_{!*}(\mathcal L[e]),
\]
where $j_{!*}$ denotes intermediate, or middle, extension. This extension is taken in the category of perverse sheaves, whose objects are constructible complexes with dimension-dependent support conditions. Here constructibility means that the cohomology sheaves restrict to finite-rank local systems on the strata of a suitable stratification. The intermediate extension is characterized by restricting to $\mathcal L[e]$ on $U$ and having no nonzero perverse subobject or quotient supported on $W\setminus U$.

Our normalization is $\IC_W(\mathcal L)|_U=\mathcal L[e]$, and its hypercohomology defines intersection cohomology by
\[
IH^k(W,\mathcal L):=\mathbb H^{k-e}(W,\IC_W(\mathcal L)).
\]
Here $\mathbb H$ denotes the cohomology of a complex of sheaves. For constant coefficients on the smooth locus, we write $\IC_W$. In particular, if $W$ is smooth, then $\IC_W=\Q_W[e]$. We denote by $\IC_W^H$ the corresponding intersection-cohomology Hodge module, whose underlying rational complex is $\IC_W$.
\end{Def}

A complex space of dimension $e$ is a rational homology manifold if its local rational homology at every point is that of $\mathbb R^{2e}$. A complex space with finite quotient singularities is also called a $V$-manifold.

The Hodge module naturally associated with a singular sector is its intersection-cohomology Hodge module. For the quotient singularities occurring in this paper, this simplifies considerably.

\begin{Prop}\label{prop:quotient-hodge-de-rham}
Let $Z$ be a complex analytic space of dimension $d_Z$ with finite quotient singularities. Then $Z$ is a rational homology manifold, and hence
\[
\IC_Z\simeq\Q_Z[d_Z].
\]
Equivalently, the intersection-cohomology Hodge module agrees with the constant Hodge module,
\[
\IC_Z^H\simeq\Q_Z^H[d_Z].
\]
Moreover, for every integer $p$, there is a canonical isomorphism
\[
\operatorname{gr}_{-p}^F
\DR(\IC_Z^H)
\simeq
\Omega_Z^{[p]}[d_Z-p].
\]
\end{Prop}

\begin{proof}
A complex space with finite quotient singularities is a $V$-manifold and hence a rational homology manifold, so its intersection complex is $\Q_Z[d_Z]$. The statement about the associated graded de Rham complex is \cite[Proposition~3.3]{Fu2024}.
\end{proof}

Thus, by \cref{prop:quotient-hodge-de-rham}, ordinary cohomology replaces intersection cohomology for the sectors considered here, while the graded pieces of the filtered de Rham complex are precisely their reflexive differential forms.

\begin{Rem}\label{rem:singular-de-rham-intrinsic}
On a singular analytic space, the filtered $\mathcal D$-module underlying a Hodge module is defined locally after embedding the space into a complex manifold. Nevertheless, the complexes
\[
\operatorname{gr}_k^F\DR(\mathcal M)
\]
are intrinsic objects of the derived category of coherent sheaves and are independent of the chosen embedding; see \cite[Remarks~2.1 and~2.2]{Fu2024}. Thus \cref{prop:quotient-hodge-de-rham} can be applied intrinsically to each sector.
\end{Rem}

Combining \cref{prop:constant-hodge-de-rham,prop:quotient-hodge-de-rham}, the two sides of the Hodge-module decomposition have parallel interpretations: the rational complexes compute the ordinary cohomology of the resolution and of its sectors, whereas applying $\operatorname{gr}_{-p}^F\DR$ extracts the sheaves $\Omega^p$ on the smooth resolution and the reflexive sheaves $\Omega^{[p]}$ on the sectors. This is the mechanism that will allow us to pass from the sector decomposition to the Dolbeault groups of $Q$.

\subsection{Saito Decomposition Theorem}

The following form of the decomposition theorem is adapted to our sectors. A relevant stratum of codimension $2r$ will contribute with shift $2r$ in ordinary cohomology and shift $(r,r)$ in Dolbeault cohomology. We will work with Guan's semismall resolution $f:Q\to Y$ as in \cref{lem:projectivesemismall,prop:partition-sectors}.

Put $d=\dim_\C Q=\dim_\C Y$. For a stratification adapted to $f$, the restriction of $f$ over each smooth stratum $S_\alpha$ is topologically locally trivial. Semismallness means that $2\dim_\C f^{-1}(s)\leq\operatorname{codim}_Y S_\alpha$ for $s\in S_\alpha$, and a stratum is relevant when equality holds. If its codimension is $2r_\alpha$, its \emph{top-fiber local system} is
\[
\mathcal L_\alpha:=(R^{2r_\alpha}f_*\Q_Q)|_{S_\alpha},\qquad (\mathcal L_\alpha)_s=H^{2r_\alpha}(f^{-1}(s),\Q).
\]
We write $Rf_*$ for derived direct image; the cohomology sheaves of $Rf_*\Q_Q$ are $R^j f_*\Q_Q$, and
\[
\mathbb H^k(Y,Rf_*\Q_Q)=H^k(Q,\Q).
\]

The notation $\Q(-r)$ denotes the one-dimensional rational Hodge structure of type $(r,r)$. A constant variation $\Q(-r)$ is the constant local system with this same Hodge structure at every point. The \emph{Tate twist} $(-r)$ means tensoring with $\Q(-r)$, so it changes a Hodge type $(a,b)$ into $(a+r,b+r)$. With the increasing filtration used above, this gives
\[
\operatorname{gr}_{-p}^F\DR(\mathcal M(-r))\simeq\operatorname{gr}_{-(p-r)}^F\DR(\mathcal M).
\]

\begin{Lem}\label{lem:analytic-decomposition}
Let $f:Q\to Y$ as above. For every connected relevant stratum $S_\alpha$, put $2r_\alpha=\operatorname{codim}_Y S_\alpha$ and let $\nu_\alpha:Z_\alpha\to Y$ be the normalization of its closure, followed by inclusion. Assume that $Z_\alpha$ has finite quotient singularities and that its top-fiber local system is the constant variation $\Q(-r_\alpha)$, of rank one and Hodge type $(r_\alpha,r_\alpha)$. Then there are additive isomorphisms
\begin{align}
H^q(Q,\Omega_Q^p)&\simeq\bigoplus_\alpha H^{q-r_\alpha}(Z_\alpha,\Omega_{Z_\alpha}^{[p-r_\alpha]}),
\label{eq:general-dolbeault}\\
H^k(Q,\C)&\simeq\bigoplus_\alpha H^{k-2r_\alpha}(Z_\alpha,\C).
\label{eq:general-topological}
\end{align}
Groups with an out-of-range index are zero. If $Z=C/G$ is a finite quotient of a compact complex manifold, its contributions can be computed as
\[
H^q(Z,\Omega_Z^{[p]})\simeq H^q(C,\Omega_C^p)^G, \qquad H^k(Z,\C)\simeq H^k(C,\C)^G.
\]
\end{Lem}
\begin{proof}
The relatively ample line bundle of $f$ allows us to apply Saito's direct-image theorem in the complex analytic setting \cite[Theorem~(0.3)]{Saito1990}. For a semismall map, the decomposition is indexed by the relevant strata, with coefficients given by the top-fiber local systems and no further cohomological shifts. In the form explained in \cite[Theorem~4.2]{Fu2024}, our assumptions give
\[
Rf_*\Q_Q^H[d]\simeq
\bigoplus_\alpha\nu_{\alpha*}\IC_{Z_\alpha}^H(-r_\alpha).
\]
Here $Rf_*$ denotes the derived direct image. The twist records the Hodge type of the fiber cohomology. The normalizations may be used because intersection complexes are compatible with finite birational maps \cite[Lemma~2.6]{Fu2024}.

Since $Z_\alpha$ has finite quotient singularities, $\IC_{Z_\alpha}=\Q_{Z_\alpha}[d-2r_\alpha]$. After complexification and shifting by $[-d]$, we obtain
\[
Rf_*\C_Q\simeq\bigoplus_\alpha
\nu_{\alpha*}\C_{Z_\alpha}[-2r_\alpha].
\]
Taking hypercohomology gives \eqref{eq:general-topological}.

For Dolbeault cohomology we use Saito's strictness theorem, which says that taking a graded de Rham piece commutes with projective direct image \cite[Theorem~2.4 and Remark~2.5]{Fu2024}. Applying $\operatorname{gr}_{-p}^F\DR$ to the left side therefore gives $Rf_*\Omega_Q^p[d-p]$. On the $\alpha$-summand, the twist replaces $p$ by $p-r_\alpha$, and the reflexive-form formula gives
\[
\operatorname{gr}_{-p}^F\DR(\IC_{Z_\alpha}^H(-r_\alpha))
\simeq\Omega_{Z_\alpha}^{[p-r_\alpha]}[d-p-r_\alpha],
\]
because $\dim Z_\alpha=d-2r_\alpha$. Cancelling the common shift $[d-p]$ yields
\[
Rf_*\Omega_Q^p\simeq\bigoplus_\alpha
\nu_{\alpha*}\Omega_{Z_\alpha}^{[p-r_\alpha]}[-r_\alpha].
\]
Taking cohomology proves \eqref{eq:general-dolbeault}. These operations are compatible with local embeddings of the singular target \cite[Remarks~2.1, 2.2 and~2.5]{Fu2024}. Only projectivity of $f$ over $Y$ is required.

Finally, for $\pi:C\to Z=C/G$, one has $\Omega_Z^{[p]}=(\pi_*\Omega_C^p)^G$. Finite maps have no higher direct images on coherent sheaves. Moreover, averaging over $G$ makes taking invariants exact, so it commutes with cohomology. This proves the Dolbeault quotient formula. For ordinary cohomology, proper base change gives $R^j\pi_*\C_C=0$ for $j>0$, since the fibers are finite. Together with $(\pi_*\C_C)^G=\C_Z$, the same argument gives the topological formula.
\end{proof}

With the notation as in \cref{prop:partition-sectors}, the identity partition $(1^m)$ contributes $Y=M/S_m$ with shift zero. More generally, a sector $(Z,r)$ contributes degree $j$ to $H^{j+2r}(Q)$ and bidegree $(a,b)$ to bidegree $(a+r,b+r)$. For example, a point sector with $r=2$ contributes one class to $H^4(Q,\C)$, in Dolbeault bidegree $(2,2)$. Only sectors with $2r\leq k$ can contribute to $H^k(Q)$, and only those with $r\leq\min(p,q)$ can contribute to $H^q(Q,\Omega_Q^p)$.

These are additive decompositions, which do not specify cup products. They also do not by themselves prove degeneration of the Frölicher spectral sequence. For this we will compare the two cohomology calculations using the following elementary criterion.

\begin{Lem}\label{lem:frolicher-criterion}
For a compact complex manifold $X$, the Frölicher spectral sequence degenerates at $E_1$ if and only if $b_k(X)=\sum_{p+q=k}h^{p,q}(X)$ for every $k$.
\end{Lem}
\begin{proof}
This is the equality case of Frölicher inequalities, see \cite[Theorem 2]{Frolicher1955}.
\end{proof}

\subsection{Invariant forms}
\label{subsec:invariant-forms}

We now recall the invariant complexes on the nilmanifold appearing in \cref{lem:holomorphic-cover}. Let $M$ denote that nilmanifold and let $\mathfrak g$ be the real Lie algebra of its simply connected nilpotent covering group. The invariant complex structure $J$ on $M$ is determined by a complex structure on $\mathfrak g$, again denoted by $J$. After complexification one has $\mathfrak g_{\C}=\mathfrak g^{1,0}\oplus\mathfrak g^{0,1}$, and hence $\Lambda^k\mathfrak g_{\C}^*=\bigoplus_{p+q=k}\Lambda^{p,q}\mathfrak g^*$. We identify $\Lambda^{p,q}\mathfrak g^*$ with the space of invariant smooth differential forms of type $(p,q)$ on $M$.

\begin{Rem}\label{rem:invariant-not-holomorphic}
An invariant form of type $(1,0)$ need not be a holomorphic section of $\Omega_M^1$. In particular, the invariant forms $y_i$ occurring in \cref{lem:holomorphic-cover} are smooth forms of type $(1,0)$, but in general $\dbar y_i\neq0$. Here ``of type $(1,0)$'' refers only to the decomposition determined by the complex structure $J$.
\end{Rem}

Since $J$ is integrable, the Chevalley--Eilenberg differential on $\Lambda^\bullet\mathfrak g_{\C}^*$ decomposes as $d=\partial+\dbar$, where $\partial$ and $\dbar$ have bidegrees $(1,0)$ and $(0,1)$. Thus $(\Lambda^{p,\bullet}\mathfrak g^*,\dbar)$ is a subcomplex of the Dolbeault complex $\mathcal A^{p,\bullet}(M)$, while $(\Lambda^\bullet\mathfrak g_{\C}^*,d)$ is a subcomplex of the smooth de Rham complex $\mathcal A^\bullet(M)$.

The structure equations of \cref{lem:holomorphic-cover} satisfy $d(\mathfrak g^{1,0})^*\subset\Lambda^{1,1}\mathfrak g^*$. Equivalently, $\mathfrak g^{1,0}$ is abelian, so the invariant complex structure is abelian. Indeed, integrability already excludes the $(0,2)$-component of the differential of an invariant $(1,0)$-form, while the vanishing of the $(2,0)$-component is equivalent to $[\mathfrak g^{1,0},\mathfrak g^{1,0}]=0$. Consequently $\partial\alpha=0$ and $d\alpha=\dbar\alpha$ for every invariant $(1,0)$-form $\alpha$.

The point of this condition is that both the de Rham and Dolbeault cohomology of $M$ are computed by invariant forms.

\begin{Prop}\label{prop:invariant-cohomology}
Let $M$ be the compact nilmanifold of \cref{lem:holomorphic-cover}, endowed with its invariant abelian complex structure. Then the inclusions
\[
(\Lambda^\bullet\mathfrak g_{\C}^*,d)\hookrightarrow(\mathcal A^\bullet(M),d),
\qquad
(\Lambda^{p,\bullet}\mathfrak g^*,\dbar)\hookrightarrow(\mathcal A^{p,\bullet}(M),\dbar)
\]
are quasi-isomorphisms. Consequently,
\[
H^k(M,\C)\simeq H^k(\Lambda^\bullet\mathfrak g_{\C}^*,d),
\qquad
H^q(M,\Omega_M^p)\simeq H^q(\Lambda^{p,\bullet}\mathfrak g^*,\dbar)
\]
for all $p,q,k$.
\end{Prop}

\begin{proof}
The de Rham statement is Nomizu's theorem \cite{Nomizu1954}. The Dolbeault statement for nilmanifolds endowed with an invariant abelian complex structure is proved by Console--Fino \cite[Section~6, Remark on abelian complex structures]{ConsoleFino2001}. The identification $H^q(\mathcal A^{p,\bullet}(M),\dbar)\simeq H^q(M,\Omega_M^p)$ is the Dolbeault theorem.
\end{proof}

We now take $S_m$-invariants. By \cref{lem:holomorphic-cover}, the $S_m$-action preserves the invariant complex structure and the space of invariant differential forms, so the two quasi-isomorphisms above are $S_m$-equivariant. Since $S_m$ is finite and we work over $\C$, averaging makes the invariants functor exact, hence taking $S_m$-invariants commutes with cohomology. Therefore
\[
H^q(M,\Omega_M^p)^{S_m}\simeq H^q\bigl((\Lambda^{p,\bullet}\mathfrak g^*)^{S_m},\dbar\bigr), \qquad H^k(M,\C)^{S_m}\simeq H^k\bigl((\Lambda^\bullet\mathfrak g_{\C}^*)^{S_m},d\bigr).
\]
Thus the invariant Dolbeault and de Rham cohomology groups are computed on the same finite-dimensional exterior algebra, endowed respectively with the differentials $\dbar$ and $d$. This is the comparison that will be applied, component by component, to the sectors.
\begin{Lem}\label{lem:sector-comparison}
Let $C$ be a connected component of $F_\lambda$, and put $G=G_{\lambda,C}$ and $Z=C/G$. There is an isomorphism of total-degree differential graded $\C$-algebras between the invariant Dolbeault and de Rham models of $C$, equivariant for $G$. Consequently, for every $j$, we get
\[
\dim H^j(Z,\C)=\sum_{a+b=j} \dim H^b(Z,\Omega_Z^{[a]}).
\]
\end{Lem}
\begin{proof}
Write $\lambda=(\lambda_1,\ldots,\lambda_l)$. The assertion is immediate if $l=1$, when $C$ is a point. For $l>1$, restrict the global coframe of \cref{lem:holomorphic-cover} to $C$. Within each block the restricted forms agree, and we denote them by $X_a,Y_a,\bar X_a,\bar Y_a$. The relations $\sum_i x_i=\sum_i y_i=0$ become $\sum_a\lambda_aX_a=\sum_a\lambda_aY_a=0$, together with their conjugates. The weights count the original coordinates in each block. Eliminating its last member leaves a coframe, since the tangent space to the fixed locus is $(V\oplus V)^{H_\lambda}$. Its equations are
\[
dX_a=d\bar X_a=0,\qquad dY_a=F_a,\qquad d\bar Y_a=-F_a,\qquad F_a=\sum_b\lambda_b X_b\bar X_b-mX_a\bar X_a.
\]
The equation $d\bar Y_a=-F_a$ follows by conjugation, including the sign from interchanging the two degree-one factors.

The real and imaginary parts of this coframe have constant two-step nilpotent structure equations. Their dual vector fields are complete because $C$ is compact. They integrate to a transitive action of the simply connected nilpotent Lie group with this Lie algebra, with discrete cocompact stabilizer. Thus $C$ is a compact nilmanifold. Its complex structure is abelian, so invariant forms compute Dolbeault cohomology by \cite[Section~6]{ConsoleFino2001} and de Rham cohomology by \cite{Nomizu1954}. These comparison maps are equivariant for $G$: the action on the restricted coframe is obtained by restricting the global permutation action. Thus elements of $G$ permute equal-sized blocks, regardless of which component $C$ was chosen.

Let $\mathcal E_\lambda$ be the exterior algebra on the four families modulo their weighted relations, graded by total degree. Thus a form of bidegree $(a,b)$ has degree $a+b$, and both $d$ and $\dbar$ increase this degree by one. On this algebra $\dbar Y_a=F_a$, and $\dbar$ vanishes on $X_a,\bar X_a,\bar Y_a$. Define an algebra automorphism $T$ by
\[
T(X_a)=X_a,\qquad T(\bar X_a)=\bar X_a,\qquad T(Y_a)=Y_a,\qquad T(\bar Y_a)=Y_a+\bar Y_a.
\]
It respects the weighted relations: for example, $\sum_a\lambda_aT(\bar Y_a)= \sum_a\lambda_aY_a+\sum_a\lambda_a\bar Y_a=0$. Its inverse fixes $X_a,\bar X_a,Y_a$ and sends $\bar Y_a$ to $\bar Y_a-Y_a$. It also commutes with permutations of equal-sized blocks, hence with the action of the component stabilizer $G$. To check the chain identity, note that $T(F_a)=F_a$ and
\[
dT(Y_a)=F_a=T\dbar Y_a,\qquad dT(\bar Y_a)=d(Y_a+\bar Y_a)=0=T\dbar\bar Y_a.
\]
Both differentials vanish on $X_a$ and $\bar X_a$. The Leibniz rule therefore gives $dT=T\dbar$ on the entire exterior algebra. It follows that
\[
(\mathcal E_\lambda^G,\dbar) \xrightarrow{\ T\ }(\mathcal E_\lambda^G,d)
\]
is an isomorphism of total-degree complexes. In degree $j$, the cohomology groups on the left and right are, respectively,
\[
\bigoplus_{a+b=j}H^b(C,\Omega_C^a)^G
\qquad\text{and}\qquad H^j(C,\C)^G.
\] \cref{lem:analytic-decomposition} identifies these with the asserted sector groups.
\end{proof}

The comparison on each sector is independent of its number of components. It therefore proves degeneration before the individual sector multiplicities or Hodge numbers are known, this yields the main result in this section.

\begin{Thm}\label{thm:all-degeneration}
For every standard $\BG_n$-type manifold $Q=Q_{n,t}$, with $n\ge2$, the Frölicher spectral sequence degenerates at $E_1$.
\end{Thm}

\begin{proof}
Let $Z_{\lambda,C}$ be the connected sectors of \cref{prop:partition-sectors}, where $\lambda$ runs over partitions of $m=n+1$ and $[C]$ over the corresponding component orbits. By \cref{lem:analytic-decomposition}, the cohomology of $Q$ decomposes additively as
\[
H^k(Q,\C)\simeq\bigoplus_{\lambda,[C]}H^{k-2r_\lambda}(Z_{\lambda,C},\C), \qquad H^q(Q,\Omega_Q^p)\simeq\bigoplus_{\lambda,[C]}H^{q-r_\lambda} \bigl(Z_{\lambda,C},\Omega_{Z_{\lambda,C}}^{[p-r_\lambda]}\bigr).
\]
As usual, a summand with a negative index is understood to be zero. Hence, writing $h^{a,b}(Z):=\dim H^b(Z,\Omega_Z^{[a]})$, we obtain
\[
b_k(Q)=\sum_{\lambda,[C]}b_{k-2r_\lambda}(Z_{\lambda,C}), \qquad \sum_{p+q=k}h^{p,q}(Q) =\sum_{\lambda,[C]}\sum_{p+q=k} h^{p-r_\lambda,q-r_\lambda}(Z_{\lambda,C}).
\]
For every connected sector $Z_{\lambda,C}$ and every integer $\ell$, by \cref{lem:sector-comparison} we get
\[
b_\ell(Z_{\lambda,C}) = \sum_{a+b=\ell}h^{a,b}(Z_{\lambda,C}).
\]
Applying this with $\ell=k-2r_\lambda$ yields
\[
b_k(Q) = \sum_{\lambda,[C]} \sum_{a+b=k-2r_\lambda} h^{a,b}(Z_{\lambda,C}).
\]
For each fixed $\lambda$, set $p=a+r_\lambda$ and $q=b+r_\lambda$. Then $a+b=k-2r_\lambda$ is equivalent to $p+q=k$, and therefore
\[
b_k(Q) = \sum_{\lambda,[C]}\sum_{p+q=k} h^{p-r_\lambda,q-r_\lambda}(Z_{\lambda,C}) = \sum_{p+q=k}h^{p,q}(Q).
\]
Thus the Frölicher inequalities are equalities in every total degree $k$. By \cref{lem:frolicher-criterion}, the Frölicher spectral sequence of $Q$ degenerates at $E_1$.
\end{proof}

\begin{Cor}\label{cor:small-deformations}
Every sufficiently small deformation of $Q_{n,t}$ has the same Hodge numbers as $Q_{n,t}$ and has Frölicher spectral sequence degenerating at $E_1$.
\end{Cor}
\begin{proof}
In a smooth proper deformation, Betti numbers are locally constant and each Hodge number is upper semicontinuous. There are only finitely many bidegrees, so we may shrink the base so that every $h^{p,q}$ is at most its value on $Q_{n,t}$. The Frölicher inequality gives $b_k\leq\sum_{p+q=k}h^{p,q}$ on every fiber, while \cref{thm:all-degeneration} gives equality on the central fiber. Hence equality holds term by term nearby, and \cref{lem:frolicher-criterion} applies again.
\end{proof}

The Kuranishi family parametrizes all sufficiently small deformations of a compact complex manifold. Smoothness of the Kuranishi space is equivalent to unobstructedness of small deformations, while universality means that every sufficiently small deformation is induced from it by a unique morphism of germs.

\begin{Cor}\label{cor:unobstructed-deformations}
Let $Q=Q_{n,t}$ be a standard $\BG_n$-type manifold, with $n\ge2$. Then the Kuranishi space of $Q$ is smooth of dimension $h^{1,1}(Q)$, and the Kuranishi family is universal. Moreover, every sufficiently small deformation of $Q$ remains holomorphic symplectic. For $n=2,3$, the Kuranishi space has dimension four.
\end{Cor}

\begin{proof}
The holomorphic symplectic form on $Q$ trivializes $K_Q$, and \cref{thm:all-degeneration} shows that the Frölicher spectral sequence of $Q$ degenerates at $E_1$. Hence \cite[Theorem~3.3 and Corollary~3.5]{AnthesCattaneoRollenskeTomassini2018} applies and yields a smooth universal Kuranishi family. Its tangent space at the origin is $H^1(Q,T_Q)$. Contraction with the symplectic form gives an isomorphism $T_Q\simeq\Omega_Q^1$, and therefore
\[
\dim\mathrm{Def}(Q)=h^1(Q,T_Q)=h^1(Q,\Omega_Q^1)=h^{1,1}(Q).
\]
By \cite[Proposition~3.1]{AnthesCattaneoRollenskeTomassini2018}, the holomorphic symplectic form persists under sufficiently small deformations, so every nearby fibre is again holomorphic symplectic. Finally, \cref{thm:hodge,thm:six-hodge} give $h^{1,1}(Q)=4$ for $n=2,3$.
\end{proof}

The smoothness of the Kuranishi space and the persistence of the holomorphic symplectic structure recover \cite[Corollary~5.10]{KurnosovVerbitsky2019}. In the present approach, however, they arise from the sectorwise comparison of de Rham and Dolbeault cohomology together with the general deformation-theoretic criterion above.

\begin{Rem}
The argument above proves only $E_1$-degeneration. The automorphism $T$ used in the sector comparison does not preserve the bigrading and is not compatible with complex conjugation, so it does not imply Hodge symmetry. Likewise, the additive sector decomposition is an isomorphism of graded vector spaces rather than of differential graded algebras or cohomology rings, and therefore gives no $\partial\bar\partial$-lemma. These limitations are essential: the resulting Hodge-theoretic behavior is compatible with the nonformality of $\BG_n$-type manifolds.   
\end{Rem}

We conclude by recording the low-degree Betti numbers, previously computed in \cite{GuanII1995,RiosOrtiz2026BGNonformal}, while making explicit how the transposition sector contributes to $b_2$ and $b_3$.

\begin{Prop}\label{prop:low}
Every $\BG_n$-type manifold with $n\ge2$ satisfies
\[
b_1=0,\qquad b_2=6,\qquad b_3=3.
\]
\end{Prop}

\begin{proof}
We first consider the standard representative $Q=Q_{n,t}$. By Guan's construction \cite{GuanII1995}, $Q$ is simply connected, hence $b_1(Q)=0$.

For the remaining two Betti numbers, only the identity sector and the transposition sector can contribute. Indeed, if $\tau=(12)$, the partition $\lambda=(2,1^{m-2})$ has shift $r_\lambda=1$, and it is the unique nonidentity partition with this property; every other nonidentity sector has shift at least two and therefore contributes only from degree four onwards. Denoting the identity sector by $Y$ and the transposition sector by $Z$, \cref{lem:analytic-decomposition} therefore gives, for $k=2,3$,
\[
H^k(Q,\C)\simeq H^k(Y,\C)\oplus H^{k-2}(Z,\C).
\]

We first determine the connectedness of $Z$. The partition $(2,1^{m-2})$ has gcd equal to one, so \cref{lem:component-bound} implies that $M^\tau$ is connected. The residual group is $S_{m-2}$, acting by permutation of the singleton blocks, and hence $Z=M^\tau/S_{m-2}$ is connected. In particular, $b_0(Z)=1$.

The identity-sector calculation of \cite[Corollary~5.5]{RiosOrtiz2026BGNonformal} gives
\[
b_2(Y)=5,\qquad b_3(Y)=0,
\]
while \cite[Lemma~5.12]{RiosOrtiz2026BGNonformal} gives $b_1(Z)=3$. We briefly recall the latter calculation. Set $w_i=\bar x_i$, $u_i=y_i+\bar y_i$ and $v_i=(y_i-\bar y_i)/2$. Then the families $x_i,w_i,u_i$ consist of closed invariant one-forms, whereas
\[
dv_i=\sum_jx_jw_j-mx_iw_i.
\]
After restriction to $M^\tau$ and passage to $S_{m-2}$-invariants, each of the four families $x,w,u,v$ contributes a one-dimensional invariant subspace in degree one. The classes represented by the $x$-, $w$- and $u$-families are closed, whereas the invariant form arising from the $v$-family is not. Since there are no degree-one boundaries, this gives $b_1(Z)=3$.

The preceding decomposition now yields
\[
b_2(Q)=b_2(Y)+b_0(Z)=5+1=6,\qquad
b_3(Q)=b_3(Y)+b_1(Z)=0+3=3.
\]
Betti numbers are invariant under smooth proper deformation. Since every $\BG_n$-type manifold belongs to the deformation type of a standard representative, the same values hold throughout the $\BG_n$ deformation class.
\end{proof}

The odd Betti number $b_3=3$ gives an immediate topological obstruction to both Kählerness and the $\partial\bar\partial$-lemma.

\begin{Cor}\label{cor:no-kahler-ddbar}
Let $X$ be the underlying smooth manifold of a $\BG_n$-type manifold, with $n\ge2$. Then $X$ admits neither a Kähler complex structure nor a complex structure satisfying the $\partial\bar\partial$-lemma.
\end{Cor}

\begin{proof}
Suppose that a complex structure on $X$ satisfies the $\partial\bar\partial$-lemma. Then de Rham cohomology admits the Hodge decomposition
\[
H^k(X,\C)=\bigoplus_{p+q=k}H^{p,q}(X),
\]
with complex conjugation identifying $H^{p,q}(X)$ and $H^{q,p}(X)$. If $k$ is odd, no pair $(p,q)$ with $p+q=k$ satisfies $p=q$, so the summands occur in distinct conjugate pairs. Hence every odd Betti number is even. This contradicts \cref{prop:low}, which gives $b_3(X)=3$. Thus $X$ admits no complex structure satisfying the $\partial\bar\partial$-lemma. Since every compact Kähler manifold satisfies the $\partial\bar\partial$-lemma, it admits no Kähler complex structure either.
\end{proof}

\section{Lagrangian Fibrations}\label{sec:lagrangian}

Let $Q=Q_{n,t}$ with $n\ge2$ and $m=n+1$, and keep the notation of \cref{sec:preliminaries}. Write $\varpi\colon M\to B_E$ for the bundle projection and $\pi\colon M\to Y$ for the quotient map. We identify the natural fibration induced by $M\to B_E$ with Bogomolov's algebraic reduction, then prove that its hyperplane pullback $h$ is nonzero and has square zero. We also use this to show that every holomorphic Lagrangian fibration on $Q$ with normal projective base is isomorphic to the natural one. The consequences for the BBF form will be developed in \cref{sec:bbf}.

Let $\pi_S\colon S\to E$ be the primary Kodaira surface of \cref{sec:preliminaries}, and let $W\subset S^{[m]}$ be the fiber over $0$ of the composite $S^{[m]}\to\Sym^mS\to\Sym^mE\to E$ of the Douady--Barlet morphism, the projection and the summation. The diagonal action of $F$ preserves $W$ with finite stabilizers, and $W_F=W/F$ is a compact complex orbifold of dimension $2n$ with a holomorphic symplectic form \cite{Bogomolov1996}. By Abel's theorem, the fiber over $0$ of $\Sym^mE\to E$ is the linear system $|m\cdot0_E|\simeq\mathbb P^n$. We write $\Pi\colon W_F\to\mathbb P^n$ for the induced map.

Bogomolov's construction gives a finite map $p_W\colon Q\to W_F$ of degree $m^2$ \cite[Section~3]{Bogomolov1996}, \cite[Theorem~3.3]{BKKY2022}, identifying his manifold with Guan's $Q_{n,t}$ \cite[Theorem~1]{GuanIII1995}, see also \cite[Remark~3.3.7]{Kuznetsova2022Thesis}. The composite
\[
\Phi=\Pi\circ p_W\colon Q\longrightarrow|m\cdot0_E|\simeq\mathbb P^n
\]
is the algebraic reduction of $Q$, and its general fiber is an abelian variety of dimension $n$ \cite[Corollary~4.3 and Proposition~3.6]{BKKY2022}.  The construction fits into the commutative diagram
\[
\begin{tikzcd}[column sep=large, row sep=large]
& W \arrow[r, hook] \arrow[d, "{/F}"'] & S^{[m]} \arrow[d, "\delta"] \\
Q \arrow[r, "p_W"] \arrow[dr, "\Phi"'] & W_F \arrow[d, "\Pi"] & \Sym^m S \arrow[d, "\Sym^m(\pi_S)"] \\
& {\lvert m\cdot0_E\rvert\simeq\mathbb P^n} \arrow[r, hook] \arrow[d] & \Sym^m E \arrow[d, "\Sigma"] \\
& \{0\} \arrow[r, hook] & E
\end{tikzcd}
\]
where $\delta$ is the Douady--Barlet morphism and $\Sigma$ is the summation map.

Since $\Phi$ is the algebraic reduction, $\Phi^*\C(\mathbb P^n)=\mathcal M(Q)$, and therefore $a(Q)=n$. The special fibers lie over the discriminant $\cD\subset\mathbb P^n$ of nonreduced divisors, a hypersurface of degree $2m$. For $n=2$ it is the curve dual to the plane cubic $E\subset|3\cdot0_E|^\vee$, a sextic with nine cusps at the divisors $3[z]$, $z\in E[3]$ \cite[Lemma~3.5.1, Propositions~3.5.7 and~3.5.8, and Example~3.5.2]{Kuznetsova2022Thesis}. Restricting the quotient map $E^m\to\Sym^mE$ to $B_E$ gives
\[
\alpha_E\colon B_E\longrightarrow|m\cdot0_E|,\qquad \mathbf z=(z_i)\longmapsto\sum_i[z_i],
\]
which identifies $B_E/S_m$ with $|m\cdot0_E|$. The map $\alpha_E\circ\varpi$ is $S_m$-invariant, hence descends to a map $\ell_Y\colon Y\to|m\cdot0_E|$.

\begin{Def}\label{def:lagrangian}
A holomorphic Lagrangian fibration on a compact holomorphic symplectic manifold $X$ of dimension $2n$, we mean a surjective holomorphic map $g\colon X\to B$ onto a normal compact complex space of dimension $n$, with connected fibers, whose general fiber is smooth and Lagrangian.   
\end{Def}

With the notation as above, the \emph{natural} Lagrangian fibration of $Q$ is $\ell=\ell_Y\circ f\colon Q\to|m\cdot0_E|\simeq\mathbb P^n$. Put $h=\ell^*c_1(\mathcal O_{\mathbb P^n}(1))\in H^2(Q,\Z)$ and $h_Y=\ell_Y^*c_1(\mathcal O_{\mathbb P^n}(1))\in H^2(Y,\Q)$, so that $h=f^*h_Y$. By the explicit description it is Lagrangian in the sense of \cref{def:lagrangian}.

The comparison with Bogomolov's fibration uses the following factorization statement. Here $a(X)=\operatorname{trdeg}_{\C}\mathcal M(X)$ denotes the algebraic dimension.

\begin{Lem}\label{lem:equidimensional-factorization}
Let $X$ be a connected compact complex manifold of dimension $d$, and let $u\colon X\to B$ be a surjective holomorphic map with connected fibers, all of pure dimension $d-r$. Assume that $B$ is a normal projective variety of dimension $r=a(X)$. Then every holomorphic map $v\colon X\to Z$ to a projective variety factors uniquely as $v=\psi\circ u$, with $\psi\colon B\to Z$ holomorphic.
\end{Lem}

\begin{proof}
Let $I\subset B\times Z$ be the image of $(u,v)$, with its reduced structure. It is an irreducible projective variety, and pullback embeds $\C(I)$ into $\mathcal M(X)$. Hence $\dim I\le a(X)=r$. The first projection is surjective onto $B$, so $\dim I=r$.

Every fiber of $X\to I$ has dimension at least $d-r$. Since $u^{-1}(b)$ has pure dimension $d-r$, its image $I_b$ under $v$ must be finite. It is connected because $u^{-1}(b)$ is connected. Thus the first projection $I\to B$ is finite and birational, and is an isomorphism since $B$ is normal. The second projection gives $\psi$, whose uniqueness follows from the surjectivity of $u$.
\end{proof}

\begin{Prop}\label{prop:lagrangian}
\begin{enumerate}
\item Over $\mathbb P^n\setminus\cD$, every fiber of $\ell$ is mapped isomorphically by $f$ onto a fiber of $\ell_Y$ contained in the smooth locus of $Y$, and is isomorphic to the complex torus $T_F\simeq F^n$. These fibers are Lagrangian for the holomorphic symplectic form $\sigma_Q$.
\item The map $\ell$ is surjective, and all its fibers are connected and of pure dimension $n$.
\item The map $\ell$ is an algebraic reduction of $Q$. For any identification of $Q$ with Bogomolov's manifold, $\Phi=\psi\circ\ell$ for an automorphism $\psi$ of $\mathbb P^n$ with $\psi(\cD)=\cD$. In particular $h=\Phi^*c_1(\mathcal O_{\mathbb P^n}(1))$.
\end{enumerate}
\end{Prop}

\begin{proof}
(1) Let $y\notin\cD$ and $\mathbf z\in\alpha_E^{-1}(y)$. The coordinates of $\mathbf z$ are pairwise distinct, so every point of $\varpi^{-1}(\mathbf z)$ has trivial stabilizer, and $S_m$ acts simply transitively on $\alpha_E^{-1}(y)$. Hence $\pi$ identifies $\varpi^{-1}(\mathbf z)\simeq T_F$ with $\ell_Y^{-1}(y)$, which lies in the smooth locus of $Y$, over which $f$ is an isomorphism. By \cref{lem:holomorphic-cover}, the $x_i$ are pulled back from $B_E$. Thus $\sigma=\sum_ix_iy_i$ vanishes on these fibers, whose dimension is $n$, so they are Lagrangian.

(2) The map $\alpha_E$ is surjective, and $S_m$ acts transitively on each of its fibers. Thus every fiber of $\ell_Y$ is the image of a single fiber of $\varpi$, hence connected. Since $f$ is a resolution of the normal space $Y$, its fibers are connected, and the same holds for $\ell$.

For the dimension, the proof of \cref{lem:component-bound} identifies the fibers of $F_\lambda\to\mathcal B_\lambda$ with translates of $T_\lambda=\ker(\sum_a\lambda_a\colon F^l\to F)$, of dimension $l-1=n-r_\lambda$. Hence $\ell_Y^{-1}(y)$ meets the stratum of type $\lambda$ in dimension at most $n-r_\lambda$. By \cref{prop:partition-sectors}, the fibers of $f$ over this stratum have dimension $r_\lambda$. Therefore $\dim\ell^{-1}(y)\le n$, while every component has dimension at least $\dim Q-\dim\mathbb P^n=n$.

(3) By (2) and $a(Q)=n$, \cref{lem:equidimensional-factorization} applied to $u=\ell$ and $v=\Phi$ gives $\Phi=\psi\circ\ell$ for a morphism $\psi\colon\mathbb P^n\to\mathbb P^n$. Since $\Phi^*\C(\mathbb P^n)=\mathcal M(Q)$, the inclusions
\[
\Phi^*\C(\mathbb P^n)\subseteq\ell^*\C(\mathbb P^n)\subseteq\mathcal M(Q)
\]
are equalities. Thus $\ell$ is an algebraic reduction and $\psi$ is birational. Writing $\psi^*\mathcal O(1)=\mathcal O(k)$, one has $k\ge1$ and $k^n=\deg\psi=1$, so $\psi$ is an automorphism.

For a general point $y\in\cD$, \cref{prop:partition-sectors} shows that $\ell^{-1}(y)$ contains rational curves over the transposition stratum. The fibers of $\Phi$ outside $\cD$ are complex tori \cite[Lemma~3.4.18]{Kuznetsova2022Thesis}, so $\psi(\cD)\subseteq\cD$. Equality follows because both are reduced hypersurfaces of the same degree.
\end{proof}

The same argument applies to arbitrary fibrations with projective base.

\begin{Thm}\label{thm:unique-projective-fibration}
Let $Q=Q_{n,t}$, with $n\ge2$, and let $g\colon Q\to B$ be a nonconstant surjective holomorphic map with connected fibers onto a normal projective variety. Then there is a unique isomorphism $\psi\colon\mathbb P^n\to B$ such that $g=\psi\circ\ell$. In particular $\dim B=n$, and $g$ is a Lagrangian fibration. Thus every holomorphic Lagrangian fibration on $Q$ with normal projective base is the natural one, up to an isomorphism of the base.
\end{Thm}

\begin{proof}
By \cref{prop:lagrangian,lem:equidimensional-factorization}, one has $g=\psi\circ\ell$ for a surjective morphism $\psi\colon\mathbb P^n\to B$. Let $A$ be an ample line bundle on $B$ and write $\psi^*A=\mathcal O_{\mathbb P^n}(k)$. Since $\psi$ is nonconstant, $k>0$. A positive-dimensional fiber of $\psi$ would contain a projective curve on which $\psi^*A$ has degree zero, contradicting the ampleness of $\mathcal O(k)$. Thus $\psi$ is finite. Its fibers are connected, being the images under $\ell$ of the fibers of $g$, so $\psi$ is birational. Normality of $B$ implies that $\psi$ is an isomorphism. The remaining assertions follow from \cref{prop:lagrangian}.
\end{proof}

For a general point $y\in\cD$, the fiber $\ell_Y^{-1}(y)$ is the quotient of $T_F$ by a transposition, which is a $\mathbb P^1$-bundle over $F^{n-1}$. The fiber $\ell^{-1}(y)$ consists of a component birational to it and a $\mathbb P^1$-bundle over the image of the fixed locus. Thus all its components are uniruled, in accordance with \cite[Lemma~3.4.18]{Kuznetsova2022Thesis}. For $n=2$, the fibers over the nine cusps of $\cD$ are described in \cref{subsec:cusp-fibers}.

\begin{Lem}\label{lem:hyperplane-pullback}
For $z\in E$, the divisors of $|m\cdot0_E|$ containing $z$ form a hyperplane $H_z$, and $\alpha_E^*H_z=\sum_ip_i^{-1}(z)$, each component with multiplicity one. Put $b=\sum_ix_i\bar x_i$. Then $[x_i\bar x_i]=[x_j\bar x_j]$ in $H^2(M,\C)$ for all $i,j$, and
\[
\pi^*h_Y=\varpi^*\alpha_E^*c_1(\mathcal O_{\mathbb P^n}(1))=\kappa\,[\,ib\,]\qquad\text{for some }\kappa>0.
\]
\end{Lem}

\begin{proof}
The hyperplane $H_z\subset\mathbb P(H^0(E,\mathcal O_E(m\cdot0_E)))$ is the zero locus of evaluation at $z$, a nonzero linear form because $\mathcal O_E(m\cdot0_E)$ is globally generated. Its inverse image is $\bigcup_ip_i^{-1}(z)$. Each component is smooth and irreducible, since $p_i$ is a surjective homomorphism with connected kernel. It contains points with pairwise distinct coordinates, where $\alpha_E$ is a local isomorphism, so its multiplicity is one and its class is $p_i^*[\mathrm{pt}]$.

By \cref{lem:holomorphic-cover}, $x_i=\varpi^*p_i^*\theta$ for a fixed nonzero holomorphic one-form $\theta$ on $E$. Since $i\theta\bar\theta$ is positive, summing the classes $\varpi^*p_i^*[\mathrm{pt}]$ gives the displayed formula. Finally, the structure equations give $x_i\bar x_i-x_j\bar x_j=m^{-1}d(y_j-y_i)$.
\end{proof}

\begin{Thm}\label{thm:h-square}
The class $h$ is nonzero, whereas $h^2=0$ in $H^4(Q,\Q)$. More precisely, $[b]\ne0$ in $H^2(Y,\C)$, while, with $w_i=\bar x_i$ and $v_i=\frac12(y_i-\bar y_i)$ as in the proof of \cref{prop:low},
\begin{equation}\label{eq:b-square}
b^2=d\,T_{xwv},\qquad T_{xwv}=\sum_ix_iw_iv_i,
\end{equation}
is exact in the invariant de Rham complex of $Y$.
\end{Thm}

\begin{proof}
The invariant de Rham complex of $Y$ has no nonzero cochains of degree one, since $V^{S_m}=0$. Hence the nonzero closed invariant form $b$ has nonzero class in $H^2(Y,\C)$. Under the identification $H^*(Y,\C)\simeq H^*(M,\C)^{S_m}$ of \cref{lem:analytic-decomposition}, \cref{lem:hyperplane-pullback} gives $h_Y=\kappa[ib]\ne0$. The map $f^*$ is injective on rational cohomology: $f$ has degree one and $Y$ is a rational homology manifold, so the projection formula gives
\[
f_*(f^*\alpha\cap[Q])=\alpha\cap[Y],
\]
and Poincaré duality on $Y$ applies. Thus $h=f^*h_Y\ne0$.

For the square, recall from \cref{prop:low} that $dw_i=0$ and $dv_i=b-mx_iw_i$. Since $x_iw_i$ is closed of even degree and $(x_iw_i)^2=0$,
\[
dT_{xwv}=\sum_ix_iw_i\,dv_i=\sum_ix_iw_i\,(b-mx_iw_i)=b^2.
\]
The primitive $T_{xwv}$ is $S_m$-invariant, so $h_Y^2=-\kappa^2[b]^2=0$, and $h^2=f^*h_Y^2=0$.
\end{proof}

\begin{Rem}\label{rem:geometric-square}
The vanishing has a geometric meaning. For $z\in E$ let $D(z,i)=(p_i\circ\varpi)^{-1}(z)$. By \cref{lem:hyperplane-pullback}, $\varpi^*\alpha_E^*H_z=\sum_iD(z,i)$, and all these divisors have the same rational cohomology class. Since $D(z,1)\cap D(z',1)=\emptyset$ for $z\ne z'$, this class has square zero. Consequently, the preimage of $H_z\cap H_{z'}$ is the nonempty union of the rationally null-homologous submanifolds $D(z,i)\cap D(z',j)$, $i\ne j$. At the level of forms, $[b]=m[x_1\bar x_1]$ and $(x_1\bar x_1)^2=0$. This is a higher-dimensional version of the vanishing of the fiber class of the Kodaira surface, where $x\bar x=d\zeta$; here $[b]$ itself remains nonzero.
\end{Rem}

\begin{Cor}\label{cor:null-homologous}
Let $g\colon Q\to B$ be a holomorphic Lagrangian fibration with normal projective base, and let $A$ be an ample line bundle on $B$. Then
\[
g^*c_1(A)=k h\ne0,\qquad (g^*c_1(A))^2=0
\]
for some positive integer $k$. If $Z\subset B$ is a smooth subvariety of codimension $c\ge2$ transverse to $g$, then $g^{-1}(Z)$ is a nonempty compact complex submanifold which is null-homologous over $\Q$. In particular, this applies to a general fiber and to the inverse image of a general codimension-two linear subspace after identifying $B$ with $\mathbb P^n$.
\end{Cor}

\begin{proof}
By \cref{thm:unique-projective-fibration}, $g=\psi\circ\ell$ for an isomorphism $\psi\colon\mathbb P^n\to B$, and $\psi^*A=\mathcal O(k)$ with $k>0$. The first assertions follow from \cref{thm:h-square}. The class of $\psi^{-1}(Z)$ is $\deg(\psi^{-1}(Z))$ times the $c$-th power of the hyperplane class. Transversality therefore gives
\[
\operatorname{PD}[g^{-1}(Z)]=\deg(\psi^{-1}(Z))\,h^c=0.
\]
The preimage is nonempty by surjectivity. A general point is a regular value by Sard's theorem, and a general linear subspace is transverse by the parametric transversality theorem.
\end{proof}

\begin{Prop}\label{prop:class-C}
Let $X$ be a compact complex manifold of Fujiki class $\mathcal C$ and $\phi\colon X\to\mathbb P^r$ a surjective holomorphic map. Then $(\phi^*c_1(\mathcal O_{\mathbb P^r}(1)))^j\ne0$ in $H^{2j}(X,\R)$ for $0\le j\le r$. In particular, a standard $\BG_n$ manifold with $n\ge2$ is not of class $\mathcal C$.
\end{Prop}

\begin{proof}
Choose a proper bimeromorphic holomorphic map $\mu\colon\hat X\to X$ with $\hat X$ compact Kähler \cite{Varouchas1989,Fujiki1978,Hironaka1975}, and let $\omega$ be a Kähler form on $\hat X$. Put $d=\dim X$, and let $\omega_{FS}$ be the Fubini--Study form. The differential of $\phi\circ\mu$ has rank $r$ on a dense open set, where $(\mu^*\phi^*\omega_{FS})^j\wedge\omega^{d-j}$ is positive, and this form is nonnegative everywhere. Thus
\[
\int_{\hat X}(\mu^*\phi^*\omega_{FS})^j\wedge\omega^{d-j}>0,
\]
which proves the nonvanishing. For the last assertion apply this to $\ell$ with $j=2\le n$ and use \cref{thm:h-square}.
\end{proof}

Non-membership in class $\mathcal C$ also follows from $b_3=3$ in \cref{prop:low}, since manifolds of class $\mathcal C$ satisfy the $\partial\bar\partial$-lemma \cite{DGMS1975}. The argument above detects the obstruction through the square of a single class of degree two. In \cref{sec:bbf} we show that the same square controls the BBF form.

\section{Beauville--Bogomolov--Fujiki form}\label{sec:bbf}

The Beauville--Bogomolov--Fujiki form is the higher-dimensional analogue of the intersection form on a holomorphic symplectic surface. For an irreducible hyperkähler manifold, it determines all top-degree intersections of classes of degree two through the Fujiki relation \cite{Beauville1983,Fujiki1987}. In the present setting the Fujiki relation still holds, and we will see how the vanishing of the square of the Lagrangian class leads to degeneracy of the form.

Let $X$ be a compact complex manifold of dimension $2n$ with a closed holomorphic symplectic form $\tau$, normalized by $\int_X(\tau\bar\tau)^n=1$. The Beauville--Bogomolov--Fujiki quadratic form on $H^2(X,\C)$ is defined by Beauville's formula \cite[Section~8]{Beauville1983}, see also \cite[(1.3)]{CattaneoTomassini2018}:
\begin{equation}\label{eq:bbf-definition}
q_\tau(\alpha)=\frac n2\int_X\alpha^2(\tau\bar\tau)^{n-1}+(1-n)\left(\int_X\alpha\tau^{n-1}\bar\tau^n\right)\left(\int_X\alpha\tau^n\bar\tau^{n-1}\right).
\end{equation}
Here products denote cup products, or wedge products of closed representatives. The formula is well defined on cohomology because $\tau$ is closed, and it is real on $H^2(X,\R)$. In dimension two it is one half of the intersection form. Its definition does not require the Kähler condition.

If $X$ is irreducible hyperkähler, $q_\tau$ is a positive multiple of a primitive integral form $q_X$, and Fujiki's theorem gives a constant $c_X>0$ such that
\[
\int_X\alpha^{2n}=c_Xq_X(\alpha)^n\qquad\text{for every }\alpha\in H^2(X,\R).
\]
Here primitive means that the integral form is not an integer multiple greater than one of another integral form. Expanding this identity for $\alpha=\sum_i t_i\alpha_i$ recovers the mixed top-degree intersections from the associated bilinear form. Neither this relation nor nondegeneracy is automatic for an arbitrary compact complex symplectic manifold.

We compute the BBF form of a standard $\BG_n$-type manifold $Q=Q_{n,t}$, with $n\ge2$ and $m=n+1$. Kurnosov--Verbitsky establish a Fujiki relation $\int_Q\alpha^{2n}=c_Qq(\alpha)^n$, where $q$ is a real quadratic form and $c_Q\ne0$ \cite[Theorem 4.3 and Corollary 5.10]{KurnosovVerbitsky2019}. We use the same letter for its complexification and write $q(\alpha,\beta)=\frac12(q(\alpha+\beta)-q(\alpha)-q(\beta))$ for its associated bilinear form.

The radical is $\operatorname{rad}(q)=\{\beta\in H^2(Q,\R):q(\beta,\alpha)=0\text{ for every }\alpha\in H^2(Q,\R)\}$. A class $\beta$ is isotropic if $q(\beta)=0$; belonging to the radical is the stronger condition that it pairs trivially with every class.

Let $\sigma_Q$ be the symplectic form induced by $\sigma=\sum_i x_i y_i$ on $M$, and put $\eta=[\sigma_Q+\bar\sigma_Q]$. We normalize $q(\eta)=1$, so that $c_Q=\int_Q\eta^{2n}>0$. Inertia will mean the numbers of positive, negative and zero directions, in that order. {The degeneracy of $q$ is forced by the vanishing $h^2=0$ of \cref{thm:h-square}, see \cref{subsec:null}. The remaining inertia is computed in \cref{lem:bbf-U}. Then in \cref{subsec:massey} we show that the annihilator of $h$ in $f^*H^2(Y,\C)$, spanned by the classes of the Massey product of \cite[Section~3]{RiosOrtiz2026BGNonformal}, is a maximal isotropic subspace, and \cref{prop:fiber-restriction} proves that all classes of degree two restrict to zero on the fibers of $\ell$ over $\mathbb P^n\setminus\cD$.}

This normalization agrees with \eqref{eq:bbf-definition} up to positive scale. Indeed, if $\tau$ is a normalized multiple of $\sigma_Q$, polarizing the Fujiki relation gives $q=q_\tau/q_\tau(\eta)$. Thus the form used below is Beauville's form with the normalization $q(\eta)=1$.

Use the coframe of \cref{lem:holomorphic-cover} and put $w_i=\bar x_i$, $u_i=y_i+\bar y_i$ and $v_i=(y_i-\bar y_i)/2$. For $B(r,s)=\sum_i r_i s_i$, set
\[
a=B(x,u),\quad b=B(x,w),\quad c=B(w,u),\quad p=B(x,v),\quad q_0=B(w,v).
\]
These are closed invariant two-forms, and their classes give $H^2(Y,\C)=\langle a,b,c,p,q_0\rangle$, where $Y=M/S_m$, see \cite[equation~(9) and Corollary~5.5]{RiosOrtiz2026BGNonformal}. Indeed, the six invariant pairings between the four copies of the standard representation span invariant degree two. The remaining pairing $B(u,v)$ is not closed for $n\ge2$, and there are no invariant degree-one forms, hence no degree-two boundaries. We use the same symbols for these forms and their classes. {In particular $b=\sum_ix_i\bar x_i$ is the form of \cref{lem:hyperplane-pullback}.}

Let $f:Q\to Y$ be the resolution, let $D$ be its reduced exceptional divisor, and put $e=c_1(\mathcal O_Q(D))$. By \cref{prop:low}, the transposition sector is connected and contributes one class in degree two. On a generic exceptional curve, $e$ has degree $-2$, whereas every pullback from $Y$ has degree zero. Consequently
\[
H^2(Q,\C)=U\oplus\C e,\qquad U=f^*H^2(Y,\C).
\]
{Here $f^*$ is injective, see the proof of \cref{thm:h-square}.}
 We also use $a,b,c,p,q_0$ for their pullbacks to $U$.

\begin{Lem}\label{lem:symplectic-class}
Let $\sigma_Y,\bar\sigma_Y\in H^2(Y,\C)$ be the classes of $\sigma$ and $\bar\sigma$. Then $[\sigma_Q]=f^*\sigma_Y$ and $[\bar\sigma_Q]=f^*\bar\sigma_Y$. In particular $\eta=f^*\eta_Y$, where $\eta_Y\in H^2(Y,\R)$ is the class of $\sigma+\bar\sigma$.
\end{Lem}

\begin{proof}
Put $\Theta=f^{-1}(f(D))$. Over $Y\setminus f(D)$ the map $f$ is an isomorphism and $\sigma_Q=f^*\sigma$, so the class $\delta=[\sigma_Q]-f^*\sigma_Y$ restricts to zero on $Q\setminus\Theta$. By semismallness, $\Theta$ is the union of $D$ and a closed analytic subset of dimension at most $2n-2$. By the exact sequence of the pair $(Q,Q\setminus\Theta)$ and Lefschetz duality $H^2(Q,Q\setminus\Theta;\C)\simeq H_{4n-2}(\Theta,\C)$, the kernel of $H^2(Q,\C)\to H^2(Q\setminus\Theta,\C)$ is spanned by the fundamental classes of the components of $\Theta$ of dimension $2n-1$. The divisor $D$ is irreducible: over the transposition stratum, whose normalization is connected by \cref{prop:low}, the map $f$ is a $\mathbb P^1$-bundle. Hence $\delta=\varepsilon e$ for some $\varepsilon\in\C$. On a generic exceptional curve $\mathbb P^1\subset D$, the form $\sigma_Q$ restricts to zero by type, the class $f^*\sigma_Y$ restricts to zero since $f$ contracts the curve, and $e$ has degree $-2$. Thus $\varepsilon=0$. Conjugation gives the second identity.
\end{proof}

\subsection{Degeneracy of the BBF form}\label{subsec:null}

The degeneracy of $q$ has a simple cause. By \cref{thm:h-square}, the class $h$ has square zero, and a Fujiki relation forces every such class into the radical.

\begin{Lem}\label{lem:fujiki-null}
Let $q'$ be a real quadratic form on $H^2(Q,\R)$ such that $\int_Q\alpha^{2n}=c'q'(\alpha)^n$ for every $\alpha\in H^2(Q,\R)$, where $c'$ is a constant. Then $c'\ne0$, and every $\beta\in H^2(Q,\R)$ with $\beta^2=0$ satisfies $q'(\beta,\alpha)=0$ for all $\alpha\in H^2(Q,\R)$.
\end{Lem}

\begin{proof}
The constant $c'$ is nonzero because $\int_Q\eta^{2n}>0$. Fix $\alpha$. Since $\beta^2=0$, the polynomial $s\mapsto\int_Q(\alpha+s\beta)^{2n}=\int_Q\alpha^{2n}+2ns\int_Q\alpha^{2n-1}\beta$ has degree at most one. It equals $c'(q'(\alpha)+2s\,q'(\alpha,\beta)+s^2q'(\beta))^n$. The coefficient of $s^{2n}$ gives $q'(\beta)=0$, and then the coefficient of $s^n$, with $n\ge2$, gives $q'(\alpha,\beta)=0$.
\end{proof}

\begin{Cor}\label{cor:h-radical}
The class $h$ lies in the radical of $q$, and more generally in the radical of every real quadratic form on $H^2(Q,\R)$ satisfying a Fujiki relation. In particular $q$ is degenerate.
\end{Cor}

\begin{proof}
Combine \cref{thm:h-square} with \cref{lem:fujiki-null}, and recall that $h\ne0$.
\end{proof}

\cref{thm:bbf-radical} below shows that the radical of $q$ is exactly $\R h$.

\begin{Rem}\label{rem:kahler-contrast}
In the Kähler case the BBF form is nondegenerate for a reason that is absent here. Let $X$ be an irreducible hyperkähler manifold of dimension $2n$ with a Kähler class $\omega$, and let $c_X$ be its Fujiki constant. The Fujiki relation and the Hodge--Riemann bilinear relations for $\omega$ show that $q$ is positive definite on $\langle\operatorname{Re}\sigma,\operatorname{Im}\sigma,\omega\rangle$ and negative definite on the $q$-orthogonal complement of $\omega$ in $H^{1,1}(X,\R)$ \cite{Beauville1983,Fujiki1987}. If $\pi\colon X\to\mathbb P^n$ is a Lagrangian fibration and $h_X=\pi^*c_1(\mathcal O_{\mathbb P^n}(1))$, then $h_X$ is isotropic, because $h_X^{n+1}=0$, but it is not in the radical: polarizing the Fujiki relation gives
\[
\int_Xh_X^n\omega^n=\frac{c_X\,2^n(n!)^2}{(2n)!}\,q(h_X,\omega)^n,
\]
and the left side is a positive multiple of the volume of a fiber. In particular $h_X^n\ne0$, in accordance with Verbitsky's theorem that $\Sym^jH^2(X,\Q)\to H^{2j}(X,\Q)$ is injective for $j\le n$ \cite{Verbitsky1996,Bogomolov1996GAFA}.

On a standard $\BG_n$ manifold there is no Kähler class, the class $h$ lies in the radical of $q$ rather than being merely isotropic, and already $h^2=0$.
\end{Rem}

\begin{Lem}\label{lem:bbf-U}
On $U$ one has $q(Aa+Bb+Cc+Pp+Rq_0)=PC-AR$. Its complex radical is $\C b$, and its restriction to $U\cap H^2(Q,\R)$ has inertia $(2,2,1)$.
\end{Lem}
\begin{proof}
Choose a real orthonormal basis of the standard representation for the symmetric pairing $\sum_i r_i s_i$. In the corresponding coframe, the two-form $\alpha=Aa+Bb+Cc+Pp+Rq_0$ is the sum of $n$ blocks
\[
\beta=A xu+B xw+C wu+P xv+R wv.
\]
All products are exterior products in the displayed order. For one block, $\beta^2=2(PC-AR)xwuv$. Thus the top exterior power of $\alpha$ is proportional to $(PC-AR)^n$, with the same nonzero volume factor for every $\alpha$. {The coefficient $B$ of $b$ does not occur, in accordance with \cref{cor:h-radical}.}

The quotient map $M\to Y$ has degree $m!$, and $f$ has degree one. Hence for a class $\alpha_Y\in H^2(Y,\C)$,
\[
\int_Q(f^*\alpha_Y)^{2n} =\int_Y\alpha_Y^{2n} =\frac1{m!}\int_M(\pi^*\alpha_Y)^{2n},
\]
where $\pi:M\to Y$ is the quotient map. These integrals use the complex fundamental classes. Now $\sigma=\frac12a+p$ and $\bar\sigma=\frac12c-q_0$, so $PC-AR$ takes the value one on $\eta${, by \cref{lem:symplectic-class}}. Taking the ratio of the integrals for $\alpha$ and $\eta$ and using the Fujiki relation gives $q(\alpha)^n=(PC-AR)^n$. The two quadratic polynomials differ by a constant $n$-th root of unity, and their value one at $\eta$ fixes that constant. This proves the formula.

Complex conjugation acts by $\bar a=c$, $\bar b=-b$ and $\bar p=-q_0$. A real class therefore has $C=\bar A$, $R=-\bar P$ and $B\in i\R$. The quadratic form on this real subspace is
\[
2\operatorname{Re}(P\bar A)
=\frac12|A+P|^2-\frac12|A-P|^2,
\]
with the additional radical direction $\R(ib)$. This gives inertia $(2,2,1)$ and radical $\C b$ after complexification.
\end{proof}

Let $Z$ be the normalization of the closure of the transposition stratum and let $\nu:Z\to Y$ be the natural map. The proof of \cref{prop:low} identifies $Z=M^\tau/S_{m-2}$, where $\tau=(12)$ and $M^\tau$ is connected. This description of the normalization agrees with \cite[Lemma~5.7 and Corollary~5.11]{RiosOrtiz2026BGNonformal}.

\begin{Lem}\label{lem:bbf-exceptional}
The exceptional class satisfies $q(U,e)=0$ and $q(e)<0$.
\end{Lem}
\begin{proof}
The generic transverse resolution over $Z$ is the minimal resolution of an $A_1$ singularity. Its exceptional curve has self-intersection $-2$. The corresponding pushforward identities are
\[
f_*(e\cap[Q])=0,\qquad f_*(e^2\cap[Q])=-2\nu_*[Z].
\]
The first vanishes because $D$ maps to a subset of complex codimension two. For the second, the generic transverse self-intersection determines the coefficient of $[Z]$, and the remaining strata have too small a dimension to support a class in $H_{4n-4}(Y,\Q)$ coming from this intersection.

The projection formula gives $\int_Q e\alpha^{2n-1}=0$ for $\alpha\in U$. Comparing coefficients of $\lambda$ in the Fujiki relation for $\alpha+\lambda e$ yields $q(\alpha)^{n-1}q(\alpha,e)=0$. Since $q|_U$ is not identically zero, the linear function $q(\alpha,e)$ vanishes on a dense open subset of $U$, and hence on all of $U$.

Let $\eta_Y\in H^2(Y,\R)$ be the class descending from $\sigma+\bar\sigma$, so that $f^*\eta_Y=\eta$ {by \cref{lem:symplectic-class}}. The fixed tangent space of a finite-order symplectic automorphism is symplectic: averaging over the automorphism shows that an invariant vector orthogonal to all invariant vectors is orthogonal to every vector. Thus $\sigma|_{M^\tau}$ is nondegenerate. The effective quotient $M^\tau\to Z$ has generic degree $(m-2)!$, so
\[
J:=\int_Z(\nu^*\eta_Y)^{2n-2} =\frac1{(m-2)!}\int_{M^\tau}  (\sigma|_{M^\tau}+\bar\sigma|_{M^\tau})^{2n-2}>0.
\]
Positivity follows from the complex orientation: on a complex symplectic manifold of dimension $2d$, the top power of $\sigma+\bar\sigma$ is $\binom{2d}{d}\sigma^d\bar\sigma^d$, a positive volume form.

The second pushforward identity now gives $\int_Q e^2\eta^{2n-2}=-2J$. Since $q(\eta)=1$ and $q(\eta,e)=0$, the coefficient of $\lambda^2$ in the Fujiki relation is $\binom{2n}{2}(-2J)=c_Qnq(e)$. Therefore
\[
q(e)=-\frac{2(2n-1)J}{\int_Q\eta^{2n}}<0.
\]
\end{proof}

\begin{Thm}\label{thm:bbf-radical}
For a standard $\BG_n$-type manifold $Q$, $n\ge2$, the BBF form has real inertia $(2,3,1)$. If $\ell:Q\to\mathbb P^n$ is its Lagrangian fibration {(\cref{def:lagrangian})} and $h=\ell^*c_1(\mathcal O_{\mathbb P^n}(1))$, then $\operatorname{rad}(q)=\R h$. For any holomorphic Lagrangian fibration $g\colon Q\to B$ with normal projective base and any ample line bundle $A$ on $B$, this radical is also $\R g^*c_1(A)$. Every $\BG_n$-type manifold carries the Fujiki form transported from a standard representative, with the same inertia and one-dimensional radical.
\end{Thm}
\begin{proof}
\cref{lem:bbf-U,lem:bbf-exceptional} give the orthogonal decomposition $H^2(Q,\C)=U\oplus\C e$. The exceptional class adds one negative direction to the real inertia $(2,2,1)$ on $U$. The radical is therefore $\R(ib)$ and the full inertia is $(2,3,1)$.

{By \cref{lem:hyperplane-pullback} and the proof of \cref{thm:h-square}, $h=\kappa f^*[ib]$ with $\kappa>0$. Hence $\operatorname{rad}(q)=\R h$, in accordance with \cref{cor:h-radical}.}

The assertion for arbitrary fibrations on $Q$ follows from \cref{cor:null-homologous}.

In a smooth proper deformation family, parallel transport preserves the cup product and complex orientation. It therefore transports the Fujiki relation and the real quadratic form, preserving its inertia and radical dimension. This applies along any chain of deformation families and proves the last assertion. There is no ambiguity of sign under monodromy: the top intersection polynomial determines $q$ up to scalar, and $q$ and $-q$ have different inertias. Along a family of holomorphic symplectic manifolds, this transported sign agrees with the normalization used above. Indeed, for a continuous local choice of $\sigma_s$, the integral $\int(\sigma_s+\bar\sigma_s)^{2n}$ stays positive, so $q([\sigma_s+\bar\sigma_s])$ cannot cross zero.
\end{proof}

The same cohomological behavior holds for arbitrary Lagrangian fibrations on sufficiently small deformations, whenever such fibrations exist.

\begin{Prop}\label{prop:deformed-lagrangian-class}
Let $X$ be a sufficiently small deformation of a standard $\BG_n$-type manifold $Q_{n,t}$, with $n\ge2$, and let $q_X$ be the transported Fujiki form. Suppose that $g\colon X\to B$ is a holomorphic Lagrangian fibration onto a normal projective variety, and let $A$ be an ample line bundle on $B$. Then $\alpha=g^*c_1(A)$ is nonzero and
\[
\operatorname{rad}(q_X)=\R\alpha,\qquad \alpha^2=0.
\]
A general fiber of $g$ is null-homologous over $\Q$. Moreover, any two holomorphic Lagrangian fibrations on $X$ with normal projective bases differ by an isomorphism of their bases.
\end{Prop}

\begin{proof}
By \cref{cor:small-deformations,cor:unobstructed-deformations}, $X$ is holomorphic symplectic and its Frölicher spectral sequence degenerates at $E_1$. Since $b_1(X)=0$, one has $H^1(X,\mathcal O_X)=0$. The exponential sequence therefore makes the first Chern class map $\operatorname{Pic}(X)\to H^2(X,\Z)$ injective. Put $L=g^*A$. If $c_1(L)$ vanished in real cohomology, some positive power of $L$ would be trivial. This is impossible: a sufficiently high power of $A$ is very ample, and its sections pull back to independent sections of the corresponding power of $L$. Thus $\alpha\ne0$.

Let $\sigma_X$ be a holomorphic symplectic form chosen continuously from the standard fiber. The real plane
\[
P_X=\langle\operatorname{Re}[\sigma_X],\operatorname{Im}[\sigma_X]\rangle
\]
is positive definite for $q_X$ after shrinking the deformation space. By \cref{thm:bbf-radical}, $q_X$ has inertia $(2,3,1)$, so its restriction to $P_X^\perp$ is negative semidefinite, with kernel $\operatorname{rad}(q_X)$.

The class $\alpha$ is represented by a closed real $(1,1)$-form. Type considerations give
\[
\int_X(\alpha+s\sigma_X)^{2n}=\int_X\alpha^{2n}
\]
for every $s\in\C$. Applying the complexified Fujiki relation shows that $q_X(\alpha,\sigma_X)=0$, and similarly $q_X(\alpha,\bar\sigma_X)=0$. Hence $\alpha\in P_X^\perp$. Since $\dim B=n$, one has $\alpha^{2n}=0$, so $q_X(\alpha)=0$. Negative semidefiniteness now implies $\alpha\in\operatorname{rad}(q_X)$, and nonvanishing shows that it spans the radical.

Let $\rho_X$ be the parallel transport of the natural class $h$ on $Q_{n,t}$. The radical is $\R\rho_X$, and $\rho_X^2=0$ because parallel transport preserves cup products and $h^2=0$ by \cref{thm:h-square}. Thus $\alpha^2=0$. Intersecting $n$ general divisors from a sufficiently high power of $A$ gives a finite set of regular values of $g$. Its inverse image is a union of smooth fibers whose class is a positive multiple of $\alpha^n=0$. These fibers are homologous, since they lie over the connected open set of regular values. Therefore the class of a general fiber is zero over $\Q$.

For uniqueness, let $g_i\colon X\to B_i$, $i=1,2$, be two such fibrations, choose ample line bundles $A_i$, and put $L_i=g_i^*A_i$. Their nonzero integral Chern classes lie on the same real line, so they are rationally proportional. They cannot lie on opposite rays. Indeed, in that case some positive tensor powers of $L_1$ and $L_2$ would have trivial tensor product, by the injectivity of the Chern class map after killing torsion. Choosing sections of sufficiently high powers with nonempty zero divisors would then give a nonzero holomorphic function on $X$ with a zero, a contradiction.

Consequently $L_1^a\simeq L_2^b$ for some positive integers $a,b$. After taking a common further power, $A_1^a$ and $A_2^b$ are very ample. Since each $g_i$ has connected fibers and normal base, $(g_i)_*\mathcal O_X=\mathcal O_{B_i}$, so the complete linear systems of $L_1^a$ and $L_2^b$ are pulled back from those of $A_1^a$ and $A_2^b$. The common line bundle therefore defines both fibrations, and identifies their bases.
\end{proof}

\begin{Rem}\label{rem:kodaira-form}
On $U/\C h$ the form $q$ induces a positive multiple of the intersection form of the Kodaira surface $S$. Put $\tilde\sigma=\sum_ix_i\bar y_i$, the fiberwise conjugate of $\sigma$. Then
\[
a=\sigma+\tilde\sigma,\qquad c=\bar\sigma+\bar{\tilde\sigma},\qquad p=\tfrac12(\sigma-\tilde\sigma),\qquad q_0=-\tfrac12(\bar\sigma-\bar{\tilde\sigma}),
\]
and the formula $PC-AR$ of \cref{lem:bbf-U} becomes
\begin{equation}\label{eq:q-kodaira}
q(\lambda_1\sigma+\lambda_2\tilde\sigma+\mu_1\bar\sigma+\mu_2\bar{\tilde\sigma}+Bb)=\lambda_1\mu_1-\lambda_2\mu_2.
\end{equation}
On $S$, with the coframe $x,\zeta$ of the proof of \cref{lem:holomorphic-cover}, the classes of $\sigma_S=x\zeta$, $\tilde\sigma_S=x\bar\zeta$ and their conjugates form a basis of $H^2(S,\C)$, while $x\bar x=d\zeta$ is exact. Pointwise $\tilde\sigma_S\bar{\tilde\sigma}_S=-\sigma_S\bar\sigma_S$, and all other products of two of these four forms vanish, so that $\int_S\alpha_S^2=2(\lambda_1\mu_1-\lambda_2\mu_2)\int_S\sigma_S\bar\sigma_S$ for $\alpha_S=\lambda_1\sigma_S+\lambda_2\tilde\sigma_S+\mu_1\bar\sigma_S+\mu_2\bar{\tilde\sigma}_S$.

Symmetrization $\alpha_S\mapsto\sum_i(\alpha_S)_i$, where $(\alpha_S)_i$ is the copy of $\alpha_S$ on the $i$-th factor of $S^m$, takes closed invariant forms on $S$ to closed forms on $R$. These are basic for the diagonal action of $F$, because $\sum_ix_i=0$, so they descend to $N_n$ and define classes on $Y$. Since $y_i=\sum_j\zeta_j-m\zeta_i$, symmetrization sends $\sigma_S$ and $\tilde\sigma_S$ to $-\frac1m\sigma$ and $-\frac1m\tilde\sigma$, and the exact form $x\bar x$ to $b$. It therefore induces an isomorphism $\operatorname{sym}\colon H^2(S,\C)\to H^2(Y,\C)/\C b$, and \eqref{eq:q-kodaira} becomes
\[
q\bigl(f^*\operatorname{sym}(\alpha_S)\bigr)=\frac{1}{2m^2\int_S\sigma_S\bar\sigma_S}\int_S\alpha_S^2,\qquad\alpha_S\in H^2(S,\C),
\]
where the left side is well defined because $f^*b$ lies in the radical. Together with \cref{lem:bbf-exceptional}, this is the analogue of Beauville's description of the BBF form of the Hilbert scheme of points on a K3 surface, which restricts to the intersection form on the second cohomology of the surface and is completed by the class of the exceptional divisor \cite{Beauville1983}. The new feature is the radical, spanned by the symmetrization $b$ of the exact form $x\bar x=d\zeta$ (\cref{rem:geometric-square}).
\end{Rem}

\begin{Rem}\label{rem:intro-guan-bbf}
Guan proved the existence of a one-dimensional kernel in \cite{GuanIII1995}. This result is recalled explicitly in \cite[Theorem~6.7]{Guan2023Survey}. The computation above recovers this degeneracy and describes the real inertia and the radical in the quotient model. In particular, the nondegeneracy asserted in \cite[Conjecture~5.11]{KurnosovVerbitsky2019} does not hold for these examples. Their Fujiki relation, used here, does not require nondegeneracy. {The nondegeneracy was also expected in \cite[Section~1.1]{BKKY2022}.}
\end{Rem}

\begin{Rem}
The one-dimensional radical also obstructs the $\partial\bar\partial$-lemma. Indeed, for a compact complex symplectic manifold $X$ satisfying that lemma, \cite[Proposition~2.10]{CattaneoTomassini2018} gives $\dim_{\C}\operatorname{rad}(q_\sigma)=2h^{2,0}(X)-2$. In particular, the radical has even dimension. {Here $q_\sigma$ is defined by Beauville's formula \cite[(1.3)]{CattaneoTomassini2018}. Polarizing a Fujiki relation shows that $q_\sigma$ is a nonzero multiple of every form satisfying it, so the two radicals coincide.} \cref{thm:bbf-radical} therefore gives another proof that a standard $\BG_n$-type manifold does not satisfy the $\partial\bar\partial$-lemma. This obstruction concerns its complex symplectic structure, whereas the nonformality theorem in \cite[Theorem~1.2 and Corollary~5.2]{RiosOrtiz2026BGNonformal} excludes every complex structure satisfying the $\partial\bar\partial$-lemma on the same underlying smooth manifold.
\end{Rem}

\subsection{Failure of hard Lefschetz}\label{subsec:lefschetz}

The radical of the BBF form also obstructs the hard Lefschetz property.  One of the distinctive features of hyperkähler manifolds is the action of their Looijenga--Lunts--Verbitsky (LLV) Lie algebra on cohomology \cite{Bogomolov1996GAFA,Verbitsky1996,LooijengaLunts1997}. This construction combines the Lefschetz $\mathfrak{sl}_2$-actions associated with classes of degree two. More precisely, for a compact complex manifold $X$ of dimension $2n$, put $L_\omega(\alpha)=\omega\cup\alpha$ and let $H$ act on $H^k(X,\C)$ as multiplication by $k-2n$. A class $\omega\in H^2(X,\C)$ satisfies hard Lefschetz if
\[
L_\omega^{2n-k}\colon H^k(X,\C)\longrightarrow H^{4n-k}(X,\C)
\]
is an isomorphism for every $0\le k\le2n$. This is equivalent to the existence of an operator $\Lambda_\omega$ of degree $-2$ such that $(L_\omega,H,\Lambda_\omega)$ is an $\mathfrak{sl}_2$-triple. The total LLV Lie algebra is generated by these triples as $\omega$ ranges over the classes satisfying hard Lefschetz. Although Kähler classes provide such triples by the hard Lefschetz theorem, the construction itself depends only on the graded cohomology ring.

For $\BG_n$-type manifolds, the odd Betti number $b_3=3$ already excludes hard Lefschetz. The radical of the BBF form gives a further obstruction in degree two, so the failure persists even after restricting to even-degree cohomology.

\begin{Cor}\label{cor:no-LLV}
Let $X$ be a $\BG_n$-type manifold, with $n\geq2$, and let $q_X$ be its transported Fujiki form. For every $\omega\in H^2(X,\C)$,
\[
\operatorname{rank}\bigl(L_\omega^{2n-3}\colon H^3(X,\C)\to H^{4n-3}(X,\C)\bigr)\le2
\]
and
\[
\operatorname{rad}(q_X)\subseteq\ker\bigl(L_\omega^{2n-2}\colon H^2(X,\C)\to H^{4n-2}(X,\C)\bigr).
\]
If $q_X(\omega)\ne0$, the latter kernel is exactly $\operatorname{rad}(q_X)$, and the map has rank five. In particular, no class satisfies hard Lefschetz, even on $H^{\mathrm{even}}(X,\C)$, and no cup-product operator $L_\omega$ belongs to a Lefschetz $\mathfrak{sl}_2$-triple with grading operator $H$.
\end{Cor}
\begin{proof}
By Poincaré duality, the map $L_\omega^{2n-3}$ has the same rank as the bilinear form
\[
B_\omega(\alpha,\beta)=\int_X\alpha\cup\beta\cup\omega^{2n-3},
\qquad \alpha,\beta\in H^3(X,\C).
\]
Graded commutativity gives $B_\omega(\beta,\alpha)=-B_\omega(\alpha,\beta)$. Since $b_3(X)=3$ by \cref{prop:low}, this alternating form has rank at most two and hard Lefschetz fails in degree three.

For degree two, write the Fujiki relation as $\int_X\gamma^{2n}=c_Xq_X(\gamma)^n$. Differentiating twice at $\gamma=\omega$ in the directions $\alpha,\beta\in H^2(X,\C)$ gives
\[
\int_X\alpha\beta\omega^{2n-2}
=\frac{c_X}{2n-1}\left(q_X(\omega)^{n-1}q_X(\alpha,\beta)+2(n-1)q_X(\omega)^{n-2}q_X(\alpha,\omega)q_X(\beta,\omega)\right).
\]
If $\alpha\in\operatorname{rad}(q_X)$, the right-hand side vanishes for every $\beta$. Poincaré duality therefore gives $\alpha\omega^{2n-2}=0$, proving the asserted inclusion.

Suppose now that $q_X(\omega)\ne0$ and $\alpha\omega^{2n-2}=0$. Taking $\beta=\omega$ in the displayed identity gives
\[
0=c_Xq_X(\omega)^{n-1}q_X(\alpha,\omega),
\]
hence $q_X(\alpha,\omega)=0$. Substituting this back into the identity shows that $q_X(\alpha,\beta)=0$ for every $\beta$, so $\alpha\in\operatorname{rad}(q_X)$. By \cref{thm:bbf-radical}, this radical has dimension one, while $b_2(X)=6$, proving the rank assertion. Since the radical is nonzero, hard Lefschetz fails in degree two for every $\omega$. The assertion about Lefschetz triples follows from the equivalence recalled above.
\end{proof}

Consequently, the usual LLV construction has no Lefschetz triples with the cohomological grading from which to generate a Lie algebra. In particular, no symplectic form on the underlying smooth manifold satisfies hard Lefschetz.

\subsection{Relations with Massey products}\label{subsec:massey}

We now determine the products of two classes of degree two, following \cite[Section~3 and Corollary~5.5]{RiosOrtiz2026BGNonformal}. In terms of the families $x,w,u,v$ introduced at the beginning of this section, consider the invariant three-forms
\[
T_{xwu}=\sum_ix_iw_iu_i,\qquad T_{xuv}=\sum_ix_iu_iv_i,\qquad T_{xwv}=\sum_ix_iw_iv_i,\qquad T_{wuv}=\sum_iw_iu_iv_i.
\]
In \cite[equation~(5)]{RiosOrtiz2026BGNonformal} the last three are denoted $\alpha$, $\beta$ and $\gamma$.

\begin{Prop}\label{prop:degree-four}
The invariant three-forms are spanned by $T_{xwu}$, $T_{xuv}$, $T_{xwv}$ and $T_{wuv}$, and
\[
dT_{xwu}=0,\qquad dT_{xuv}=ab,\qquad dT_{xwv}=b^2,\qquad dT_{wuv}=bc.
\]
In $H^4(Y,\C)$, the fifteen products of two of the classes $a,b,c,p,q_0$ span a subspace of dimension $12$, and the linear relations among them are spanned by $ab=0$, $b^2=0$ and $bc=0$. The annihilator of $b$ in $H^2(Y,\C)$ is $\langle a,b,c\rangle$. The same statements hold after applying the injective ring homomorphism $f^*$.
\end{Prop}

\begin{proof}
The description of the invariant three-forms and the identities for their differentials are given in \cite[equation~(6) and the proof of Lemma~3.1]{RiosOrtiz2026BGNonformal}; we recall the argument. The differentials are computed as in the proof of \cref{thm:h-square}, using $dx_i=dw_i=du_i=0$ and $dv_i=b-mx_iw_i$. An invariant three-form has one factor from each of three different families. Indeed, the representations $\Lambda^kV$, $1\le k\le m-1$, are pairwise nonisomorphic irreducibles \cite{FultonHarris1991}, so $(\Lambda^3V)^{S_m}=(\Lambda^2V\otimes V)^{S_m}=0$, while $(V^{\otimes3})^{S_m}$ is spanned by the symmetric tensor $\sum_ir_is_it_i$. Thus the exact invariant four-forms are spanned by $ab$, $b^2$ and $bc$.

The invariant algebra is graded by $\Z^4$, according to the number of factors from the families $x,w,u,v$, and $d$ is homogeneous of degree $(1,1,0,-1)$. Hence a form is exact if and only if each of its multihomogeneous components is exact. The products $ab$, $b^2$ and $bc$ lie in the multidegrees $(2,1,1,0)$, $(2,2,0,0)$ and $(1,2,1,0)$, which contain no other product. Every other multidegree contains a single product, except $(1,1,1,1)$, which contains $aq_0$ and $cp$. Choose orthonormal coordinates $X_k,W_k,U_k,V_k$, $1\le k\le n$, in the four families, as in the proof of \cref{lem:bbf-U}. Then $a=\sum_kX_kU_k$, $b=\sum_kX_kW_k$, $c=\sum_kW_kU_k$, $p=\sum_kX_kV_k$ and $q_0=\sum_kW_kV_k$. A product of two of these classes is $\sum_{k,l}m_{kl}$, where $m_{kl}$ is a product of four coordinates. It vanishes if a coordinate repeats, which can only happen when $k=l$. For two different classes, the nonzero $m_{kl}$ are pairwise distinct monomials. For a square, $m_{kl}=m_{lk}$, and the $m_{kl}$ with $k<l$ are pairwise distinct. As $n\ge2$, every product is a nonzero form. Finally $X_1U_1W_2V_2$ occurs in $aq_0$ but not in $cp=\sum_{k,l}W_kU_kX_lV_l$, so $aq_0$ and $cp$ are independent. Hence the fifteen products are linearly independent forms, of which exactly $ab$, $b^2$ and $bc$ are exact. If $b(\alpha_1a+\alpha_2b+\alpha_3c+\alpha_4p+\alpha_5q_0)$ is exact, then so is $\alpha_4bp+\alpha_5bq_0$, which forces $\alpha_4=\alpha_5=0$. The map $f^*$ is injective by the proof of \cref{thm:h-square}, and it is multiplicative.
\end{proof}

In particular $\Sym^2H^2(Q,\C)\to H^4(Q,\C)$ is not injective, in contrast with Verbitsky's theorem recalled in \cref{rem:kahler-contrast}.

The three relations $ab=dT_{xuv}$, $b^2=dT_{xwv}$ and $bc=dT_{wuv}$ say that $b$ annihilates $a$, $b$ and $c$. In this situation Cavalcanti, Fernández and Muñoz define a Massey product $\langle b;a,b,c\rangle\subset H^8(Y,\C)$ with central class $b$: it is the set of classes of the forms $\mathfrak M$ below, for all choices of the primitives \cite[Definition~2.3]{CFM2008}. A nontrivial Massey product of this kind obstructs formality \cite[Theorem~2.10]{CFM2008}. We compute it in the invariant model, which is a model of $Y$ over $\C$ \cite[Corollary~5.5]{RiosOrtiz2026BGNonformal}, so that this criterion applies to it.

\begin{Prop}\label{prop:massey}
\begin{enumerate}
\item In the invariant model, the Massey product $\langle b;a,b,c\rangle$ consists of the single class of the closed form
\[
\mathfrak M=a\,T_{xwv}T_{wuv}+b\,T_{wuv}T_{xuv}+c\,T_{xuv}T_{xwv}.
\]
\item Its central class is the hyperplane class: $f^*[ib]=\kappa^{-1}h$, where $\kappa>0$ is the constant of \cref{lem:hyperplane-pullback}.
\end{enumerate}
\end{Prop}

\begin{proof}
(1) By \cref{prop:degree-four}, and since forms of even degree commute with all forms,
\[
d\mathfrak M=a\,(b^2T_{wuv}-T_{xwv}\,bc)+b\,(bc\,T_{xuv}-T_{wuv}\,ab)+c\,(ab\,T_{xwv}-T_{xuv}\,b^2)=0.
\]
With $\xi_1=T_{xuv}$, $\xi_2=T_{xwv}$ and $\xi_3=T_{wuv}$, so that $ba=d\xi_1$, $bb=d\xi_2$ and $bc=d\xi_3$, the form of \cite[Definition~2.3]{CFM2008} is $a\xi_2\xi_3-\xi_1b\xi_3+\xi_1\xi_2c=\mathfrak M$. The representatives of $a$, $b$ and $c$ are unique, since there are no invariant one-forms, and two choices of primitives differ by closed invariant three-forms. By \cref{prop:degree-four} these are the multiples of $T_{xwu}$, which is exact: $d\sum_iu_iv_i=-\sum_iu_i(b-mx_iw_i)=mT_{xwu}$. Replace the primitives one at a time. If $\xi_3$ is replaced by $\xi_3+d\theta$, then $\mathfrak M$ changes by
\[
a\xi_2\,d\theta+b\,d\theta\,\xi_1=d\bigl(\theta\,b\,\xi_1-a\,\xi_2\,\theta\bigr),
\]
because $d\xi_1=ab$ and $d\xi_2=b^2$. The other two cases are analogous. Hence $\langle b;a,b,c\rangle=\{[\mathfrak M]\}$, as in \cite[Corollary~3.2]{RiosOrtiz2026BGNonformal}.

(2) By the proof of \cref{thm:h-square}, $h_Y=\kappa[ib]$, so $h=f^*h_Y=\kappa f^*[ib]$.
\end{proof}

Thus the three exactness relations $ab=dT_{xuv}$, $b^2=dT_{xwv}$ and $bc=dT_{wuv}$ have a direct geometric interpretation. Since $b$ represents the Lagrangian class $h$ up to a nonzero scalar, the middle relation is precisely the vanishing $h^2=0$ of \cref{thm:h-square}, while the other two show that $h$ also annihilates the classes $a=\sigma+\tilde\sigma$ and $c=\bar\sigma+\bar{\tilde\sigma}$, in the notation of \cref{rem:kodaira-form}. These three relations define the Massey product with central class $b$: taking $T_{xuv}$, $T_{xwv}$ and $T_{wuv}$ as primitives of $ba$, $b^2$ and $bc$, respectively, gives the closed representative $\mathfrak M=aT_{xwv}T_{wuv}+bT_{wuv}T_{xuv}+cT_{xuv}T_{xwv}$. Its cohomology class is nonzero. Indeed, $\mathfrak M\,a^{n-2}q_0^{n-2}$ is a nonzero multiple of a volume form, whereas no top-degree form is exact in the invariant model \cite[Proposition~4.1 and Theorem~4.5]{RiosOrtiz2026BGNonformal}. Thus the quotient $Y$ is nonformal. The passage from this quotient-model obstruction to nonformality of every $\BG_n$-type manifold is established by the nonzero-degree argument in the proof of \cite[Theorem~1.2]{RiosOrtiz2026BGNonformal}. In complex dimension four this recovers Guan's relative Massey product description for the $\BG_2$ manifold: the distinguished central class is precisely the generator of the radical of the BBF form, which by \cref{thm:bbf-radical} is a multiple of $h$. Thus the same Lagrangian class governs both the degeneracy of the BBF form and the nonformality obstruction in every dimension.

\begin{Cor}\label{cor:isotropic}
The annihilator of $h$ in $U$ for the cup product is
\[
f^*\langle a,b,c\rangle=\bigl\langle\,[\sigma_Q]+f^*[\tilde\sigma],\ [\bar\sigma_Q]+f^*[\bar{\tilde\sigma}],\ h\,\bigr\rangle,
\]
and it is a maximal totally isotropic subspace of $H^2(Q,\C)$ for $q$. On the four classes $a,c,p,q_0$ the only nonzero values of the bilinear form are $q(a,q_0)=-\frac12$ and $q(c,p)=\frac12$. Thus $\langle a,c\rangle$ and $\langle p,q_0\rangle$ span complementary isotropic planes of $U/\C b$.
\end{Cor}
\begin{proof}
By \cref{prop:massey}, $h$ is a nonzero multiple of $f^*b$. Since $f^*$ is an injective ring homomorphism, \cref{prop:degree-four} gives
\[
\{\alpha\in U:h\alpha=0\}=f^*\langle a,b,c\rangle.
\]
Moreover, \cref{rem:kodaira-form} gives $a=\sigma+\tilde\sigma$ and $c=\bar\sigma+\bar{\tilde\sigma}$, while \cref{lem:symplectic-class} identifies $f^*[\sigma]$ with $[\sigma_Q]$. Applying $f^*$ to these identities and replacing $f^*b$ by its nonzero multiple $h$ gives the displayed description of the annihilator.

As before, we use $a,b,c,p,q_0$ also for their pullbacks to $U$. By \cref{lem:bbf-U},
\[
q(Aa+Bb+Cc+Pp+Rq_0)=PC-AR.
\]
The coefficient of a product of two coordinates in a quadratic form is twice the pairing of the corresponding basis vectors. Thus the only nonzero pairings among $a,c,p,q_0$, up to symmetry, are
\[
q(a,q_0)=-\frac12,\qquad q(c,p)=\frac12.
\]
All four classes have square zero for $q$, and $q(a,c)=q(p,q_0)=0$. Consequently $q$ vanishes identically on each of the planes $\langle a,c\rangle$ and $\langle p,q_0\rangle$. Their images are complementary in $U/\C b$, since the images of $a,c,p,q_0$ form a basis of this quotient.

Let $L=f^*\langle a,b,c\rangle$. The preceding calculation shows that every pairing between $a$ and $c$ vanishes. Since $\C b=\C h$ is contained in the radical, every pairing between two elements of $L$ therefore vanishes. Thus $L$ is a three-dimensional totally isotropic subspace.

It remains to prove maximality in $H^2(Q,\C)$. By \cref{thm:bbf-radical}, the radical is exactly $\C h$, so $q$ induces a nondegenerate bilinear form on the five-dimensional space
\[
V=H^2(Q,\C)/\C h.
\]
If $W\subset V$ is totally isotropic, then $W\subset W^\perp$. Nondegeneracy gives $\dim W^\perp=5-\dim W$, hence $2\dim W\le5$ and therefore $\dim W\le2$. For any totally isotropic subspace $L'\subset H^2(Q,\C)$, its image $\overline{L'}$ in $V$ is totally isotropic, and
\[
\dim L'=\dim\overline{L'}+\dim(L'\cap\C h)\le2+1=3.
\]
Since $L$ has dimension three, it is maximal totally isotropic.
\end{proof}

For a Lagrangian fibration $\pi\colon X\to B$ on an irreducible hyperkähler manifold and a general smooth fiber $F$, Matsushita showed that the restriction map $H^2(X,\Q)\to H^2(F,\Q)$ has rank one. Its image is generated by a polarization class of $F$, since the restriction of a Kähler class on $X$ is nonzero and lies in this one-dimensional image after extension of scalars to $\R$. For a standard Bogomolov--Guan manifold, \cref{prop:fiber-restriction} below shows that the corresponding restriction map has rank zero. Thus, although the general fiber is an abelian variety, none of its polarization classes is the restriction of a class in $H^2(Q,\Q)$. In particular, no line bundle on $Q$ restricts to an ample line bundle on the general fiber.

\begin{Prop}\label{prop:fiber-restriction}
Let $g\colon Q\to B$ be a holomorphic Lagrangian fibration with normal projective base. For a general fiber $F_g$, the restriction map $H^2(Q,\C)\to H^2(F_g,\C)$ is zero.
\end{Prop}

\begin{proof}
By \cref{thm:unique-projective-fibration}, it suffices to consider $g=\ell$ and $y\in\mathbb P^n\setminus\cD$. Since $H^2(Q,\C)=U\oplus\C e$, it suffices to consider $e$ and the classes $a,b,c,p,q_0$. By \cref{prop:lagrangian}(1), the fiber $\ell^{-1}(y)$ does not meet $D$, so $e$ restricts to zero, and $f$ and $\pi$ identify $\ell^{-1}(y)$ with a fiber $\varpi^{-1}(\mathbf z)$ of $M$, where $\mathbf z\in\alpha_E^{-1}(y)$. Every term of the forms $a,b,c,p,q_0$ contains a factor $x_i$ or $\bar x_i$, and such factors vanish on the fibers of $\varpi$.
\end{proof}

The restriction statement also allows the uniqueness theorem to be applied to smooth Kähler bases.

\begin{Cor}\label{cor:unique-kahler-base}
Let $Q=Q_{n,t}$, with $n\ge2$, and let $g\colon Q\to B$ be a nonconstant surjective holomorphic map with connected fibers onto a smooth compact Kähler manifold. Then there is an isomorphism $\psi\colon\mathbb P^n\to B$ such that $g=\psi\circ\ell$. In particular, every holomorphic Lagrangian fibration on $Q$ with smooth compact Kähler base is the natural one, up to an isomorphism of the base.
\end{Cor}

\begin{proof}
Let $\omega_B$ be a Kähler form on $B$, let $F$ be a general fiber of $\ell$, and choose a Kähler form $\omega_F$ on the torus $F$. By \cref{prop:fiber-restriction}, the restriction of $g^*[\omega_B]$ to $F$ is zero. Hence
\[
\int_F(g|_F)^*\omega_B\wedge\omega_F^{n-1}=0.
\]
The integrand is nonnegative and would be positive on an open set if $g|_F$ were nonconstant. Thus $g$ is constant on a general fiber of $\ell$.

Let $I\subset\mathbb P^n\times B$ be the image of $(\ell,g)$, with its reduced structure. It is a compact irreducible analytic space of dimension $n$. As in the proof of \cref{lem:equidimensional-factorization}, equidimensionality and connectedness of the fibers of $\ell$ show that the first projection is finite and birational, hence an isomorphism. Therefore $g=\psi\circ\ell$ for a surjective holomorphic map $\psi\colon\mathbb P^n\to B$.

Write $\psi^*[\omega_B]=cH$, where $H=c_1(\mathcal O_{\mathbb P^n}(1))$. Since $\psi$ is nonconstant, integration against a positive $(n-1)$-st power of the Fubini--Study form gives $c>0$. A positive-dimensional fiber of $\psi$ would contain a projective curve on which $\psi^*[\omega_B]$ has zero integral, contradicting $c>0$. Thus $\psi$ is finite. Its fibers are connected because those of $g$ are connected, so it is birational and hence an isomorphism.
\end{proof}

\begin{Rem}\label{rem:twist-dichotomy}
Abasheva announces in \cite[Remark~6.3.13]{Abasheva2026v5} that the standard $\BG_n$ manifold is a Shafarevich--Tate twist, in the sense of \cite{AbashevaRogovST}, of the Lagrangian fibration $\ell_K\colon\Kum_n(E\times F)\to\mathbb P^n$. Her dichotomy \cite[Theorem~C]{Abasheva2026v5} distinguishes twists which preserve $b_2$ and carry a class restricting to an ample class on a smooth fiber from twists with $b_2$ smaller by one, all of whose classes of degree two restrict trivially to a smooth fiber. Twists of the second kind are not of class $\mathcal C$ \cite[Theorem~D]{Abasheva2026v5}. \cref{prop:fiber-restriction} and $b_2=6=b_2(\Kum_n)-1$ place the standard $\BG_n$ manifold in the second case, directly in Guan's model. For a twist of the second kind satisfying a Fujiki relation, the form is necessarily degenerate as soon as the pullback of the hyperplane class is nonzero: the Fujiki constant is nonzero, as in \cref{lem:fujiki-null}, so this class is isotropic, and by the polarization formula of \cref{rem:kahler-contrast} its pairing with a class $\alpha$ is determined by $\int h^n\alpha^n$, the integral of $\alpha^n$ over a smooth fiber, which vanishes.
\end{Rem}

\section{Euler characteristics}\label{sec:euler}

Let $n\ge2$, $Q=Q_{n,t}$, $m=n+1$ and $Y=M/S_m$ as in \cref{sec:preliminaries}. We compute the topological and holomorphic Euler characteristics from the quotient construction. Neither computation requires the full Hodge diamond or even the degeneration of the Frölicher spectral sequence.

For generalized Kummer manifolds, these invariants are already known. The topological Euler characteristic was computed by G\"ottsche; an elementary proof using a Lagrangian fibration is given in \cite{Gulbrandsen2007}. The holomorphic Euler characteristic is $\chi(\mathcal O_{\Kum_n(A)})=n+1$, as follows from the description of holomorphic forms on an irreducible hyperkähler manifold \cite{Beauville1983}. We show that both values are shared by $\BG_n$-type manifolds in every dimension. This numerical agreement complements the differences in low-degree cohomology and in the BBF form established in \cref{prop:low,thm:bbf-radical}.

\begin{Lem}\label{lem:all-sectors}
For $g\in S_m$, let $c(g)$ be its number of cycles, including cycles of length one, and set $r_g=m-c(g)$. Let $C(g)$ be its centralizer and set $Z_g=M^g/C(g)$, allowing a disconnected quotient. Then
\[
H^k(Q,\C)\simeq\bigoplus_{[g]}H^{k-2r_g}(Z_g,\C),\qquad H^q(Q,\mathcal O_Q)\simeq H^q(M,\mathcal O_M)^{S_m}.
\]
The sum is over conjugacy classes in $S_m$, and groups with a negative degree are zero.
\end{Lem}
\begin{proof}
Let $H_g$ be the product of the symmetric groups on the cycles of $g$. Every stabilizer is a Young subgroup by \cref{prop:partition-sectors}. If such a subgroup contains $g$, each cycle of $g$ lies in one of its blocks, so it contains $H_g$. Therefore $M^g=M^{H_g}$.

The group $N_{S_m}(H_g)/H_g$ permutes equal-sized blocks. The centralizer $C(g)$ induces the same group on $M^{H_g}$ since rotations within a cycle act trivially, while cycles of equal length can be exchanged. Consequently $M^g/C(g)\simeq M^{H_g}/N_{S_m}(H_g)$. By \cref{prop:partition-sectors}, this is the disjoint union of the normalized stratum closures associated with the component orbits of this partition type. It therefore includes exactly one copy of each connected sector, with shift $r_g$.

\cref{lem:analytic-decomposition} gives the first formula. For $H^q(Q,\mathcal O_Q)$, every nonidentity sector has a positive shift in holomorphic degree and contributes zero. Only $g=1$ remains, giving the second formula by the finite quotient description in the same lemma.
\end{proof}

The first formula reduces the topological Euler characteristic to a sum over fixed loci. Positive-dimensional fixed loci will have zero Euler characteristic, leaving a finite count. The second formula will give the holomorphic Euler characteristic from invariant antiholomorphic forms.

\begin{Thm}\label{thm:top-euler}
Every $\BG_n$-type manifold $X$, $n\ge2$, satisfies
\[
e(X)=e(\Kum_n(A))=(n+1)^3\sigma_1(n+1),
\]
where $\sigma_1(m):=\sum_{d\mid m}d$.
\end{Thm}
\begin{proof}
We first work with the standard representative $Q$. All shifts in \cref{lem:all-sectors} are even, so $e(Q)=\sum_{[g]}e(M^g/C(g))$. For a finite group $G$ acting on a compact manifold $W$, the Euler-characteristic form of Burnside's formula is $e(W/G)=|G|^{-1}\sum_{h\in G}e(W^h)$. It follows by applying the usual orbit-counting identity to the cells of an equivariant triangulation. Writing $M^{g,h}=M^g\cap M^h$, we obtain
\[
e(Q)=\sum_{[g]}\frac1{|C(g)|}\sum_{h\in C(g)}e(M^{g,h})     =\frac1{m!}\sum_{gh=hg}e(M^{g,h}).
\]

Let $H=\langle g,h\rangle$, and let $a$ be the number of its orbits on $\{1,\ldots,m\}$. For the standard representation $V$, one has $\dim V^H=a-1$. The equivariant coframe of \cref{lem:holomorphic-cover} therefore restricts to a global real coframe of rank $4(a-1)$ on each nonempty component of $M^H$. If $a>1$, each such component is a compact positive-dimensional manifold with a nowhere-vanishing vector field. Its Euler characteristic is zero by the Poincaré--Hopf theorem.

If $a=1$, a stabilizer containing $H$ is a transitive Young subgroup and hence equals $S_m$. Thus $M^H=M^{S_m}$. We count these full-stabilizer points using the holomorphic bundle model. Their base coordinates are $(b,\ldots,b)$ with $b\in E[m]$. At such a point, $\cN|_{(b,\ldots,b)}=\cL_0|_b^{\otimes m} \simeq\cL|_b$. The lines $\cP_i$ therefore have canonical units, compatible with the product trivialization and the $S_m$-action. In these coordinates, the fixed points of the fiber $T_F$ are $(z,\ldots,z)$ with $z\in F[m]$. Hence $M^{S_m}$ consists of $\#E[m]\cdot\#F[m]=m^4$ points.

It remains to count the ordered commuting pairs acting transitively on $m$ labelled points. Such a pair is a transitive $\Z^2$-action, hence isomorphic to the action on $\Z^2/\Lambda$ for an index-$m$ subgroup $\Lambda\subset\Z^2$. The group of equivariant automorphisms of this action is $\Z^2/\Lambda$, of order $m$. Thus each subgroup gives $m!/m=(m-1)!$ labelled actions. In Hermite normal form, $\Lambda$ has basis $(u,0),(v,w)$ with $u,w>0$, $uw=m$ and $0\leq v<u$. There are consequently $\sum_{u\mid m}u=\sigma_1(m)$ such subgroups, and $(m-1)!\sigma_1(m)$ transitive pairs. Only these pairs contribute to the Euler sum, giving
\[
e(Q)=\frac{m^4}{m!}(m-1)!\sigma_1(m)
     =m^3\sigma_1(m).
\]
Betti numbers are constant in smooth proper deformation families, which extends the formula to every $\BG_n$-type manifold.
\end{proof}

In particular, $e(\BG_2)=108$ and $e(\BG_3)=448$. The formula agrees with the generalized Kummer value \cite[Theorem~1.1]{Gulbrandsen2007}, with index $m=n+1$ in that paper.

In the next section, we will recover $e(\BG_2)=108$ directly from the explicit Lagrangian fibration over $\mathbb P^2$; see \cref{subsec:cusp-fibers}. Only the fibers over the nine cusps of the discriminant contribute, each with Euler characteristic $12$. This gives a separate geometric derivation of the fourfold value and relates the sector count above to the singular fibers of the fibration.

\begin{Rem}\label{rem:fixed-count}
The count $\#M^{S_m}=m^4$ used above can be checked without the bundle model. As in the proof of \cref{thm:top-euler}, an equivariant triangulation \cite{Illman1978} gives $e(Y)=\frac1{m!}\sum_{g\in S_m}e(M^g)$, and only the $(m-1)!$ cycles of length $m$ contribute, each with $M^g=M^{S_m}$. Hence $e(Y)=\#M^{S_m}/m$. On the other hand, by Nomizu's theorem $e(Y)$ is the Euler characteristic of the invariant complex, which averaging the alternating traces on the exterior algebra expresses as $\frac1{m!}\sum_{g\in S_m}\det(1-g\mid V)^4$. The determinant vanishes unless $g$ is a cycle of length $m$, in which case it equals $m$. Thus $e(Y)=m^3$ and $\#M^{S_m}=m^4$. For $m=3$ this agrees with \cref{lem:deRham}.
\end{Rem}

\begin{Prop}\label{prop:holomorphic-euler}
For a standard $\BG_n$-type manifold $Q$, let $\sigma_Q$ denote its holomorphic symplectic form. Then $H^{2j}(Q,\mathcal O_Q)=\C\cdot[\bar\sigma_Q^{\,j}]$ for $0\leq j\leq n$, and $H^q(Q,\mathcal O_Q)=0$ for odd $q$. Every $\BG_n$-type manifold has $\chi(\mathcal O)=n+1$.
\end{Prop}
\begin{proof}
By \cref{lem:all-sectors}, it suffices first to compute $H^q(M,\mathcal O_M)^{S_m}$. The invariant Dolbeault complex computes this group, as explained in {\cref{subsec:invariant-forms}}. Its differential vanishes in holomorphic degree zero, since $\dbar\bar x_i=\dbar\bar y_i=0$. The antiholomorphic coframe consists of two copies of $V$, so the desired space is $(\Lambda^q(V\oplus V))^{S_m}$.

The representations $\Lambda^aV$, $0\leq a\leq m-1$, are pairwise nonisomorphic self-dual irreducibles, corresponding to the hook partitions $(m-a,1^a)$ \cite{FultonHarris1991}. By Schur's lemma, $\dim(\Lambda^aV\otimes\Lambda^bV)^{S_m}$ is one if $a=b$ and zero otherwise. The decomposition
\[
\Lambda^q(V\oplus V)= \bigoplus_{a+b=q}\Lambda^aV\otimes\Lambda^bV
\]
therefore has a one-dimensional invariant subspace for $q=0,2,\ldots,2n$ and no invariants for odd $q$.

To identify the generator on $Q$, observe that $[\bar\sigma_Q^{\,j}]$ is nonzero for $1\leq j\leq n$. Otherwise $\bar\sigma_Q^{\,j}=\dbar\alpha$ for a form $\alpha$ of type $(0,2j-1)$, and closedness of $\sigma_Q$ would give, with the $\partial\alpha$ term vanishing by bidegree,
\[
0=\int_Q d\bigl(\alpha\bar\sigma_Q^{\,n-j}\sigma_Q^n\bigr)  =\int_Q\bar\sigma_Q^{\,n}\sigma_Q^n\ne0.
\]
The last integral is nonzero because $\sigma_Q^n$ is a nowhere-vanishing holomorphic volume form. These classes thus span the even groups already computed, and the degree-zero generator is the constant function one. Taking the alternating sum gives $\chi(\mathcal O_Q)=n+1$. This holomorphic index is constant in smooth proper deformation families, proving the final assertion.
\end{proof}

Together, \cref{thm:top-euler,prop:holomorphic-euler} give
\[
e(\BG_n)=e(\Kum_n(A)),\qquad \chi(\mathcal O_{\BG_n})=\chi(\mathcal O_{\Kum_n(A)})=n+1.
\]
These equalities hold throughout the corresponding deformation classes. The fourfold and sixfold calculations below will make explicit how they coexist with different Betti numbers and Hodge diamonds.

\begin{Rem}\label{rem:abasheva-deev}
A shorter proof would follow from the interpretation announced in \cite[Remark~6.3.13]{Abasheva2026v5}: the standard $\BG_n$-type manifold is a Shafarevich--Tate twist of $\Kum_n(E\times F)$ over $\mathbb P^n$. Once this identification is established, the {Lagrangian fibrations $\ell\colon Q\to\mathbb P^n$ of \cref{def:lagrangian} and $\ell_K\colon\Kum_n(E\times F)\to\mathbb P^n$, which sends a subscheme to the image in $E$ of its Hilbert--Chow cycle,} are locally isomorphic over the whole base, including the singular fibers. Their fiber Euler functions coincide, and integration with respect to compactly supported Euler characteristic gives
\[
e(Q)=\int_{\mathbb P^n}e({\ell^{-1}(b)})\,d\chi_c(b)  =\int_{\mathbb P^n}e({\ell_K^{-1}(b)})\,d\chi_c(b) =e(\Kum_n(E\times F)).
\]
The Kummer formula then proves \cref{thm:top-euler}. The result given in \cite[Proposition 6.2.1]{Abasheva2026v5} would also give $\chi(\mathcal O_Q)=n+1$. Thus the announced identification would recover both Euler characteristics, and deformation invariance extends them to every $\BG_n$-type manifold.
\end{Rem}

\section{Bogomolov-Guan fourfolds}
\label{sec:fourfold}

Let $Q=Q_{2,t}$ be a standard $\BG_2$ manifold. Write $M=M_{2,t}$, $Y=M/S_3$, and let $\pi:M\to Y$ and $f:Q\to Y$ be the quotient and resolution maps. We use the root-lattice covering of \cref{sec:preliminaries}. We first identify the three types of sectors and compute their Dolbeault cohomology. We also compute ordinary cohomology independently, checking \cref{thm:all-degeneration} in this dimension. In every diamond below, the row of total degree $k$ is ordered from $h^{k,0}$ on the left towards $h^{0,k}$ on the right, omitting bidegrees outside the dimension of the space. We first recall the partition sectors relevant in this particular case.

\begin{Prop}\label{prop:geometry}
For the standard resolution $f:Q\to Y$:
\begin{enumerate}
\item $f$ is projective in the complex analytic sense and is semismall.
\item The normalization of the codimension-two stratum closure is one primary Kodaira surface $K$. The generic exceptional fiber is $\mathbb P^1$ and its top-cohomology local system is constant.
\item There are $81$ point strata, all with full stabilizer $S_3$. The fiber over each is $\Hilb^3(\C^2)_0$, an irreducible surface. There are no other relevant strata.
\end{enumerate}
\end{Prop}

\begin{proof}
The properties of $f$ in (1) follow from \cref{{lem:projectivesemismall}} and \cref{prop:partition-sectors} shows that the stabilizers are Young subgroups. Thus the only types are $(1,1,1)$, $(2,1)$ and $(3)$, with shifts $r=0,1,2$, respectively. In particular, a cyclic subgroup of order three cannot occur as a stabilizer.

Let $\tau=(12)$ and put $K=M^\tau$. The connectedness argument in \cref{prop:low} applies here: $K$ is a principal $F$-bundle over $\{(a,a,-2a):a\in E\}\simeq E$, with connected fiber $\{(z,z,-2z):z\in F\}\simeq F$. On $K$, the coframe of \cref{lem:holomorphic-cover} restricts as $x_1=x_2=\xi$ and $y_1=y_2=\eta$, and
\[
d\xi=0,\qquad d\eta=3\xi\bar\xi.
\]
The restrictions $\xi,\eta$ form a coframe, since the fixed tangent space has complex dimension two. These are the structure equations of a primary Kodaira surface. Moreover $N_{S_3}(\langle\tau\rangle)/\langle\tau\rangle$ is trivial. Consequently $K$ itself is the normalization of the transposition-stratum closure. The transverse quotient is $\C^2/\{\pm1\}$, whose resolution has fiber $\mathbb P^1$.

It remains to count the full-stabilizer points. In the model of \cref{lem:holomorphic-cover}, their base coordinates are $(a,a,a)$ with $a\in E[3]$. At such a point, $\cN|_{(a,a,a)}=\cL_0|_a^{\otimes3}\simeq \cL|_a$, so every $\cP_i$ has a canonical unit. Using these units identifies the fiber with $T_F$ and its $S_3$-action with permutation of the coordinates. Its fixed points are $(z,z,z)$ with $z\in F[3]$. There are therefore $\#E[3]\cdot\#F[3] = 9\cdot9 = 81$ such points. Each is fixed by all of $S_3$, and hence gives a distinct point of $Y$. The fiber of the centred Hilbert--Chow resolution there is $\Hilb^3(\C^2)_0$, irreducible of dimension two \cite{Briancon1977}. This accounts for every possible partition.
\end{proof}

The shifts of the three sector types give the following additive
decompositions. A cohomology group with an index outside its usual
range is understood to be zero.

\begin{Thm}\label{thm:sectors}
For every $p,q$, there is an isomorphism
\begin{equation}
H^q(Q,\Omega_Q^p)\simeq H^q(M,\Omega_M^p)^{S_3}\oplus H^{q-1}(K,\Omega_K^{p-1})\oplus E^{p,q},
\label{eq:dolbeault-sectors}
\end{equation}
where $E^{2,2}=\C^{81}$ and all other $E^{p,q}$ vanish. There is also a topological decomposition
\begin{equation}
H^k(Q,\C)\simeq H^k(Y,\C)\oplus H^{k-2}(K,\C) \oplus\begin{cases}\C^{81},&k=4,\\0,&k\ne4.\end{cases}
\label{eq:topological-sectors}
\end{equation}
\end{Thm}

\begin{proof}
Apply \cref{lem:analytic-decomposition} to the strata in \cref{prop:geometry}. Their normalizations are $Y$, the Kodaira surface $K$, and $81$ points, with $r=0,1,2$, respectively. The top-fiber local systems are constant of rank one. The quotient identity $\Omega_Y^{[p]}=(\pi_*\Omega_M^p)^{S_3}$ gives the identity sector. The shifts in \eqref{eq:general-dolbeault} and \eqref{eq:general-topological} give the two formulas.
\end{proof}

We now describe the invariant Dolbeault complex explicitly. Choose invariant $(1,0)$-forms $x_1,x_2,y_1,y_2$ such that $dx_i=0$ and
\begin{align}
dy_1&=-x_1\bar x_1+2x_2\bar x_2+x_1\bar x_2+x_2\bar x_1,
\label{eq:dy1}\\
dy_2&=2x_1\bar x_1-x_2\bar x_2+x_1\bar x_2+x_2\bar x_1.
\label{eq:dy2}
\end{align}
On each of the ordered pairs $(x_1,x_2)$, $(y_1,y_2)$ and their conjugates, the generators $\tau=(12)$ and $g=(123)$ act by pullback as
\begin{equation}
\tau:(z_1,z_2)\longmapsto(z_2,z_1),\qquad
g:(z_1,z_2)\longmapsto(z_2,-z_1-z_2).
\label{eq:action}
\end{equation}
It is convenient to introduce $x_3=-x_1-x_2$ and $y_3=-y_1-y_2$. In these variables the structure equations take the symmetric form
\begin{equation}
dy_i=\sum_{j=1}^3x_j\bar x_j-3x_i\bar x_i,
\qquad i=1,2,3.
\label{eq:symmetric-equations}
\end{equation}

\begin{Lem}\label{lem:S3-equivariant-dolbeault}
The action of $S_3$ on the invariant coframe commutes with $d$ and $\dbar$. Moreover, the complex structure is abelian and
\[
\dbar y_i=dy_i,\qquad
\dbar x_i=\dbar\bar x_i=\dbar\bar y_i=0.
\]
\end{Lem}

\begin{proof}
With the convention $x_1+x_2+x_3=y_1+y_2+y_3=0$, the formulas in \eqref{eq:action} are precisely the permutation action on the three indices: $\tau$ exchanges the first two indices, whereas $g$ acts cyclically. Since the right-hand side of \eqref{eq:symmetric-equations} is equivariant under simultaneous permutation of the indices, the $S_3$-action commutes with $d$, and hence with its $(0,1)$-component $\dbar$. Finally, $d\Lambda^{1,0}\subset\Lambda^{1,1}$, so the complex structure is abelian. The displayed formulas for $\dbar$ follow immediately from the structure equations.
\end{proof}

Let $\mathfrak g$ denote the Lie algebra of invariant vector fields and set
\[
A^{p,q}:=\bigl(\Lambda^{p,q}\mathfrak g_\C^*\bigr)^{S_3}.
\]
Thus $(A^{\bullet,\bullet},\dbar)$ is the $S_3$-invariant part of the Dolbeault complex of invariant forms.

\begin{Prop}\label{prop:invariant-dolbeault-complex}
The inclusion of invariant forms induces natural isomorphisms
\begin{equation}
H^{p,q}_{\dbar}(M)^{S_3}
\simeq
H^{p,q}(A^{\bullet,\bullet},\dbar).
\label{eq:invariant-complex}
\end{equation}
\end{Prop}

\begin{proof}
Since the complex structure is abelian, the invariant-forms theorem for nilmanifolds with abelian complex structure gives
\[
H^{p,q}_{\dbar}(\mathfrak g_\C)\simeq H^{p,q}_{\dbar}(M);
\]
see \cite[Section~6, Remark on abelian complex structures]{ConsoleFino2001}. The inclusion of invariant forms is $S_3$-equivariant by \cref{lem:S3-equivariant-dolbeault}. Since $S_3$ is finite and we work over $\C$, taking invariants is exact: the Reynolds operator $6^{-1}\sum_{\gamma\in S_3}\gamma$ is a projection onto the invariant subcomplex. Hence cohomology commutes with taking $S_3$-invariants, which gives \eqref{eq:invariant-complex}.
\end{proof}

We next record two invariant tensors that will be used repeatedly. For ordered pairs of one-forms $r=(r_1,r_2)$ and $s=(s_1,s_2)$, and for three such pairs $r,s,t$, put
\begin{align}
B(r,s)&:=2r_1s_1+r_1s_2+r_2s_1+2r_2s_2,
\label{eq:B}\\
T(r,s,t)&:=r_2s_1t_1+r_1s_2t_1+r_1s_1t_2
+r_1s_2t_2+r_2s_1t_2+r_2s_2t_1.
\label{eq:T}
\end{align}
The order of the one-forms in each summand is part of the definition.

\begin{Lem}\label{lem:BT-invariants}
The tensors defining $B$ and $T$ are $S_3$-invariant for the representation \eqref{eq:action}. In particular,
\[
\sigma:=B(x,y)
\]
is an $S_3$-invariant closed nondegenerate holomorphic two-form.
\end{Lem}

\begin{proof}
The invariance of the coefficient tensors is checked on the generators $\tau$ and $g$ of $S_3$. Alternatively, after introducing the third components, one has
\[
B(x,y)=\sum_{i=1}^3x_i y_i,
\]
which makes its invariance immediate. The same direct substitution in \eqref{eq:T} proves the invariance of $T$.

To check closedness, write $S=\sum_jx_j\bar x_j$. By \eqref{eq:symmetric-equations},
\[
d\sigma=-\sum_{i=1}^3x_i\,dy_i
=-\sum_{i=1}^3x_i(S-3x_i\bar x_i)
=-\Bigl(\sum_{i=1}^3x_i\Bigr)S=0.
\]
Since $\sigma$ is of type $(2,0)$, it is therefore holomorphic. Finally, relative to the bases $(x_1,x_2)$ and $(y_1,y_2)$, its coefficient matrix is
\[
\begin{pmatrix}2&1\\1&2\end{pmatrix},
\]
whose determinant is $3$. Hence $\sigma$ is nondegenerate.
\end{proof}

It remains to determine the dimensions of the invariant cochain spaces before computing the differential. We recall the appropriate form of Molien's construction.

\begin{Def}[Exterior Molien series]\label{def:exterior-molien}
Let $G$ be a finite group and let $V$ be a finite-dimensional complex $G$-module. The \emph{exterior Molien series} of $V$ is
\[
\operatorname{Mol}^{\wedge}_{G,V}(u)
:=
\sum_{p=0}^{\dim V}
\dim\bigl(\Lambda^pV\bigr)^G\,u^p.
\]
More generally, for two $G$-modules $U$ and $V$, we call
\[
\operatorname{Mol}^{\wedge}_{G;U,V}(s,t)
:=
\sum_{p,q}
\dim\bigl(\Lambda^pU\otimes\Lambda^qV\bigr)^G\,s^pt^q
\]
the \emph{bigraded exterior Molien series}. Since exterior powers vanish above the dimensions of the representations, these series are in fact polynomials.
\end{Def}

The exterior version of Molien's formula is standard, see \cite[\S2.2, Exercise~(5)]{Sturmfels2008}. See also \cite{DontenBuryWeber2018} for the closely related extended Molien series involving symmetric and exterior powers.

\begin{Lem}[Exterior Molien formula]\label{lem:exterior-molien}
For finite-dimensional $G$-modules $U$ and $V$,
\begin{equation}
\operatorname{Mol}^{\wedge}_{G;U,V}(s,t)
=
\frac{1}{|G|}
\sum_{\gamma\in G}
\det(I+s\gamma|_U)\,
\det(I+t\gamma|_V).
\label{eq:molien-formula}
\end{equation}
\end{Lem}

\begin{proof}
For every finite-dimensional $G$-module $E$, the Reynolds operator $|G|^{-1}\sum_{\gamma\in G}\gamma$ is the projection onto $E^G$, and therefore
\[
\dim E^G=\frac1{|G|}\sum_{\gamma\in G}\tr(\gamma|_E).
\]
Applying this with $E=\Lambda^pU\otimes\Lambda^qV$ and using multiplicativity of the trace on tensor products gives
\[
\dim\bigl(\Lambda^pU\otimes\Lambda^qV\bigr)^G
=
\frac1{|G|}\sum_{\gamma\in G}
\tr(\gamma|\Lambda^pU)\,
\tr(\gamma|\Lambda^qV).
\]
Summing over $p$ and $q$ and using
\[
\sum_p\tr(\gamma|\Lambda^pU)s^p
=\det(I+s\gamma|_U)
\]
gives \eqref{eq:molien-formula}.
\end{proof}

We apply this to the $S_3$-modules
\[
W^{1,0}:=\langle x_1,x_2,y_1,y_2\rangle,
\qquad
W^{0,1}:=\langle\bar x_1,\bar x_2,\bar y_1,\bar y_2\rangle.
\]
Set
\[
C^{p,q}:=
\dim\bigl(\Lambda^pW^{1,0}\otimes\Lambda^qW^{0,1}\bigr)^{S_3},
\qquad
C(s,t):=\sum_{p,q}C^{p,q}s^pt^q.
\]
Thus $C(s,t)=\operatorname{Mol}^{\wedge}_{S_3;W^{1,0},W^{0,1}}(s,t)$, and \cref{lem:exterior-molien} gives
\begin{equation}
C(s,t)
=
\frac16\sum_{\gamma\in S_3}
\det(I+s\gamma|_{W^{1,0}})
\det(I+t\gamma|_{W^{0,1}}).
\label{eq:molien-dolbeault}
\end{equation}
The coefficient $C^{p,q}$ is therefore the dimension of the invariant Dolbeault \emph{cochain} space in bidegree $(p,q)$. It should not be confused with the Hodge number $h^{p,q}$, which is obtained only after taking the cohomology of $(A^{\bullet,\bullet},\dbar)$, or equivalently after subtracting the ranks of the incoming and outgoing differentials.

\begin{Prop}\label{prop:invariant}
The invariant cochain dimensions and differential ranks are
\[
(\dim A^{p,q})_{p,q=0}^4=
\begin{pmatrix}
1&0&1&0&1\\0&4&2&4&0\\1&2&11&2&1\\0&4&2&4&0\\1&0&1&0&1
\end{pmatrix},
\]
\[
\bigl(\rank(\dbar:A^{p,q}\to A^{p,q+1})\bigr)
=\begin{pmatrix}
0&0&0&0\\0&1&1&0\\0&2&2&0\\0&1&1&0\\0&0&0&0
\end{pmatrix},
\]
In both matrices, rows are indexed by $p$ and columns by $q$, and
the rank matrix omits the zero outgoing maps for $q=4$.
Writing $h^{p,q}_{\mathrm{id}}=\dim H^q(M,\Omega_M^p)^{S_3}$,
the identity sector has Hodge diamond
\begin{equation}
\begin{array}{ccccccccc}
&&&&1&&&&\\
&&&0&&0&&&\\
&&1&&3&&1&&\\
&0&&0&&0&&0&\\
1&&3&&7&&3&&1\\
&0&&0&&0&&0&\\
&&1&&3&&1&&\\
&&&0&&0&&&\\
&&&&1&&&&
\end{array}.
\label{eq:identity-table}
\end{equation}
\end{Prop}

\begin{proof}
Both coframe spaces $W^{1,0}$ and $W^{0,1}$ are isomorphic to $V\oplus V$, where $V$ is the standard two-dimensional representation of $S_3$. On $V$, the identity has eigenvalues $(1,1)$, a transposition has eigenvalues $(1,-1)$, and a three-cycle has eigenvalues $(\zeta,\zeta^2)$, where $\zeta^3=1$ and $\zeta\ne1$. Consequently their determinants $\det(I+sg|V)$ are $(1+s)^2$, $1-s^2$ and $1-s+s^2$, respectively. Each determinant is squared because the coframe contains two copies of $V$.

There are one identity, three transpositions and two three-cycles. Averaging over these conjugacy classes gives the exterior Molien series
\begin{align}
C(s,t)=\frac16\bigl[&(1+s)^4(1+t)^4 +3(1-s^2)^2(1-t^2)^2\nonumber\\
&+2(1-s+s^2)^2(1-t+t^2)^2\bigr].
\label{eq:molien}
\end{align}
The coefficient of $s^pt^q$ is the entry $C^{p,q}$ of the first matrix, with rows indexed by $p$ and columns by $q$. For example, the coefficient of $st$ is $(16+2\cdot4)/6=4$, so $C^{1,1}=4$.

We compute the ranks as follows. The space $A^{1,1}$ has basis
\[
B(x,\bar x),\quad B(y,\bar x),\quad
B(x,\bar y),\quad B(y,\bar y).
\]
The first three are $\dbar$-closed, whereas
\begin{equation}
\dbar B(y,\bar y)=3T(x,\bar x,\bar y)\ne0.
\label{eq:rank11}
\end{equation}
Thus $\rank\dbar^{1,1}=1$. The two independent invariant forms $T(x,y,\bar x)$ and $T(x,y,\bar y)$ form a basis of $A^{2,1}$, and
\begin{align}
\dbar T(x,y,\bar x)&=-6x_1x_2\bar x_1\bar x_2,\label{eq:rank21a}\\ \dbar T(x,y,\bar y)&=3x_1x_2   (\bar x_2\bar y_1-\bar x_1\bar y_2).\label{eq:rank21b}
\end{align}
The images are nonzero and independent: the first contains two $\bar x$ factors, whereas every monomial of the second contains one $\bar x$ and one $\bar y$. Thus $\rank\dbar^{2,1}=2$.

Wedge multiplication by $\sigma$ identifies the rows $p=1$ and $p=3$ as complexes. The wedge pairing $A^{p,q}\times A^{4-p,4-q}\to A^{4,4}\simeq\C$ is perfect. Indeed, the full exterior pairing is perfect and its top form is $S_3$-invariant, and averaging a partner of an invariant form preserves the pairing. Moreover $\dbar$ vanishes in bidegree $(4,3)$: a form of holomorphic degree four contains both $x_1,x_2$, so replacing any $y_i$ by $dy_i$ gives a repeated $x$ factor. The Leibniz rule therefore makes the differentials transposes, up to sign, in complementary bidegrees. Hence
\[
\rank\dbar^{p,q}=\rank\dbar^{4-p,3-q}.
\]
The nonzero ranks now follow from \eqref{eq:rank11}--\eqref{eq:rank21b}. The rows $p=0,4$ have zero differential, and $A^{2,0}$ is generated by $\sigma$, so the remaining maps vanish as stated. Finally,
\[
h^{p,q}_{\mathrm{id}}=\dim A^{p,q} -\rank\dbar^{p,q}-\rank\dbar^{p,q-1}.
\]
\end{proof}

\subsection{Hodge diamond}

\begin{Lem}\label{lem:deRham}
The invariant de Rham algebra $\mathcal A=(\Lambda^\bullet\mathfrak g^*_{\C})^{S_3}$ computes the ordinary cohomology of $Y$. Its Betti numbers are
\begin{equation}
(b_k(Y))_{k=0}^8=(1,0,5,0,15,0,5,0,1).
\label{eq:betti-Y}
\end{equation}
\end{Lem}
\begin{proof}
Nomizu's theorem \cite{Nomizu1954}, followed by exactness of finite-group invariants, gives the calculation. Invariant smooth forms on the quotient orbifold resolve its constant sheaf, so this computes the ordinary cohomology of the coarse space $Y$, not Chen--Ruan cohomology. The cochain dimensions in degrees $0,\ldots,8$ are
\[
(1,0,6,4,21,4,6,0,1).
\]
Conjugating \eqref{eq:dy1}--\eqref{eq:dy2} gives $d\bar y_i=-dy_i$. The differential ranks are
\[
(0,0,1,3,3,1,0,0,0).
\]
Here the degree-two kernel is spanned by $B(x,y),B(x,\bar x),B(y,\bar x),B(x,\bar y),B(\bar x,\bar y)$, while the sixth invariant two-form has nonzero differential. In degree three the images of $T(x,y,\bar x)$, $T(x,y,\bar y)$, and $T(\bar x,y,\bar y)$ are independent. The first is the nonzero $(2,2)$-form $-6x_1x_2\bar x_1\bar x_2$. The second has a nonzero $(3,1)$ component and no $(1,3)$ component. The third has a nonzero $(1,3)$ component and no $(3,1)$ component, as follows by conjugating and using \eqref{eq:rank21b}. Consequently the degree-three differential has rank at least three. It has rank at most three because $\dim\mathcal A^3=4$ and its kernel contains the one-dimensional image of $d:\mathcal A^2\to \mathcal A^3$. Its rank is therefore exactly three. The invariant top-degree pairing is perfect by averaging, as above. Since $\mathcal A^7=0$, the Leibniz rule gives $\rank d^k=\rank d^{7-k}$, determining the remaining ranks. This proves \eqref{eq:betti-Y}.
\end{proof}

\begin{Lem}\label{lem:kodaira-cohomology}
The Kodaira surface sector has Hodge diamond
\[
\begin{array}{ccccc}
&&1&&\\
&1&&2&\\
1&&2&&1\\
&2&&1&\\
&&1&&
\end{array}
\]
and Betti numbers $(1,3,4,3,1)$.
\end{Lem}
\begin{proof}
Use the coframe $\xi,\eta$ from \cref{prop:geometry}.
Its invariant complex structure is abelian, so its invariant forms compute Dolbeault cohomology, as in {\cref{subsec:invariant-forms}}. The only nonzero differential on generators is $\dbar\eta=3\xi\bar\xi$. In holomorphic degrees zero and two the differential is zero, and the cochain dimensions are $(1,2,1)$. In holomorphic degree one, the cochain dimensions are $(2,4,2)$. The two differential maps both have rank one, as witnessed by
\[
\dbar\eta=3\xi\bar\xi,\qquad \dbar(\eta\bar\eta)=3\xi\bar\xi\bar\eta.
\]
The remaining basis elements have zero differential. The cohomology dimensions are therefore $(1,2,1)$ in each holomorphic degree, giving the displayed diamond.

For ordinary cohomology, Nomizu's theorem applies to the same coframe with $d\bar\eta=-3\xi\bar\xi$. The de Rham cochain dimensions are $(1,4,6,4,1)$. The maps out of degrees one and two have rank one, with images spanned respectively by $\xi\bar\xi$ and $\xi\bar\xi(\eta+\bar\eta)$, and all other maps vanish. Subtracting incoming and outgoing ranks gives $(1,3,4,3,1)$, independently of the Dolbeault calculation.
\end{proof}

\begin{Thm}\label{thm:hodge}
The standard $\BG_2$ manifold has Hodge diamond
\begin{equation}
\begin{array}{ccccccccc}
&&&&1&&&&\\
&&&0&&0&&&\\
&&1&&4&&1&&\\
&0&&1&&2&&0&\\
1&&4&&90&&4&&1\\
&0&&2&&1&&0&\\
&&1&&4&&1&&\\
&&&0&&0&&&\\
&&&&1&&&&
\end{array}.
\label{eq:full-table}
\end{equation}
Its Betti numbers are $(b_k(Q))_{k=0}^8=(1,0,6,3,100,3,6,0,1)$.
\end{Thm}
\begin{proof}
Apply \eqref{eq:dolbeault-sectors} to the identity-sector diamond \eqref{eq:identity-table} and the Kodaira surface diamond in \cref{lem:kodaira-cohomology}. Shift the latter by $(1,1)$ and add $81$ in bidegree $(2,2)$. This gives \eqref{eq:full-table}. In particular, $h^{2,2}(Q)=7+2+81=90$, $h^{2,1}(Q)=1$ and $h^{1,2}(Q)=2$.

Independently, \eqref{eq:topological-sectors} gives $b_k(Q)=b_k(Y)+b_{k-2}(K)+81\delta_{k,4}$. The two de Rham calculations above therefore give the stated Betti numbers. In particular $b_4(Q)=15+4+81=100$. The diagonal sums $\sum_{p+q=k}h^{p,q}(Q)$ agree with these Betti numbers for every $k$. \cref{lem:frolicher-criterion} now proves $E_1$-degeneration.
\end{proof}

The symplectic form and Serre duality give $h^{p,q}=h^{4-p,q}=h^{4-p,4-q}$. Ordinary Hodge symmetry fails: $h^{2,1}=1$ and $h^{1,2}=2$. This asymmetry comes from the Kodaira surface contribution.

\begin{Rem}\label{rem:bg2-salamon}
    We can also note that
$\mathrm{BG}_2$ does not satisfy Salamon's relation $b_4=46+10b_2-b_3$, valid for compact
hyperk\"ahler fourfolds \cite{Salamon1996}, which would give $103$.
\end{Rem}

\subsection{Euler characteristic via its Lagrangian fibration}\label{subsec:cusp-fibers}

In this final subsection, we recover the Euler characteristic of \cref{thm:top-euler} in the case $n=2$ directly from the fibers of the algebraic reduction. We include this computation because it is completely explicit and uses the geometry of $\BG_2$-type manifolds in an essential way.

We use Bogomolov's description recalled in \cref{sec:lagrangian}, with $m=3$: the orbifold $W_F=W/F$, the map $\Pi\colon W_F\to|3\cdot0_E|\simeq\mathbb P^2$, the finite map $p_W\colon Q\to W_F$ of degree $9$, and $\Phi=\Pi\circ p_W$. Write $W_y=\Pi^{-1}(y)$ and $Q_y=\Phi^{-1}(y)$, and let $\cD\subset\mathbb P^2$ be the discriminant sextic. By \cref{prop:lagrangian}(3), $Q_y=\ell^{-1}(\psi^{-1}(y))$, where the automorphism $\psi$ preserves $\cD$ and hence its cusps. Therefore the results below also describe the fibers of $\ell$. Apart from the results of \cite{Bogomolov1996,Kuznetsova2022Thesis} quoted below, the computation is elementary. In this subsection the group law of $F$ is written additively, and $F[3]\simeq(\Z/3)^2$ is its $3$-torsion subgroup.

The points of $|3\cdot0_E|$ are of three kinds: the reduced divisors form $\mathbb P^2\setminus\cD$, the divisors $2[z]+[-2z]$ with $3z\ne0$ are the smooth points of $\cD$, and the nine divisors $3[z]$ with $z\in E[3]$ are its cusps \cite[Propositions~3.5.7 and~3.5.8, and Example~3.5.2]{Kuznetsova2022Thesis}.

\begin{Lem}\label{lem:orbifold-points}
A point of $W$ has nontrivial stabilizer in $F$ if and only if it is a reduced subspace $Z=\{s,s+g,s+2g\}$ with $g\in F[3]\setminus\{0\}$, and its stabilizer is then $\langle g\rangle\simeq\Z/3$. These points form $36$ orbits of $F$: four in each fiber $W_{3[z]}$ with $z\in E[3]$, one for each subgroup of order three of $F[3]$, and none in the other fibers. The map $p_W$ is a finite étale covering of degree $9$ over the complement of these $36$ points, and each of them has exactly three preimages in $Q$, at each of which $p_W$ has local degree three.
\end{Lem}

\begin{proof}
An element $g\ne0$ of $F$ fixing a subspace $Z$ permutes its support freely, since $F$ acts freely on $S$, and it preserves multiplicities. The support has at most three points, so it is a single orbit of $\langle g\rangle$ of size three, $g$ has order three, and $Z$ is reduced. Such a $Z$ lies over $3[\pi_S(s)]$, which belongs to $|3\cdot0_E|$ exactly when $\pi_S(s)\in E[3]$. As $F$ acts simply transitively on the fibers of $\pi_S$, the orbit of $Z$ is determined by $\pi_S(s)$ and by the subgroup $\langle g\rangle$. This gives $9\cdot4=36$ orbits.

By \cite[Corollary~3.7 and Remark~3.9]{Bogomolov1996}, see the proof of \cite[Lemma~3.4.18]{Kuznetsova2022Thesis}, the branch locus of $p_W$ is the union of the loci of subspaces of the form $\{s,s+g,s+2g\}$, and for $m=3$ these are the $36$ points above. Over their complement $W_F$ is smooth, and $p_W$ is a finite étale covering of degree $9$. 

Let $Z$ represent one of these points, and write $[Z]\in W_F$ for its class. The holomorphic slice theorem identifies a neighborhood of $[Z]$ with $B/\langle g\rangle$, where $B\subset\C^4$ is a ball on which $\langle g\rangle$ acts linearly. A subspace near $Z$ fixed by $g$ again has the form $\{s',s'+g,s'+2g\}$ with $\pi_S(s')=\pi_S(s)$, hence lies in the orbit $F\cdot Z$. Thus $g$ fixes only the center of $B$, and $(B\setminus\{0\})/\langle g\rangle$ has fundamental group $\langle g\rangle$. Let $v\in p_W^{-1}([Z])$, and let $N_v$ be the connected component of $p_W^{-1}(B/\langle g\rangle)$ containing $v$. After shrinking $B$, we may assume $p_W^{-1}([Z])\cap N_v=\{v\}$. Then $N_v\setminus\{v\}\to(B\setminus\{0\})/\langle g\rangle$ is a connected covering, hence of the form $(B\setminus\{0\})/H$ for a subgroup $H\subseteq\langle g\rangle$, and the same holds over every smaller concentric ball. Thus $v$ has arbitrarily small punctured neighborhoods with fundamental group $H$, and the inclusions between them induce isomorphisms of fundamental groups. Since $v$ is a smooth point of the fourfold $Q$, these inclusions factor through punctured balls around $v$, which are simply connected. Hence $H=1$, and $N_v\to B/\langle g\rangle$ has degree three. As $p_W$ has degree nine, $p_W^{-1}([Z])$ consists of three points.
\end{proof}

\begin{Prop}\label{prop:cusp-fibers}
Let $y\in|3\cdot0_E|$.
\begin{enumerate}
\item If $y\notin\cD$, then $Q_y$ is an abelian surface, so $e(Q_y)=0$.
\item If $y=2[z]+[-2z]$ with $3z\ne0$, then $e(W_y)=e(Q_y)=0$.
\item If $y=3[z]$ with $z\in E[3]$, then $e(W_y)=4$ and $e(Q_y)=12$.
\end{enumerate}
Consequently $e(Q)=9\cdot12=108$.
\end{Prop}

\begin{proof}
(1) This is \cite[Lemma~3.4.18]{Kuznetsova2022Thesis}, see also \cite[Proposition~3.3.9 and Lemma~3.4.17]{Kuznetsova2022Thesis} and \cref{prop:lagrangian}(1) and~(3).

(2) A subspace over $y$ is the disjoint union of a subspace $Z'$ of length two supported on $\pi_S^{-1}(z)$ and a point of $\pi_S^{-1}(-2z)$. Since $F$ acts simply transitively on $\pi_S^{-1}(-2z)$, the fiber $W_y$ is isomorphic to the space $S^{[2]}_z$ of length-two subspaces of $S$ supported on $\pi_S^{-1}(z)$. The Douady--Barlet map $S^{[2]}_z\to\Sym^2(\pi_S^{-1}(z))\simeq\Sym^2F$ is an isomorphism over the reduced pairs, and its fiber over a double point $2s$ is $\mathbb P(T_sS)\simeq\mathbb P^1$. As $\Sym^2F$ is a $\mathbb P^1$-bundle over $F$ and the diagonal is a copy of $F$, we get $e(S^{[2]}_z)=(e(\Sym^2F)-e(F))+e(F)\,e(\mathbb P^1)=0$. By \cref{lem:orbifold-points}, $p_W$ restricts to a covering of degree nine over $W_y$, so $e(Q_y)=9\,e(W_y)=0$.

(3) Identify $\pi_S^{-1}(z)$ with $F$ by an $F$-equivariant isomorphism, and stratify $W_y$ by the multiplicities of the support.

\emph{Reduced subspaces.} These are classes of reduced triples in $F$ modulo translation. The summation map $\Sym^3F\to F$ is equivariant, translation by $u$ acting on the target by $3u$. Hence $\Sym^3F/F\simeq\mathbb P_F/F[3]$, where $\mathbb P_F=|3\cdot0_F|\simeq\mathbb P^2$ is the linear system of divisors with sum zero. The nonreduced divisors in $\mathbb P_F$ form the curve $F^\vee=\{2[u]+[-2u]\}$, which is homeomorphic to $F$. So the reduced stratum is $(\mathbb P_F\setminus F^\vee)/F[3]$, and $e_c(\mathbb P_F\setminus F^\vee)=3$. A reduced divisor with sum zero is fixed by a translation $g\ne0$ exactly when it is a coset of $\langle g\rangle$ in $F[3]$. There are $4\cdot3=12$ such divisors, forming four $F[3]$-orbits, with stabilizers of order three. Hence the reduced stratum has compactly supported Euler characteristic $(3-12)/9+4=3$. It contains exactly the four points of \cref{lem:orbifold-points} lying in $W_y$.

\emph{Subspaces with one double point.} Here $Z=Z_2\sqcup\{s'\}$ with $Z_2$ of length two at $s\ne s'$, and $Z_2$ corresponds to a line in $T_sS$. Normalizing $s$ by translation, this stratum is $(F\setminus\{0\})\times\mathbb P(T_sS)\simeq(F\setminus\{0\})\times\mathbb P^1$, with $e_c=-2$.

\emph{Punctual subspaces.} Normalizing the support, this stratum is the punctual Hilbert scheme $\Hilb^3(\C^2)_0$, irreducible of dimension two \cite{Briancon1977}, with Euler characteristic $3$ \cite{EllingsrudStromme1987}.

Adding up, $e(W_y)=3-2+3=4$. By \cref{lem:orbifold-points}, $p_W$ is a covering of degree nine over the complement of the four orbifold points of $W_y$, and each of these has three preimages. Hence $e(Q_y)=9\,(4-4)+4\cdot3=12$.

Finally, the function $y\mapsto e(Q_y)$ is constructible, and integration with respect to the Euler characteristic \cite{Viro1988}, \cite[Section~9.7]{KashiwaraSchapira1990} gives $e(Q)=\int_{\mathbb P^2}e(Q_y)\,d\chi(y)$. By (1)--(3) the integrand vanishes outside the nine cusps, where it equals $12$.
\end{proof}

\begin{Rem}\label{rem:cusp-picture}
The value $108$ agrees with \cref{thm:top-euler} and with the sector decomposition $e(Q)=e(Y)+e(K)+81=27+0+81$ of \cref{thm:sectors}, where $e(Y)=27$ and $e(K)=0$ by \cref{lem:deRham,lem:kodaira-cohomology}. \Cref{rem:cusp-guan} below shows how these three terms are distributed over the cusp fibers.

The three strata of $W_{3[z]}$ have irreducible closures of dimension two, which are the irreducible components of this fiber. The punctual component is a projective cone over a twisted cubic \cite{Briancon1977}. The closure $\Sigma_z$ of the reduced stratum consists of the subspaces of length three of the elliptic curve $\pi_S^{-1}(z)$, so that $\Sigma_z\simeq\Sym^3F/F\simeq\mathbb P_F/F[3]$. Each $g\in F[3]\setminus\{0\}$ acts on $\mathbb P_F\simeq\mathbb P^2$ with exactly three fixed points, so its linear lift has three distinct eigenvalues, which differ by cube roots of unity. Hence $g$ acts on the tangent space at each fixed point with eigenvalues $e^{2\pi i/3}$ and $e^{-2\pi i/3}$. The quotient map is étale in codimension one, so $\Sigma_z$ is a Du Val del Pezzo surface of degree $9/9=1$ whose four singular points, of type $A_2$, are the orbifold points of \cref{lem:orbifold-points} in $W_{3[z]}$.

Moreover $p_W^{-1}(\Sigma_z)\simeq\mathbb P_F$. Indeed, by the proof of \cref{lem:orbifold-points}, near each of the twelve points of $Q$ over the singular points $[Z]$ of $\Sigma_z$ the map $p_W$ is the quotient map $B\to B/\langle g\rangle$, and the preimage of $\Sigma_z$ is the smooth $\langle g\rangle$-invariant surface germ formed by the subspaces of $\pi_S^{-1}(z)$ near $Z$, on which $g$ acts freely away from the center. Hence $p_W^{-1}(\Sigma_z)$ is smooth. Over the smooth locus of $\Sigma_z$, which is the quotient of $\mathbb P_F$ minus twelve points by the free action of $F[3]$, it is a covering of degree nine on which each subgroup of order three acts freely. Such a covering is a single free $F[3]$-orbit, so $p_W^{-1}(\Sigma_z)$ is the normalization of $\Sigma_z$ in the function field of $\mathbb P_F$, that is, $\mathbb P_F$ itself.

The number $e(Q_{3[z]})=12$ can thus be split in two ways. First, $12=4\cdot3$ counts the points of $Q$ over the four $A_2$-points of $\Sigma_z$, since the rest of $W_{3[z]}$ has Euler characteristic zero and is covered nine times. Correspondingly, $e(Q)=3\cdot36$ counts the points of $Q$ over the $36$ orbifold points of $W_F$. Second, $12=e(\mathbb P^2)+9\,(e(W_{3[z]})-e(\Sigma_z))=3+9$, where the projective plane $p_W^{-1}(\Sigma_z)$ maps $9:1$ onto $\Sigma_z$.

\end{Rem}

\begin{Rem}\label{rem:cusp-guan}
In Guan's model the cusp fibers are equally explicit. Let $z\in E[3]$. As in the proof of \cref{prop:geometry}, the canonical units identify the fiber of $M$ over $(z,z,z)\in B_E$ with $T_F$, the group $S_3$ permuting the coordinates. Hence $\ell_Y^{-1}(3[z])\simeq T_F/S_3$, which is the linear system $\mathbb P_F\simeq\mathbb P^2$ by Abel's theorem. Its points with nontrivial stabilizer form the curve $F^\vee$, and the nine points with stabilizer $S_3$ are the cusps $3[u]$, $u\in F[3]$, of $F^\vee$. Over $F^\vee$ minus its cusps, $f$ has fibers $\mathbb P^1$, and over the cusps it has fibers $\Hilb^3(\C^2)_0$ (\cref{prop:geometry}). Hence
\[
e\bigl(\ell^{-1}(3[z])\bigr)=e_c(\mathbb P_F\setminus F^\vee)+2\,e_c(F^\vee\setminus\{\text{cusps}\})+9\cdot3=3-18+27=12.
\]
The other fibers of $\ell_Y$ are $T_F$, or its quotient by a transposition whose fixed locus is an elliptic curve. Both have Euler characteristic zero, and so do the corresponding fibers of $\ell$, as $f$ has fibers $\mathbb P^1$ over the image of that curve. Hence $e(Q)=9\cdot12$.

This is the fiberwise form of the sector count $e(Q)=e(Y)+e(K)+81$ of \cref{thm:sectors}. By the proof of \cref{prop:geometry}, the Kodaira surface $K=M^\tau$, with $\tau=(12)$, is a principal $F$-bundle over $\{(a,a,-2a):a\in E\}\simeq E$. Hence $\ell_Y$ maps $K$ onto $\cD$, all fibers of this map are elliptic curves, and over the cusp $3[z]$ the fiber $\{(u,u,-2u):u\in F\}$ maps bijectively onto $F^\vee$. Grouping the terms of the count above according to the sectors, we obtain
\[
e\bigl(\ell^{-1}(3[z])\bigr)=e(\mathbb P_F)+e(F^\vee)+9=3+0+9,
\]
where the three terms come from $Y$, from $K$ and from the nine points over $3[z]$ with stabilizer $S_3$. Over the other points of $\mathbb P^2$ each sector contributes zero. Summing over the nine cusps gives $e(Y)=27$, $e(K)=0$ and $81=9\cdot9$.
\end{Rem}

\section{Bogomolov-Guan sixfolds}\label{sec:sixfold}

Let $Q=Q_{3,t}$ be a standard $\BG_3$ manifold, and let $M=M_{3,t}$. We compute the Hodge diamond of $Q$ from the five partition sectors of the resolution $Q\to M/S_4$. Only the sectors of types $(1^4)$ and $(2,1,1)$ require new Dolbeault calculations, while the others are torus quotients, a primary Kodaira surface and points.

\subsection{Sector data}
We retain the notation $F_\lambda=M^{H_\lambda}$, $G_\lambda=N_{S_4}(H_\lambda)/H_\lambda$ and $Z_\lambda=F_\lambda/G_\lambda$ from \cref{sec:preliminaries}. For a partition $\lambda=(\lambda_1,\ldots,\lambda_l)$ of four, the shift is $r_\lambda=4-l$. On each fixed component, the block coframe satisfies the weighted zero-sum relations and
\begin{equation}\label{eq:block}
dX_a=0,\qquad dY_a=\sum_b\lambda_bX_b\bar X_b-4X_a\bar X_a, \qquad d\bar Y_a=-dY_a.
\end{equation}
Products of forms denote wedge products.

\begin{Prop}\label{prop:six-geometry}
The sector data are
\[
\begin{array}{ccccc}
\lambda & r_\lambda & \text{fixed components} & G_\lambda & Z_\lambda \\ \hline
(1^4)   & 0 & 1   & S_4 & M/S_4 \\
(2,1,1) & 1 & 1   & S_2 & L/S_2 \\
(2,2)   & 2 & 16  & S_2 & \bigsqcup_{j=1}^{16} T_j/\{\pm 1\} \\
(3,1)   & 2 & 1   & 1   & K \\
(4)     & 3 & 256 & 1   & 256 \text{ points}
\end{array}
\]
Here $L$ is a nilmanifold of complex dimension four, each $T_j$ is a complex two-dimensional torus, and $K$ is a primary Kodaira surface. The residual group preserves each component in the table.
\end{Prop}
\begin{proof}
By \cref{prop:partition-sectors}, the five partitions exhaust the sectors. The fixed loci of types $(2,1,1)$ and $(3,1)$ are connected by \cref{lem:component-bound}. Their restricted coframes have constant nilpotent structure equations, so they are compact nilmanifolds, as explained in \cref{subsec:invariant-forms}. The residual involution on $L$ exchanges the singleton blocks. For $(3,1)$, the relations give $X_2=-3X_1$, $Y_2=-3Y_1$ and $dY_1=8X_1\bar X_1$, which identifies $K$ as a primary Kodaira surface.

The full-stabilizer points lie over $(a,a,a,a)$ with $a\in E[4]$. At each such point of the base, $\cN=\cL_0|_a^{\otimes4}\simeq\cL|_a$, so the four lines $\cP_i$ have canonical units compatible with the product trivialization. The fixed fiber tuples are therefore $(z,z,z,z)$ with $z\in F[4]$. This gives $\#E[4]\cdot\#F[4]=256$ points, each fixed individually by $S_4$.

For $(2,2)$, put $H=\langle(12),(34)\rangle$. The weighted relations and \eqref{eq:block} give $X_2=-X_1$, $Y_2=-Y_1$ and $dY_1=0$. Thus every connected component of $M^H$ is a complex two-dimensional torus, and there are at most sixteen components by \cref{lem:component-bound}. The block exchange $s=(13)(24)$ acts by $-1$ on its coframe. On a component containing an $s$-fixed point, choosing that point as origin identifies $s$ with inversion, which then has exactly sixteen fixed points.

The $s$-fixed points of $M^H$ are precisely the full-stabilizer points. Indeed, their stabilizer contains the transitive group $\langle H,s\rangle$, and a Young subgroup containing this group must be $S_4$. The 256 points therefore require at least $256/16=16$ components. The bound is attained, every component is preserved by $s$, and the quotient on each is inversion. Finally, \cref{prop:partition-sectors} identifies these quotients with the sector normalizations and gives multiplicity one for each connected sector.
\end{proof}

\begin{Rem}\label{rem:two-block}
For a partition $\lambda=(\lambda_1,\lambda_2)$ of $m$ with two parts, the weighted relation gives $X_2=-\lambda_1X_1/\lambda_2$, and the block equations $dY_a=\sum_b\lambda_bX_b\bar X_b-mX_a\bar X_a$ become
\[
dY_1=\frac{m(\lambda_1-\lambda_2)}{\lambda_2}\,X_1\bar X_1.
\]
Hence the fixed components are complex tori if $\lambda_1=\lambda_2$, and primary Kodaira surfaces otherwise. This explains the tori of type $(2,2)$ above and the Kodaira surfaces of types $(2,1)$ and $(3,1)$, with coefficients $3$ and $8$.
\end{Rem}

For a sector $Z$, write $\mathcal H_Z(s,t)=\sum_{p,q}h^{p,q}(Z)s^pt^q$, using reflexive forms when $Z$ is singular. Let $C_Z^{p,q}$ be the dimension of its invariant Dolbeault cochains and $R_Z^{p,q}$ the rank of the outgoing differential. With $C_Z(s,t)=\sum C_Z^{p,q}s^pt^q$ and $R_Z(s,t)=\sum R_Z^{p,q}s^pt^q$, one has
\begin{equation}\label{eq:six-cochain-hodge}
\mathcal H_Z(s,t)=C_Z(s,t)-(1+t)R_Z(s,t).
\end{equation}
This is the identity $h^{p,q}=C_Z^{p,q}-R_Z^{p,q}-R_Z^{p,q-1}$. The invariant complexes compute the sector Dolbeault groups by \cref{subsec:invariant-forms}.

\begin{Prop}\label{prop:six-linear}
For the identity sector, the cochain polynomial is
\[
C_{1111}(s,t)=\frac1{24}\bigl( P_1(s)P_1(t)+6P_2(s)P_2(t)+3P_{22}(s)P_{22}(t) +8P_3(s)P_3(t)+6P_4(s)P_4(t)\bigr),
\]
where
\[
\begin{split}
P_1(s)&=(1+s)^6,\qquad
P_2(s)=(1+s)^2(1-s^2)^2,\\
P_{22}(s)&=(1+s)^2(1-s)^4,\qquad
P_3(s)=(1+s^3)^2,\qquad
P_4(s)=(1-s)^2(1+s^2)^2.
\end{split}
\]
Its outgoing ranks are determined by 
\begin{equation}\label{eq:six-id-ranks}
\bigl(R_{1111}^{p,q}\bigr)_{ \substack{1\leq p\leq3\\0\leq q\leq2}} =\begin{pmatrix}0&1&1\\0&2&4\\0&1&8\end{pmatrix}.
\end{equation}
For the $(2,1,1)$ sector,
\[
C_{211}(s,t)=\frac12\bigl((1+s)^4(1+t)^4 +(1-s^2)^2(1-t^2)^2\bigr),
\]
and the outgoing ranks are determined by
\begin{equation}\label{eq:six-211-ranks}
\bigl(R_{211}^{p,q}\bigr)_{
\substack{1\leq p\leq2\\0\leq q\leq1}}
=\begin{pmatrix}1&3\\1&5\end{pmatrix}.
\end{equation}
In each case, if $d$ is the complex dimension of the sector, the rows $p=0,d$ and the column $q=d$ vanish, and the remaining entries follow from
\[
R^{p,q}=R^{d-p,q}=R^{p,d-1-q}.
\]
\end{Prop}
\begin{proof}
For the identity sector, the holomorphic and antiholomorphic coframe spaces are each $V\oplus V$, where $V$ is the standard three-dimensional representation of $S_4$. The five conjugacy classes have sizes $1,6,3,8,6$. Their determinants $\det(1+sg\mid V\oplus V)$ are the polynomials $P_1,P_2,P_{22},P_3,P_4$ above. The exterior Molien formula, as in \cref{sec:fourfold}, gives $C_{1111}$.

For the ranks, order the generators as $(x_1,x_2,x_3,y_1,y_2,y_3,w_1,w_2,w_3,z_1,z_2,z_3)$, where $w_i=\bar x_i$ and $z_i=\bar y_i$. The fourth member of each family is minus the sum of the first three. The only nonzero generator differentials are
\[
\dbar y_i=\sum_{j=1}^4x_jw_j-4x_iw_i.
\]
Let $D_{p,q}$ denote this derivation on the increasing exterior monomials, and let $\mathcal R_{p,q}=\sum_{g\in S_4}g$ act on those monomials. Then $R_{1111}^{p,q}=\rank(D_{p,q}\mathcal R_{p,q})$. Thus the calculation consists of averaging each monomial, applying the displayed differential with the exterior signs, and row reducing its images over $\Q$. In the indicated bidegrees this gives \eqref{eq:six-id-ranks}. For example, the invariant $(1,1)$-cochains have basis
\[
\sum_i x_iw_i,\qquad \sum_i y_iw_i,\qquad
\sum_i x_iz_i,\qquad \sum_i y_iz_i.
\]
Only the last has nonzero differential, $-4\sum_i x_iw_iz_i$, giving $R_{1111}^{1,1}=1$.

For $(2,1,1)$, set $U=X_1$, $V=X_1+X_2$, $A=Y_1$ and $B=Y_1+Y_2$. The residual involution fixes $U,A$ and negates $V,B$, and acts in the same way on their conjugates. Hence its Molien polynomial is $C_{211}$. Equation~\eqref{eq:block} becomes
\[
\dbar A=2V\bar V,\qquad \dbar B=4(U\bar V+V\bar U),
\]
with all other generator differentials zero. The invariant monomials are exactly those containing an even number of factors from $V,B,\bar V,\bar B$. Applying the two displayed differentials to this basis gives \eqref{eq:six-211-ranks}. For instance, $\dbar(UA)=-2UV\bar V$ and $\dbar(VB)=4UV\bar V$, so $R_{211}^{2,0}=1$.
It remains to explain why these small arrays suffice. Multiplication by powers of the invariant holomorphic symplectic form identifies the complexes in rows $p$ and $d-p$. The perfect invariant top-form pairing and the Leibniz rule give $R^{p,q}=R^{d-p,d-1-q}$. Together these give the stated symmetries. The row $p=0$ has zero differential, hence so does $p=d$. The column $q=d$ vanishes for degree reasons.
\end{proof}

\subsection{Hodge Diamond}

The remaining sectors require no new rank calculations. Inversion on a two-dimensional torus acts on bidegree $(p,q)$ by $(-1)^{p+q}$, so
\[
\mathcal H_{22}(s,t)=(1+s^2)(1+t^2)+4st.
\]
By \cref{lem:kodaira-cohomology}, the Kodaira sector has $\mathcal H_{31}(s,t)=(1+s+s^2)(1+t)^2$. The coefficient eight in its structure equation is nonzero and can be removed by rescaling the generator. Combining the sector contributions gives the full Hodge diamond.

\begin{Thm}\label{thm:six-hodge}
The standard $\BG_3$ manifold has Hodge diamond
\begin{equation}\label{eq:six-hodge}
\begin{array}{ccccccccccccc}
&&&&&&1&&&&&&\\
&&&&&0&&0&&&&&\\
&&&&1&&4&&1&&&&\\
&&&0&&1&&2&&0&&&\\
&&1&&4&&27&&4&&1&&\\
&0&&1&&8&&10&&2&&0&\\
1&&4&&27&&344&&27&&4&&1\\
&0&&2&&10&&8&&1&&0&\\
&&1&&4&&27&&4&&1&&\\
&&&0&&2&&1&&0&&&\\
&&&&1&&4&&1&&&&\\
&&&&&0&&0&&&&&\\
&&&&&&1&&&&&&
\end{array}.
\end{equation}
Its Betti numbers are
\begin{equation}\label{eq:six-betti}
(b_0,\ldots,b_{12})=(1,0,6,3,37,21,408,21,37,3,6,0,1).
\end{equation}
\end{Thm}
\begin{proof}
Apply \cref{lem:analytic-decomposition} to the sector data in \cref{prop:six-geometry}. In terms of the Hodge polynomials, the result is
\[
\mathcal H_Q=\mathcal H_{1111}+st\,\mathcal H_{211}
+(st)^2(16\mathcal H_{22}+\mathcal H_{31})+256(st)^3.
\]
For the first two terms, use \eqref{eq:six-cochain-hodge} and \cref{prop:six-linear}. For the last three, use the polynomials just computed. Extracting coefficients gives \eqref{eq:six-hodge}. For example,
\[
\begin{split}
h^{2,2}(Q)&=6+4+16+1=27,\\
h^{2,3}(Q)&=2+6+0+2=10,\\
h^{3,3}(Q)&=12+10+16\cdot4+2+256=344.
\end{split}
\]
The cohomological comparison established in \cref{sec:cohomology} identifies the diagonal sums with the Betti numbers, giving \eqref{eq:six-betti}.
\end{proof}

The diamond has the symplectic and Serre symmetries $h^{p,q}=h^{6-p,q}=h^{6-p,6-q}$, but not Hodge symmetry: for example, $h^{2,3}=10$ and $h^{3,2}=8$.

\section{Characteristic numbers and cobordism}
\label{sec:chern}

For Kähler holomorphic symplectic manifolds, Hodge and characteristic-number calculations are available for the known deformation types: the two series of Beauville and the two sporadic examples of O'Grady. The Hodge numbers are treated in \cite{GottscheSoergel1993,MongardiRapagnettaSacca2018,deCataldoRapagnettaSacca2021}, and characteristic numbers and Riemann--Roch polynomials in \cite{EllingsrudGottscheLehn2001,Nieper2002,RiosOrtiz2024RR}. Göttsche's formula \cite{Gottsche1990} gives the Euler characteristics of Hilbert schemes of points on a K3 surface, while \cref{thm:top-euler} gives the generalized Kummer Euler characteristic for $\BG_n$ in every dimension. We now show that, for $n=2,3$, this agreement extends to all Chern and Pontryagin numbers.

We use the complex orientation throughout. In formulas for characteristic numbers, integration over the manifold is understood. For instance, $c_2^2$ denotes $\int_X c_2(T_X)^2$. When writing identities between characteristic classes, we retain the usual cohomological meaning of the same symbols. Put $\chi^p(X)=\sum_q(-1)^q h^{p,q}(X)$ and $\chi_y(X)=\sum_p\chi^p(X)y^p$.

The standard representatives studied in \cref{sec:fourfold,sec:sixfold} determine all Chern numbers in their respective dimensions. These numbers, and the indices $\chi^p$, are constant in smooth proper deformation families. The conclusions below therefore hold for every $\BG_2$ manifold and every $\BG_3$ manifold, and constancy of the individual Hodge numbers is not needed.

\subsection{Riemann--Roch and Chern numbers}

\begin{Lem}\label{lem:RR}
For a compact holomorphic symplectic fourfold $X$, $\chi(\mathcal O_X)=(3c_2^2-c_4)/720$. For a compact holomorphic symplectic sixfold $X$,
\begin{align}
\chi(\mathcal O_X)&=(10c_2^3-9c_2c_4+2c_6)/60480,
\label{eq:td6}\\
\chi(\Omega_X^1)&=(10c_2^3-51c_2c_4-82c_6)/10080.
\label{eq:omega1}
\end{align}
\end{Lem}
\begin{proof}
The symplectic structure reduces the structure group of $T_X$ to $\operatorname{Sp}(n,\C)$, and hence to its maximal compact subgroup $\operatorname{Sp}(n)$. Its odd Chern classes vanish integrally, and its formal Chern roots occur in pairs $x_1,-x_1,\ldots,x_n,-x_n$. The analytic Riemann--Roch theorem gives
\[
\chi_y(X)=\int_X\operatorname{td}(T_X)
          \operatorname{ch}\Bigl(\sum_p y^p\Omega_X^p\Bigr).
\]
It applies to compact complex manifolds without a Kähler assumption \cite[Section~4]{AtiyahSinger1968III}. One root pair contributes
\[
\frac{x^2(1+ye^{-x})(1+ye^x)}{4\sinh^2(x/2)} =(1+y)^2\left(1-\frac{x^2}{12}+\frac{x^4}{240}  -\frac{x^6}{6048}\right)+yx^2+O(x^8)
\]
to this integrand. Multiply over the root pairs and extract the top cohomological degree. Since $c(T_X)=\prod_i(1-x_i^2)$, the elementary symmetric functions in the $x_i^2$ are $-c_2,c_4,-c_6$. The constant coefficient in $y$ gives $\chi(\mathcal O_X)$, and the linear coefficient gives $\chi(\Omega_X^1)$. This yields the stated formulas.
\end{proof}

\begin{Thm}\label{thm:chern}
The nonzero Chern numbers of $\BG_2$ and $\BG_3$ manifolds are
\[
\begin{aligned}
\BG_2:&\qquad c_2^2=756,\quad c_4=108,\\
\BG_3:&\qquad c_2^3=30208,\quad c_2c_4=6784,\quad c_6=448.
\end{aligned}
\]
Every Chern number containing an odd-indexed Chern class is zero. These are exactly the Chern numbers of $\Kum_2(A)$ and $\Kum_3(A)$.
\end{Thm}
\begin{proof}
The Betti numbers in Theorems~\ref{thm:hodge} and~\ref{thm:six-hodge} give $e(Q_{2,t})=108$ and $e(Q_{3,t})=448$ by alternating summation. Gauss--Bonnet therefore gives $c_4=108$ and $c_6=448$. The $p=0$ entries of the two Hodge diamonds give $\chi(\mathcal O_{Q_{2,t}})=3$ and $\chi(\mathcal O_{Q_{3,t}})=4$. For $Q_{2,t}$, \cref{lem:RR} now gives $3c_2^2-108=2160$, so $c_2^2=756$.

For $Q=Q_{3,t}$, the entries with $p=1$ in \eqref{eq:six-hodge} give $\chi(\Omega_Q^1)=-4+2-4+2-4=-8$, while $\chi(\mathcal O_Q)=4$. Put $a=c_2^3$ and $b=c_2c_4$. Equations~\eqref{eq:td6}--\eqref{eq:omega1} become $10a-9b=241024$ and $10a-51b=-43904$. Subtracting yields $42b=284928$, hence $b=6784$ and $a=30208$. This determines every Chern number since the odd classes vanish. The numbers extend to every $\BG_2$ and $\BG_3$ manifold by deformation invariance. The comparison is with the independent table in \cite[Appendix~A]{Nieper2002}: the variety denoted there by $A^{[[n+1]]}$ is $\Kum_n(A)$ in our notation. Its dimension-four and dimension-six entries are precisely the numbers above.
\end{proof}

\subsection{Pontryagin numbers, genera, and cobordism}

\begin{Prop}\label{prop:pontryagin}
The complete Pontryagin numbers are
\[
\begin{aligned}
\BG_2:&\qquad p_1^2=3024,\quad p_2=972,\\
\BG_3:&\qquad p_1^3=-241664,\quad p_1p_2=-87552,\quad p_3=-14464.
\end{aligned}
\]
All Stiefel--Whitney numbers of these two manifolds vanish.
\end{Prop}
\begin{proof}
The underlying real tangent bundle has total Pontryagin class $p(T_X^{\R})=\prod_i(1+x_i^2)^2$. Using $c(T_X)=\prod_i(1-x_i^2)$ gives the class identities
\[
p_1=-2c_2,\qquad p_2=c_2^2+2c_4,\qquad p_3=-2c_2c_4-2c_6.
\]
Consequently the degree-eight numbers are $p_1^2=4c_2^2$ and $p_2=c_2^2+2c_4$. In degree twelve they are $p_1^3=-8c_2^3$, $p_1p_2=-2c_2^3-4c_2c_4$ and $p_3=-2c_2c_4-2c_6$. Substitution of \cref{thm:chern} gives the stated values. These lists contain every Pontryagin number in the respective dimensions. Since $w_{2i}=c_i\bmod2$ and $w_{2i+1}=0$, the only possibly nonzero Stiefel--Whitney numbers are reductions of the even Chern numbers above. All of those integers are even and the proposition follows.
\end{proof}

\begin{Cor}\label{cor:genera}
The holomorphic genera are
\[
\begin{split}
\chi_y(\BG_2)&=3-6y+90y^2-6y^3+3y^4,\\
\chi_y(\BG_3)&=4-8y+44y^2-336y^3+44y^4-8y^5+4y^6.
\end{split}
\]
The signatures are $84$ and $-256$, respectively. The middle-degree real cup-product pairings on $H^4(\BG_2,\R)$ and $H^6(\BG_3,\R)$ have inertias $(92,8)$ and $(76,332)$, where the entries count positive and negative directions. Both manifolds are spin, with $\widehat A$-genera $3$ and $4$.
\end{Cor}
\begin{proof}
For a standard representative, alternating summation in each fixed holomorphic degree of the Hodge diamonds gives the two $\chi_y$ polynomials. For example, $\chi^2(\BG_3)=1-1+27-10+27-1+1=44$ and $\chi^3(\BG_3)=-4+8-344+8-4=-336$. Riemann--Roch then extends these indices to every manifold in the corresponding deformation class.

For the signatures, apply Hirzebruch's signature theorem \cite[Section~19]{MilnorStasheff1974} to the Pontryagin numbers of \cref{prop:pontryagin}. The polynomials in real dimensions eight and twelve give
\[
\begin{split}
\operatorname{sign}(\BG_2)
 &=\frac{7\cdot972-3024}{45}=84,\\
\operatorname{sign}(\BG_3)
 &=\frac{2(-241664)-13(-87552)+62(-14464)}{945}=-256.
\end{split}
\]
These values also equal $\chi_y(X)$ at $y=1$, that is $\sum_p\chi^p(X)$, in agreement with the characteristic-genus identity \cite{Hirzebruch1966}. Poincaré duality makes the middle pairings nondegenerate. Their dimensions are $b_4(\BG_2)=100$ and $b_6(\BG_3)=408$, so their inertias are $((100+84)/2,(100-84)/2)=(92,8)$ and $((408-256)/2,(408+256)/2)=(76,332)$.

Finally, $w_2(T_X^{\R})=c_1(T_X)\bmod2=0$, so the complex orientation admits a spin structure. The characteristic-class identity $\operatorname{td}(T_X)=e^{c_1(T_X)/2}\widehat A(T_X^{\R})$ reduces to $\operatorname{td}(T_X)=\widehat A(T_X^{\R})$. Thus $\widehat A(X)=\chi(\mathcal O_X)$, giving $3$ and $4$.
\end{proof}

\begin{Cor}
For $n=2,3$, a $\BG_n$-type manifold and $\Kum_n(A)$ represent the same unitary cobordism class.
\end{Cor}
\begin{proof}
Chern numbers determine the integral unitary cobordism class by the theorem of Milnor and Novikov, see \cite{BuoncristianoHacon1983}. Thus the assertion follows from \cref{thm:chern}.
\end{proof}

\section{Higher-dimensional characteristic numbers and the Kummer conjecture}\label{sec:conjecture}

\subsection{The characteristic-number conjecture}

The equality of Chern numbers in \cref{thm:chern}
suggests the following.

\begin{Conj}\label{ques:chern}
For every $n\ge2$, a $\BG_n$-type manifold and $\Kum_n(A)$ have the same Chern numbers, or equivalently the same Pontryagin numbers.
\end{Conj}

The equivalence follows from the vanishing of odd Chern classes and the identities
\[
p_j(T_X^{\R})=(-1)^j\left(2c_{2j}(T_X)  +\sum_{a=1}^{j-1}c_{2a}(T_X)c_{2j-2a}(T_X)\right).
\]
These give an invertible triangular change of variables over $\Q$ between the Pontryagin classes and the even Chern classes. Consequently equality of all top-degree monomial integrals in either set implies equality in the other. The conjecture compares characteristic numbers, whereas comparing individual classes would require a specified cohomological correspondence.

The conjecture holds for $n=2,3$ by \cref{sec:chern}. The equalities of the topological and holomorphic Euler characteristics in all dimensions, discussed in {\cref{sec:euler}}, give further evidence.

\begin{Rem}\label{rem:bai-twists}
$\BG_n$-type manifolds are non-Kähler. Its announced Shafarevich-Tate twist description is recalled in \cref{rem:abasheva-deev}. The second Betti numbers of $\BG_n$-type and $\Kum_n(A)$ are six and seven, respectively, and their real BBF inertias are $(2,3,1)$ and $(3,4,0)$. Under that identification, the twisting class cannot be torsion, since a torsion twist of a projective Lagrangian fibration is projective \cite[Theorem~A]{Abasheva2026v5}. {Alternatively, a torsion class has zero image in $H^2(\mathbb P^n,R^1\ell_{K*}\Q)$, so it lies in the first case of Abasheva's dichotomy \cite[Theorem~C]{Abasheva2026v5}, and the twist would have $b_2=7$. The standard $\BG_n$ manifold realizes the second case, see \cref{rem:twist-dichotomy}.}
\cref{ques:chern} asks whether the Chern and Pontryagin numbers nevertheless remain unchanged for this non-Kähler twist. The cohomological comparisons above do not identify tangent-bundle characteristic classes. The Riemann--Roch calculation below also shows that the $\chi_y$-genus alone does not determine all Chern numbers in complex dimension eight.
\end{Rem}

\subsection{A procedure for higher dimensions}

The sector construction provides a way to investigate higher dimensions, with an additional characteristic-number calculation when the Riemann--Roch equations do not have full rank.

\begin{Alg}\label{alg:sectors}
Fix a standard $\BG_n$-type manifold $Q=Q_{n,t}$ and put $m=n+1$. Use the bundle $M\to B_E$ and the resolution $f:Q\to M/S_m$ of \cref{sec:preliminaries}.
\begin{enumerate}
\item \emph{Determine the fixed components.} For each partition $\lambda=(\lambda_1,\ldots,\lambda_l)$ of $m$, let $F_\lambda=M^{H_\lambda}$ and $\delta_\lambda=\gcd(\lambda_1,\ldots,\lambda_l)$. The local-section argument of \cref{lem:component-bound} applies in every dimension: $F_\lambda$ is a principal bundle over $\ker(\sum\lambda_a:E^l\to E)$ with structure group $\ker(\sum\lambda_a:F^l\to F)$. Each of these kernels has $\delta_\lambda^2$ components. Determine the monodromy on the fiber components to obtain the connected components of $F_\lambda$. The number $\delta_\lambda^4$ is an upper bound, and equality requires a proof using the global bundle or its lattice.

\item \emph{Form the quotient sectors.}
Compute the action of $G_\lambda=N_{S_m}(H_\lambda)/H_\lambda$ on the components and choose one representative $C$ per orbit. With $G_{\lambda,C}$ the stabilizer of $C$, the corresponding sector is $Z_{\lambda,C}=C/G_{\lambda,C}$ with shift $r_\lambda=m-l$. \cref{prop:partition-sectors} gives constant rank-one top-fiber local systems. Count multiplicities after taking component orbits, as in the $(2,2)$ calculation for the $\BG_3$ manifold.

\item \emph{Compute the sector cohomology.}
Restrict the coframe to block coordinates on $C$. Each family satisfies $\sum_a\lambda_a Z_a=0$, and
\[
\dbar Y_a=\sum_b\lambda_bX_b\bar X_b-mX_a\bar X_a, \qquad dY_a=\dbar Y_a,\qquad d\bar Y_a=-dY_a.
\]
The remaining generator differentials vanish. The restricted coframe makes every positive-dimensional $C$ a compact nilmanifold with abelian complex structure, so invariant forms compute its Dolbeault and de Rham cohomology as in Sections~\ref{sec:fourfold} and~\ref{sec:sixfold}. Determine the action of $G_{\lambda,C}$ on these forms and reduce the invariant Dolbeault complex over $\Q$. \cref{lem:sector-comparison} gives the total Betti numbers, and an independent de Rham reduction is a useful numerical check. This gives $h^{a,b}(Z)=\dim H^b(Z,\Omega_Z^{[a]})$ and $b_j(Z)$. A point contributes only in bidegree $(0,0)$ and degree zero.

\item \emph{Assemble the cohomology.} \cref{lem:analytic-decomposition} gives
\[
\begin{aligned}
h^{p,q}(Q)&=\sum_{\lambda,[C]}
 h^{p-r_\lambda,q-r_\lambda}(Z_{\lambda,C}),\\
b_k(Q)&=\sum_{\lambda,[C]}
 b_{k-2r_\lambda}(Z_{\lambda,C}).
\end{aligned}
\]
Omit sectors with $r_\lambda>\min(p,q)$ or $2r_\lambda>k$ in the respective calculations. The equality $b_k(Q)=\sum_{p+q=k}h^{p,q}(Q)$ holds by \cref{thm:all-degeneration}. Use it to check any independently computed Betti numbers.

\item \emph{Form the characteristic-number equations.}
Compute $\chi^p(Q)=\sum_q(-1)^q h^{p,q}(Q)$ and $e(Q)=\sum_k(-1)^k b_k(Q)$. The potentially nonzero Chern numbers are indexed by partitions $\mu$ of $n$: put $C_\mu=\int_Q\prod_i c_{2\mu_i}(T_Q)$. Expand the root-pair formula of \cref{lem:RR} to top degree and equate the coefficients of $y$ with the computed $\chi^p(Q)$. This gives a rational linear system for the $C_\mu$. The Gauss--Bonnet equation $C_{(n)}=e(Q)$ is already its evaluation at $y=-1$, so it gives no additional independent row.

\item \emph{Supply the missing characteristic data.} If the system has full rank, solve for every Chern number and use the Chern--Pontryagin identities above. Otherwise determine enough additional independent characteristic integrals or indices $\chi(Q,\mathcal E)$ to obtain full rank. Possible bundles include tensor products of exterior powers of $\Omega_Q^1$, whose Chern characters are known polynomials in the Chern classes of $T_Q$. Their indices require further geometric information: the additive sector decomposition computes neither arbitrary tensor bundle indices nor the cup product on $Q$. Until this input is supplied, retain the affine solution space of the Riemann--Roch system.
\end{enumerate}
\end{Alg}

\subsection{The first undetermined characteristic number}

The final step is already necessary in complex dimension eight (compare with Rozansky-Witten invariants' story in \cite[Section~5.2]{Sawon2000}). In the ordered Chern-number coordinates $(c_2^4,c_2^2c_4,c_4^2,c_2c_6,c_8)$, the $\chi_y$ system has rank four and kernel generated by $(3,2,1,0,0)$. To see the rank, put $z=x^2$, $A=(1+y)^2$ and $B=y-A/12$. One root-pair factor is
\[
A+Bz+\frac{A}{240}z^2-\frac{A}{6048}z^3 +\frac{A}{172800}z^4+O(z^5).
\]
In a product of four factors, the coefficients of the five symmetric monomial types of total $z$-degree four span $A^4,A^3B,A^2B^2,B^4$. These are independent because $B/A$ is nonconstant. Changing to the Chern-monomial basis preserves the rank and gives the stated kernel. Thus one additional independent characteristic number, for example $c_4^2$, would determine all five numbers once the $\chi_y$-genus of the $\BG_4$ manifold has been computed.

\section*{Acknowledgments}

NK thanks Fedor Bogomolov for his interest in this work and for sharing his ideas. He is grateful to Misha Verbitsky, Alexandra Kuznetsova, and Egor Yasinsky for discussions on the geometry of Bogomolov--Guan manifolds. The Euler characteristic of the fourfold was first obtained in discussions with Kuznetsova and Yasinsky. He also thanks Dima Kaledin for explaining the work of de Cataldo and Migliorini. Both authors thank Justin Sawon for helpful discussions and the London Mathematical Society for supporting the 2024 workshop on Islay, Scotland, where they first discussed Bogomolov--Guan manifolds.

\section*{Disclosure of AI assistance}

AI tools were used to assist with the revision of the exposition, the comparison of the sector complexes, and the preparation and checking of the exact arithmetic calculations.

\bibliographystyle{alpha}
\bibliography{references}

\end{document}